\documentclass[twoside]{article}

\usepackage[accepted]{aistats2024}

\usepackage[round]{natbib}

\usepackage[utf8]{inputenc} 
\usepackage[T1]{fontenc}    
\usepackage{hyperref}       
\usepackage{url}            
\usepackage{booktabs}       
\usepackage{amsfonts}       
\usepackage{nicefrac}       
\usepackage{amsmath}
\usepackage{amssymb}
\usepackage[colorinlistoftodos,bordercolor=orange,backgroundcolor=orange!20,linecolor=orange,textsize=scriptsize]{todonotes}
\usepackage{enumitem}
\usepackage{pifont}
\usepackage{algorithm,algorithmicx,algpseudocode}
\usepackage{wrapfig}
\usepackage[flushleft]{threeparttable}
\usepackage{makecell}
\usepackage{multirow}
\usepackage{caption}
\usepackage{xcolor}
\definecolor{ForestGreen}{RGB}{34,139,34}

\allowdisplaybreaks
\usepackage{mathtools}
\newcommand{\eqdef}{\vcentcolon=}

\usepackage{colortbl}
\definecolor{bgcolor2}{rgb}{0.8,1,1}
\definecolor{bgcolor}{rgb}{0.8,1,0.8}

\def\R{\mathbb{R}}

\def\C{\mathcal C}
\def\Q{\mathcal Q}
\def\X{\mathcal X}

\def\R{\mathbb R}

\def\EE{\mathbb E}

\def\e{\varepsilon}
\def\la{\langle}
\def\ra{\rangle}

\newcommand{\cO}{\mathcal{O}}

\newcommand{\EndProof}[1]{#1 \hfill$\square$}

\def\<#1,#2>{\langle #1,#2\rangle}

\newtheorem{theorem}{Theorem}[section]

\newtheorem{corollary}[theorem]{Corollary}
\newtheorem{lemma}[theorem]{Lemma}
\newtheorem{definition}[theorem]{Definition}
\newtheorem{assumption}[theorem]{Assumption}

\runningauthor{R. Nazykov, A. Shestakov, V. Solodkin, A. Beznosikov, G. Gidel, A. Gasnikov}

\begin{document}

\setlength{\abovedisplayskip}{2pt}
\setlength{\belowdisplayskip}{2pt}

\twocolumn[

\aistatstitle{Stochastic Frank-Wolfe: Unified Analysis and Zoo of Special Cases}

\vspace{-0.3cm}

\aistatsauthor{ Ruslan Nazykov$^*$ \\ MIPT \& Yandex \And 
Aleksandr Shestakov$^*$ \\ MIPT \& Yandex
\And  
Vladimir Solodkin$^*$ \\ MIPT \& Yandex
\AND
Aleksandr Beznosikov \\ MIPT \& MBZUAI \& Yandex
\And 
Gauthier Gidel \\ Mila \& Universit\'e de Montr\'eal
\And 
Alexander Gasnikov \\ University Innopolis \\ ~MIPT \& IITP RAS
}
\aistatsaddress{ } 
]
%

\begin{abstract}

The Conditional Gradient (or Frank-Wolfe) method is one of the most well-known methods for solving constrained optimization problems appearing in various machine learning tasks.  The simplicity of iteration and applicability to many practical problems helped the method to gain popularity in the community. In recent years, the Frank-Wolfe algorithm received many different extensions, including stochastic modifications with variance reduction and coordinate sampling for training of huge models or distributed variants for big data problems. In this paper, we present a unified convergence analysis of the Stochastic Frank-Wolfe method that covers a large number of particular practical cases that may have completely different nature of stochasticity, intuitions and application areas. Our analysis is based on a key parametric assumption on the variance of the stochastic gradients. But unlike most works on unified analysis of other methods, such as SGD, we do not assume an unbiasedness of the real gradient estimation. We conduct analysis for convex and non-convex problems due to the popularity of both cases in machine learning. With this general theoretical framework, we not only cover rates of many known methods, but also develop numerous new methods. This shows the flexibility of our approach in developing new algorithms based on the Conditional Gradient approach. We also demonstrate the properties of the new methods through numerical experiments.
 
\end{abstract}

\section{INTRODUCTION}\label{sec:intro}
\vspace{-0.2cm}

In this paper, we are interested in the constrained optimization problem
\begin{equation}
\label{eq:main}
\min_{\X \subset \R^d} f(x),
\end{equation}
where $\X$ is a convex set. This problem is a cornerstone of applied mathematics, including machine learning. The problem \eqref{eq:main} is at the heart of model training, from classical regressions \citep{shalev2014understanding} to neural networks \citep{goodfellow2016deep}. There are many approaches for solving \eqref{eq:main}. When projection onto a set is expensive (e.g. projection onto the nuclear norm-ball require a full singular value decomposition) or not possible at all (e.g., dual structural SVMs \citep{pmlr-v28-lacoste-julien13}), the Frank-Wolfe method \citep{frank1956algorithm}, also known as Conditional Gradient (see a big survey \citep{braun2022conditional} for more details), is a good option for dealing with \eqref{eq:main}. This approach is based on considering a linear minimization problem on $\X$. The Frank-Wolfe algorithm is one of the classical optimization methods, but it is still relevant even now. Particularly it finds applications in submodular optimization \citep{bachSub2011}, multi-class classification \citep{hazan2016variance}, vision \citep{miech2017learning,bojanowski2014weakly}, group fused lasso \citep{bleakley2011group}, reduced rank nonparametric regression \citep{NIPS2012_f2201f51}, trace-norm based tensor completion (\citep{6138863}), variational inference \citep{krishnan2015barrier} and routing \citep{LEBLANC1975309}, among others. 

Current world reality encourages avoiding the deterministic setting of \eqref{eq:main} and favoring the various stochastic ones. For instance, we often meet problems \eqref{eq:main} with an expectation form of target function: $f(x) = \EE_{\xi \sim \mathcal{D}}[f(x, \xi)]$. Here $\mathcal{D}$ is usually associated with unknown distribution, in terms of machine learning, it corresponds to nature of the data. In such a setting, it is impossible to compute the full gradient, but, despite the fact that data distribution is unknown, we can sample from $\mathcal{D}$ and replace the expectation form with approximation via Monte Carlo: $f(x) = \tfrac{1}{n}\sum_{i=1}^n f_i(x)$ with $n$ samples. 
However, modern application problems become increasingly larger and more computationally complex. Therefore, even for this sum-type problem, computing the full gradient is expensive and should be avoided. We can consider absolutely different randomization techniques in methods to be computationally efficient, e.g., in SGD-type methods stochasticity can be achieved by choosing data batches  \citep{roux2012stochastic, NIPS2014_ede7e2b6, NIPS2013_ac1dd209, pmlr-v70-nguyen17b}, computing gradient coordinates \citep{doi:10.1137/100802001, doi:10.1137/16M1060182, richtarik2013optimal, doi:10.1080/10556788.2016.1190360}, or even by using compression operators \citep{seide2014-bit, NIPS2017_6c340f25, NEURIPS2018_3328bdf9, mishchenko2019distributed}, client samplings \citep{cho2020client, nguyen2020fast, ribero2020communication, chen2020optimal} in distributed and federated settings \citep{konevcny2016federated, DBLP:journals/corr/abs-1912-04977}. Following the trend, we consider a stochastic version of the Conditional Gradient method:
\vspace{0.05cm}
\begin{equation}
\label{eq:stoch_fw}
\begin{split}
s^k = \arg\min_{s \in \mathcal{X}} \< g^k, s - x^k>, \\
x^{k+1} = (1-\eta_k)x^k + \eta_ks^k,
\end{split}
\end{equation}
where $g^k$ is some stochastic estimator of the real gradient $\nabla f(x^k)$.\\
Over its long history, the Frank-Wolfe method received a huge number of different modifications in the form \eqref{eq:stoch_fw}, most of them in the last decade. Here one can note modifications related to variance reduction \citep{reddi2016stochastic, hazan2016variance, pmlr-v80-qu18a, pmlr-v97-yurtsever19b, pmlr-v119-gao20b, negiar2020stochastic, lu2021generalized, weber2022projection, Beznosikov2023SarahFM}, coordinate randomization \citep{pmlr-v28-lacoste-julien13, pmlr-v89-sahu19a}, and distributed computation \citep{bellet2015distributed, pmlr-v48-wangd16, hou2022distributed}. However, all of these separate practical variants of the Conditional Gradient method have different intuitions of convergence, various formal techniques of proving and do not always cover the same cases of assumptions on the target function $f$. Moreover, there remains a rather large gap in what can still be done in the creation of Frank-Wolfe's modfications. Namely, there are advanced SGD-type methods that have not yet been adapted for use within the Frank-Wolfe iteration. These include some new coordinate approaches and approaches to finite-sum problems. Meanwhile, these two techiques help to address the main bottleneck of distributed algorithms -- expensive communications. All of mentioned issues lead us to the two key questions of this paper:
\vspace{-0.4cm}
\begin{quote} 
    \textit{
    1. Can we conduct a novel general analysis of the Stochastic Frank-Wolfe unifying special cases and providing the ability to design new extensions?
    \vspace{0.05cm}
    \\
    2. What new stochastic modifications of the classical Conditional Gradient can we possibly invent based on this unified analysis?
    }
\end{quote}
\vspace{-0.2cm}

\vspace{-0.2cm}
\subsection{Our contribution}
\vspace{-0.2cm}

$\bullet$ \textbf{Unified analysis of Stochastic Frank-Wolfe.} 
We propose a general assumption on the stochastic estimator $g^k$ from Stochastic Frank-Wolfe \eqref{eq:stoch_fw} -- see Assumption~\ref{as:key}. Below we note in more details that our assumption is broad and encompasses many special cases, in particular, those that could not be analyzed in a unified way before. Under Assumption~\ref{as:key} we present general convergence results for the problem \eqref{eq:main}.

$\bullet$ \textbf{Convex and non-convex cases.} Motivated by various applications primarily from machine learning, we provide the unified analysis in the convex (Theorem \ref{th:convex} and non-convex (Theorem \ref{th:nonconvex}) cases of the target function $f$. This is also interesting for special cases, since the authors of some papers do not give an analysis in both setups. 

$\bullet$ \textbf{Without assumptions of unbiasedness.} In our key Assumption~\ref{as:key}, we bound the variance $\EE[\| g^k - \nabla f(x^k)\|^2]$ using universal letter constants and an additional auxiliary sequence. Similar assumptions are made in papers on analysis of the SGD family methods \citep{pmlr-v108-gorbunov20a, li2020unified, DBLP:journals/corr/abs-2006-11573}. But these works also additionally assume that the stochastic gradient $g^k$ is unbiased, i.e., $\EE[g^k \mid x^k] = \nabla f(x^k)$. We avoid this assumption, it extends the class of methods that can be considered under our assumptions compared to the works around the SGD-type methods. For example, it allows to prove the convergence of distributed methods with biased/greedy compression \citep{DBLP:journals/corr/abs-1909-05350, NEURIPS2021_231141b3} or SARAH-based variance reduced methods \citep{nguyen2017sarah}.

$\bullet$ \textbf{Vast number of new special cases.}
The previous point already gives an indication of the breadth and flexibility of the approach. Our general theoretical framework allows us to analyze different variants of the classical Frank-Wolfe method. Guided by algorithmic advances for solving unconstrained minimization problems we present
a new method with coordinate randomization (\texttt{SEGA FW}), a new variance-reduced method (\texttt{L-SVRG FW}, \texttt{SARAH FW}, \texttt{SAGA FW}), new distributed methods with unbiased compression (\texttt{DIANA FW}, \texttt{MARINA FW}) and biased compression (\texttt{EF21 FW}), and others. Although the SGD-type analogs of these methods are known for solving primarily unconstrained minimization problems \citep{NEURIPS2018_fc2c7c47, pmlr-v117-kovalev20a, pmlr-v139-gorbunov21a, NEURIPS2021_231141b3}, they were never integrated into the Frank-Wolfe iteration for solving projection-free free constrained problems. We also demonstrate that our general theorems allow to obtain convergence for methods that are combinations of the two basic approaches, i.e, \texttt{SAGA SARAH FW} and others. Moreover, we presents absolutely new methods that are not found in the literature on SGD. This algorithm uses special coordinate randomization. \texttt{JAGUAR} is a new coordinate method. \texttt{Q-L-SVRG FW} is a new distributed method with unbiased compression. This method is based on non-distributed \texttt{L-SVRG FW}: instead of the randomness of choosing a batch/term number, randomness from compression is used. \texttt{PP-L-SVRG FW} is a new distributed method with client sampling also based on non-distributed \texttt{L-SVRG FW}.

$\bullet$ \textbf{Sharp rates for known special cases.}
For the known methods fitting our framework our general theorems either recover the best rates known for these methods. These methods include \texttt{SARAH FW}, \texttt{SAGA SARAH FW} \citep{Beznosikov2023SarahFM}.

$\bullet$ \textbf{Numerical experiments.} In numerical experiments, we illustrate the most important properties of the new methods. The results corroborate our theoretical findings.

Throughout the paper, we provide necessary comparisons with closely related work.

\vspace{-0.2cm}
\subsection{Technical preliminaries}
\vspace{-0.2cm}

\textbf{Notations.} We use $\la x,y \ra \eqdef \sum_{i=1}^d x_i y_i$ to denote standard inner product of $x,y\in\R^d$, where $x_i$ corresponds to the $i$-th component of $x$ in the standard basis in $\R^d$. It induces $\ell_2$-norm in $\R^n$ in the following way: $\|x\| := \sqrt{\la x, x \ra}$. Operator $\EE[\cdot]$ denotes full mathematical expectation and operator $\EE[\cdot | x^k]$ express conditional mathematical expectation w.r.t. all randomness coming from the $k$th iteration of \eqref{eq:stoch_fw}. We introduce $f^*$ as a solution of the problem \eqref{eq:main}, i.e. a global minimum of $f$ on the convex set $\X$. For the non-convex function $f$ the solution $f^*$ may not be unique. We also define $\Delta_0 := f(x^0) - f^*$, where $x^0$ is a starting point of \eqref{eq:stoch_fw}.

Throughout the paper, we assume that the target $f$ from \eqref{eq:main} satisfies the following assumptions.
\begin{assumption} \label{as:lip}
The function $f: \X \to \R$, is $L$-smooth on $\X$, i.e., there exists a constant $L > 0$ such that 
$\| \nabla f(x) - \nabla f(y)\| \leq L \|x-y \|$ ~for all $x,y \in \X$.
\end{assumption}

\begin{assumption} \label{as:conv}
The function $f: \X \to \R$, is convex, i.e., $f(x) \geq f(y) + \< \nabla f(y) , x - y>$ ~for all $x,y \in \X$.
\end{assumption}
While we always need the assumption on smoothness of $f$, we abandon the convexity in one of the main theorems. The next assumption is also key for the design and analysis of Frank-Wolfe-type methods.
\begin{assumption} \label{as:set}
The set $\X$ is convex and compact with a diameter $D$, i.e., for any $x,y \in \X$,
\begin{equation*}
\| x - y \| \leq D.
\end{equation*}
\end{assumption}
For some particular cases of the method \eqref{eq:stoch_fw}, we need to introduce additional objects and assumptions on them. This will be done in the corresponding sections.

\vspace{-0.2cm}
\section{MAIN THEOREMS}
\vspace{-0.2cm}
In this section, we first present the central part of our approach that allows us to conduct a general analysis of the algorithms (Assumption \ref{as:key}), then we provide convergence analysis for both convex and non-convex cases.
\vspace{-0.2cm}
\subsection{Unified assumption}
\vspace{-0.2cm}
First, we introduce the central part of our approach, all subsequent analysis is based on the following assumption on the stochastic gradients $g^k$ :
\begin{assumption} \label{as:key}
Let $\{x^k\}_{k=0}^K$ be the iterates produced by Stochastic Frank-Wolfe (see \eqref{eq:stoch_fw}).
Let there exist constants $A, B, C, E \geq 0$, $\rho_1, \rho_2 \in (0; 1]$ and a
(possibly) random sequence $\{\sigma_k\}_{k\geq0}$ such that the following inequalities hold
{\small
\begin{equation}
\label{grad_est}
\begin{split}
        \hspace{-0.1cm}\EE[\|g^k - \nabla f(x^k)\|^2\,|\,x^k] \leq& (1 - \rho_1) \|g^{k-1} - \nabla f(x^{k-1})\|^2 \\
        &{ +  A \sigma^2_{k-1} +\eta^{2}_{k-1} B D^2 + C,}
\end{split}
\end{equation}
\begin{equation}\label{rand_seq}
\EE[\sigma_k^2~|~x^k] \leq (1 - \rho_2) \sigma^2_{k-1} +  \eta^{2}_{k - 1} E D^2.
\end{equation}
}
\end{assumption}
The inequality $\ref{grad_est}$ bounds the second moment of stochastic estimation $g^k$. The sequence $\lbrace\sigma_k^2\rbrace_{k\geq0}$ is needed to capture the variance, which can be reduced during the algorithm's work process. Constants $\rho_1,\rho_2$ show how quick this reduction is regarding the previous iteration. $B,E$ provide the information of convergence depending on the previous step size and set's geometry (Assumption \ref{as:set}). Finally, constant $C$ stands for the remaining noise that cannot be reduced as $\sigma_k$.

Proposed assumption takes into account specificity of the Frank-Wolfe analysis in terms of upper bound containing component $D^2$ as irremovable part of such type inequalities. 

\vspace{-0.2cm}
\subsection{Convergence results}
\vspace{-0.2cm}
\textbf{Convex case}. The following theorem describes the convergence rate of stochastic Frank-Wolfe (\ref{eq:stoch_fw}) based methods under the convexity of $f$:
\begin{theorem}
\label{th:convex}
Let Assumptions~\ref{as:lip},~\ref{as:conv},~\ref{as:set} and~\ref{as:key}  be satisfied. For any $K$ choose step sizes $\{\eta_k\}_{k\geq1}$ as follows:
\begin{align*}
&\text{if ~K } \leq \text{ d}, && \eta_k = \tfrac{1}{d},
\\
&\text{if ~K } > \text{ d ~and ~k } < \text{k}_0, && \eta_k = \tfrac{1}{d},
\\
&\text{if ~K } > \text{ d ~and ~k } \geq \text{k}_0, && \eta_k = \tfrac{2}{2d + k - k_0},
\end{align*}
where $d = \tfrac{2}{\min(\rho_1, \rho_2)}$, $k_0 = \left\lceil\tfrac{K}{2}\right\rceil$.
Then the output of Stochastic Frank-Wolfe after $K$ iterations satisfies
\begin{equation*}
    \begin{split}
        \EE [r_{K+1}] = \mathcal{O}\Big(&\textstyle{ r_0 \exp\left(-\frac{K}{2d}\right)
        + \frac{LD^2}{K+d}}\\ 
        &\textstyle{+ \frac{D^2}{K+d}\sqrt{\frac{B\rho_2 + AE}{\rho_1\rho_2}} + \sqrt{\frac{KD^2}{K +d}\frac{C\rho_2}{\rho_1\rho_2}}\Big),}
    \end{split}
\end{equation*}
where  the Lyapunov function $r_k$ is defined by
$
    r_k := f(x^k) - f^* + M_1\|\nabla f(x^k) - g^k\|^2 + M_2\sigma_{k}^{2}
$
with $M_1,M_2 > 0$.
\end{theorem}
The proof is provided in Section \ref{sec:uni}. Note, that with zero noises ($C=0$ in Assumption \ref{as:key}) this theorem reflects sublinear convergence $\mathcal{O}\left(\tfrac{1}{K}\right)$.

\textbf{Non-convex case.} 
To obtain our convergence results for the non-convex objective function we introduce the \textit{Frank-Wolfe gap} function \citep{pmlr-v28-jaggi13} as a convergence criterion:
\begin{equation*}
\begin{split}
    \textbf{\textit{gap}}(y) = \max\limits_{x\in \mathcal{X}}\langle\nabla f(y),y-x\rangle
\end{split}
\end{equation*}
Such type of criterion is standard for analyzing the convergence of constrained optimization algorithms in the non-convex case \citep{reddi2016stochastic}. It is shown in \citep{lacostejulien2016convergence} that FW gap is an affine invariant generalization of standard convergence criterion $\|\nabla f(y)\|$ and therefore a meaningful measure of non-stationarity. In terms of FW gap we derive the following general convergence result:
\begin{theorem}
\label{th:nonconvex}
Let the Assumptions ~\ref{as:lip},~\ref{as:set} and~\ref{as:key} be satisfied. Then, there exist constants $M_1$, $M_2$ such that for any $K$ there exist constant $\lbrace\eta_k\rbrace_{k\geq 1} \equiv \tfrac{1}{\sqrt{K}}$ for \eqref{eq:stoch_fw}, thus
\begin{equation*}
    \begin{split}
        \EE\Bigl[\min\limits_{0\leq k \leq K-1} \textbf{gap}(x^k)\Bigr] = \mathcal{O}\Big(&\textstyle{\frac{r_0}{\sqrt{K}} + \frac{D^2}{\sqrt{K}}\Big[L + \sqrt{\frac{B\rho_2 + AE}{\rho_1\rho_2}}\Big]}\\ 
        & \textstyle{+ \sqrt{D^2\frac{C\rho_2}{\rho_1\rho_2}}\Big),}
    \end{split}
\end{equation*}
where $r_0 := f(x^0) - f^* + M_1 \|g^0 - \nabla f(x^0) \|^2 + M_2\sigma_0^2$ with $M_1,M_2 > 0$.
\end{theorem}
See the proof in Section \ref{sec:uni}. Note, that in the case of zero noises ($C=0$ in Assumption \ref{as:key}) this result establishes as $\mathcal{O}\left(\tfrac{1}{\sqrt{K}}\right).$

\vspace{-0.2cm}
\section{WIDE VARIETY OF SPECIAL METHODS}\label{sec:zoo}
\vspace{-0.2cm}

In this section, we fulfill the promises made in the introduction and show how many existing and new techniques fit our framework. Due to space restrictions, the comparison and the full listing of algorithms are described particularly in Section \ref{sec:comparison}.
\vspace{-0.3cm}
\subsection{Stochastic methods} \label{sec:stoch_main}
\vspace{-0.2cm}
As already mentioned in Section \ref{sec:intro}, in modern applications we often deal with finite-sum optimization problems (so-called empirical risk minimization):
\begin{equation}
    \label{eq:finsum}
    \textstyle{\min_{x\in\mathcal{X}}~f(x) := \frac{1}{n}\sum_{i=1}^n f_i(x).}  
\end{equation}
An important detail of this setting is that calling the full gradient of $f$ is expensive, only small batches $\tfrac{1}{b} \sum_{i=1}^b \nabla f_i (x)$ can be typically used. Therefore, for the theoretical analysis, we  need not only the smoothness of the function $f$, but also of all summands $f_i$.
\begin{assumption} \label{as:lip:local}
Each function $f_i: \X \to \R$ is $L_i$-smooth on $\X$, i.e., there exist constants $\{L_i\}>0$ such that
$\|\nabla f_i(x) - \nabla f_i(y)||\leq L_i||x-y||$ ~for all $x,y \in \X.$
We also define $\widetilde{L}$ as $\widetilde{L}^2 = \frac{1}{n}\sum_{i=1}^n L_i^2$.
\end{assumption}

\textbf{L-SVRG FW.}
One of the most popular stochastic algorithms for \eqref{eq:finsum} with $\X \equiv \R^d$ is SVRG \citep{NIPS2013_ac1dd209}. We consider its loopless variant \citep{pmlr-v117-kovalev20a} called L-SVRG that uses SVRG idea but is a bit more friendly for theoretical analysis. In more details, we need to compute $g^k$ as follows:
\begin{align}
\label{eq:lsvrg}
        w^{k+1} =& 
            \begin{cases}
                x^k, & \text{ with probability }~~~ p,\\
                w^k, & \text{ with probability } 1-p,
            \end{cases}
        \\
        \notag
        g^{k+1} =& \textstyle{\tfrac{1}{b}\sum_{i\in S_k} [\nabla f_i(x^{k+1}) - \nabla f_i(w^{k+1})] + \nabla f(w^{k+1})},
\end{align}
where batches of indexes $S_k$ size of $b$ are generated uniformly and independently.
The essence of this approach is that the probability $p$ is taken close to zero, then the full gradient at the point $w^{k+1}$ are computed quite rarely and in most cases we use the approximation $g^{k+1}$ via stochastic gradients on mini-batches of random indexes $S_k$ size of $b$. 
\begin{lemma}
    \label{lem:lsvrg_conv}
    Under assumptions \ref{as:lip}, \ref{as:conv}, \ref{as:lip:local} the algorithm \eqref{eq:stoch_fw} + \eqref{eq:lsvrg} satisfies assumption ~\ref{as:key} with: $\rho_1 = 1,$ $A = \frac{\widetilde{L}^2}{b}\left(1-\frac{p}{2}\right)$, $B = \frac{8\widetilde{L}^2}{pb}$, $C = 0$, $\sigma_k^2 = \|x^k-w^k\|^2$, $\rho_2 = \frac{p}{2}$, $E =\frac{8}{p}$.
\end{lemma}

Using this lemma, one can get the convergence of \eqref{eq:stoch_fw} + \eqref{eq:lsvrg} in both convex and non-convex cases. 
\begin{corollary}\label{crl:lsvrg}
    For the algorithm \eqref{eq:stoch_fw}+\eqref{eq:lsvrg} in the convex and non-convex cases accordingly the following convergences take place:
    \begin{eqnarray*}
        &\hspace{-0.3cm}
        \textstyle{
        \mathbb{E}[\text{{\small $f(x^{K}) - f^*$}}] =\mathcal{O}\left(\Delta_0\exp\left(-\frac{Kp}{8}\right) + \frac{ L D^2}{K}\left[1 + \frac{\widetilde{L}}{L}\frac{1}{p\sqrt{b}}\right]\right),
        }
        \\
        &\hspace{-0.3cm}
        \textstyle{
        \EE\Bigl[\min\limits_{0 \leq k \leq K-1} \textbf{gap}(x^k)\Bigr] = \mathcal{O}\left(\frac{\Delta_0}{\sqrt{K}} + \frac{L D^2}{\sqrt{K}}\left[1 + \frac{\widetilde{L}}{L}\frac{1}{p\sqrt{b}}\right]\right).
        }
    \end{eqnarray*}
\end{corollary}
See more details in Section \ref{sec:lsvrg}.

\textbf{SARAH FW.}
Another common algorithm for solving the unconstrained version of \eqref{eq:finsum} is SARAH \citep{nguyen2017sarah}. In particular, it has better theoretical results in both convex \citep{nguyen2017sarah}, non-convex \citep{pmlr-v139-li21a} target functions and also bits SVRG on practice. As in the previous method we look at the loopless version of SARAH \citep{pmlr-v139-li21a}:
    {\small 
\begin{equation}
\label{eq:sarah}
\hspace{-0.33cm}
    g^{k+1} =
    \begin{cases}
        \nabla f(x^{k+1}),~\hspace{1.7cm}\text{ with probability }p, \\
        \textstyle{g^{k} + \frac{1}{b}{\sum_{i\in S_k}}\left[\nabla f_i(x^{k+1}) - \nabla f_i(x^{k})\right]},~\text{oth.},
    \end{cases}
\end{equation}
}

where batches of indexes $S_k$ size of $b$ are generated uniformly and independently.
The methods main idea here is very close to SVRG since it also computes the full gradient only with small probability $p$. However, the approximation of the gradient $g^{k+1}$ in SARAH is done not by the old point $w^{k+1}$ as in SVRG, but more smoothly using the current and previous points: $x^{k+1}$ and $x^{k}$.
\begin{lemma}
\label{lem:sarah}
Under Assumptions \eqref{as:lip}, \eqref{as:conv}, \eqref{as:lip:local} the algorithm \eqref{eq:stoch_fw}+\eqref{eq:sarah} satisfies Assumption \ref{as:key} with:
$\rho_1 = p$, $A = 0$, $B = \frac{1-p}{b}\widetilde{L}^2$, $C = 0$, $\sigma_k = 0$, $\rho_2 = 1$, $E = 0$.
\end{lemma}

\begin{corollary}
\label{crl:sarah}
For the algorithm \eqref{eq:stoch_fw}+\eqref{eq:sarah} in the convex and non-convex cases respectively the following convergences take place:
\begin{eqnarray*}
    &\hspace{-0.3cm}\mathbb{E}[\text{{\small $f(x^{K}) - f^*$}}] = \textstyle{\mathcal{O}\left(\Delta_0 \exp\left(-\frac{Kp}{4}\right) + \frac{LD^2}{K}\left[1 + \frac{\widetilde{L}}{L} \frac{1}{\sqrt{pb}}\right]\right),}
    \\
    &\hspace{-0.3cm}\mathbb{E}\Bigl[\min\limits_{0 \leq k \leq K-1} \textbf{gap}(x^k)\Bigr] = \textstyle{\mathcal{O}\Bigl(\frac{\Delta_0}{\sqrt{K}} + \frac{LD^2}{\sqrt{K}}\left[1 + \frac{\widetilde{L}}{L}\frac{1}{\sqrt{pb}}\right]\Bigr).}
\end{eqnarray*}
\end{corollary}
The details are provided in Section \ref{sec:SARAH}. It is easy to see that with the same $p$ the results for \texttt{SARAH FW} is better than for \texttt{L-SVRG FW} (Corollary \ref{crl:lsvrg}).

\textbf{SAGA FW.}
The final algorithm of provided lineup for \eqref{eq:finsum} is SAGA \citep{NIPS2014_ede7e2b6}:
\begin{align}
\label{eq:saga}
        y_i^{k+1} &=
            \begin{cases}
            \nabla f_i(x^k), & \text{ for } i\in S_k, \\
            y_i^k, & \text{ for }i \notin S_k,
            \end{cases}
            \\
            \notag
            g^{k+1}  =&\textstyle{ \frac{1}{b}{\sum_{i\in S_k}}[\nabla f_i(x^{k+1}) -  y_i^{k+1}] + \frac{1}{n}\sum_{j=1}^n y_j^{k+1}},
\end{align}
where batches of indexes $S_k$ size of $b$  are generated uniformly and independently.
The essence of the SAGA technique is different from SVRG and SARAH. In this case, we do not compute full gradients even rarely, but we need to store an additional set of vectors $\{y_i^k\}_{i=1}^n$. The vector $y_i$ stores information about the last gradient of the function $f_i$ that was computed during the operation of the algorithm. Thus one can state that we collect "delayed" full gradient in $\frac{1}{n}\sum_{j=1}^n y_k^j$. 
\begin{lemma}
\label{lem:saga}
Under Assumptions \ref{as:lip}, \ref{as:conv}, \ref{as:lip:local} the algorithm \eqref{eq:stoch_fw}+\eqref{eq:saga} satisfies Assumption \ref{as:key} with:
$\rho_1 = 1$, $A = \frac{1}{b}\left(1 + \frac{b}{2n}\right)$, $B = \frac{2\widetilde{L}^2}{b}\left(1 + \frac{2n}{b}\right)$, $C = 0$, $\sigma_k^2 = \frac{1}{n}\sum_{j = 1}^n\|\nabla f_j(x^k) - y_j^{k+1}\|^2$, $\rho_2 = \frac{b}{2n}$, $E = \frac{2n}{b}\widetilde{L}^2$.
\end{lemma}

\begin{corollary}
\label{crl:saga}
For the algorithm \eqref{eq:stoch_fw}+\eqref{eq:saga} in the convex and non-convex cases the following convergences take place:
\begin{eqnarray*}
    &\hspace{-0.3cm}\mathbb{E}[\text{{\small $f(x^{K}) - f^*$}}] = \mathcal{O}\left(\Delta_0\exp\left(-\frac{Kb}{8n}\right) + \frac{LD^2}{K}\left[1 + \frac{\widetilde{L}}{L}\frac{n}{b\sqrt{b}}\right]\right),
    \\
    &\hspace{-0.3cm}\mathbb{E}\Bigl[\min\limits_{0 \leq k \leq K-1} \textbf{gap}(x^k)\Bigr] = \mathcal{O}\left(\frac{\Delta_0}{\sqrt{K}} + \frac{LD^2}{\sqrt{K}}\left[1 + \frac{\widetilde{L}}{L}\frac{n}{b\sqrt{b}}\right]\right).
\end{eqnarray*}
\end{corollary}
The full statement together with its proof can be found in Section \ref{sec:SAGA}. 
One more method (\texttt{SAGA SARAH FW}) for the stochastic setting will be analyzed in Section \ref{sec:combination_main}. A comparison of the presented and already existing methods is provided in Section \ref{sec:comparison_stoch}. In particular, there we analyze Corollaries \ref{crl:lsvrg}--\ref{crl:saga} with the substituted optimal parameter value of $p$.
\vspace{-0.2cm}
\subsection{Coordinate methods}\label{sec:coord_main}
\vspace{-0.2cm}
Previous approaches reduce the cost of gradient computing by selecting small batches, but there are other strategies. In particular, wide range of algorithms use random sampling of coordinates for gradient evaluation \citep{doi:10.1137/100802001, richtarik2013optimal, qu2016coordinate}. This technique can also significantly decrease the computational cost. Then, in this section, we focus on methods, where gradient estimator stochastically depend on function's partial derivatives.

\textbf{SEGA FW.} The original algorithm developed in \citep{NEURIPS2018_fc2c7c47} covers a general setting, instead of which we use a slightly more simplified version. Particularly, we update $g^k$ as follows:
\begin{equation}
\label{eq:sega}
    \begin{split}
        h^{k+1} &= h^k + e_{i_k}(\nabla_{i_k} f(x^k) - h_{i_k}^k), \\
        g^{k+1} &= d(\nabla_{i_k} f(x^{k+1}) - h_{i_k}^{k+1})e_{i_k} + h^{k+1},
    \end{split}
\end{equation}
where coordinate $i_k$ is chosen uniformly and independently. 
The idea of this approach is in some sense close to SAGA. We also have some memory buffer, but unlike SAGA, where we save the last calculated gradient on $i$th batch, here in $h_i$ we save the last partial derivative $\nabla_i f$ calculated for the $i$th coordinate. 
\begin{lemma} 
\label{lem:sega}
Under Assumptions \ref{as:lip}, \ref{as:conv} the algorithm \eqref{eq:stoch_fw} +\eqref{eq:sega} satisfies Assumption \ref{as:key} with:
$\rho_1 = 1$, $A = d$, $B = d^2L^2$, $C = 0$, $\sigma_k^2 =\|h^{k+1} - \nabla f(x^k)\|^2$, $\rho_2 = \frac{1}{2d}$, $E = 3L^2d$.
\end{lemma}

\begin{corollary}
\label{crl:sega}
For the algorithm (\ref{eq:stoch_fw})+(\ref{eq:sega}) in the convex and non-convex cases the following convergences take place:
\begin{eqnarray*}
    &\mathbb{E}[\text{{\small $f(x^{K}) - f^*$}}] =\mathcal{O}\Big(\Delta_0 \exp\Big(-\frac{K}{8d}\Big) + \frac{LD^2}{K} \cdot d\sqrt{d} \Big),
    \\
    &\EE\Big[\min\limits_{0 \leq k \leq K-1} \textbf{gap}(x^k)\Big] = \mathcal{O}\Big(\frac{\Delta_0}{\sqrt{K}} + \frac{LD^2}{\sqrt{K}}\cdot d\sqrt{d} \Big).
\end{eqnarray*}
\end{corollary}
See details in Section \ref{sec:SEGA}. Despite the fact that SEGA has recommended itself as an effective coordinate method, theoretical estimation of demanded steps to converge has undesired component $d\sqrt{d}$ (while the original algorithm has only $d$). Therefore, it makes sense to introduce a new algorithm, which has theoretically better convergence.

\textbf{JAGUAR.} An important feature of SEGA is the fact that it uses unbiased gradient estimation: $\EE[g^k \mid x^k] = \nabla f(x^k)$. On the one hand it is good and helps to simplify the theoretical analysis. But experience shows that stochastic methods with biased gradient approximation can outperform unbiased ones. The example of SARAH (biased) and SVRG (unbiased) supports this. Thus, we propose to consider the following form of $g^k$:
\begin{equation}
    \label{eq:jaguar}
    g^{k+1} = e_{i_k}(\nabla_{i_k} f(x^{k+1}) - g_{i_k}^{k}) + g^{k},
\end{equation}
where coordinate $i_k$ is chosen uniformly and independently. 
\begin{lemma}
\label{lem:yaguar}
Under Assumptions \ref{as:lip}, \ref{as:conv} the algorithm \eqref{eq:stoch_fw}+\eqref{eq:jaguar} satisfies Assumption \ref{as:key} with:
$\rho_1 = 1, A = 0, ~ B = 3dL^2,~C = 0,~$
$\sigma_k^2 = 0 ,~\rho_2 = 1, E =0.$
\end{lemma}

\begin{corollary}
\label{crl:yaguar}
For the algorithm \eqref{eq:stoch_fw}+\eqref{eq:jaguar} in the convex and non-convex cases the following convergences take place:
\begin{eqnarray*}
    &\mathbb{E}[\text{{\small $f(x^{K}) - f^*$}}] =
    \mathcal{O}\left(\Delta_0 \exp\left(-\frac{K}{8d}\right) + \frac{LD^2}{K} \cdot d\right),
    \\
    &\EE\left[\min\limits_{0 \leq k \leq K-1} \textbf{gap}(x^k)\right] = \mathcal{O}\left(\frac{\Delta_0}{\sqrt{K}} + \frac{LD^2}{\sqrt{K}} \cdot d \right).
\end{eqnarray*}
\end{corollary}
One can found the full statement together with its proof in Section \ref{sec:JAGUAR}. A comparison of the presented and already existing methods is presented in Section \ref{sec:comparison_coord}.
\vspace{-0.2cm}
\subsection{Distributed methods with compression} \label{sec:distr_main}
\vspace{-0.2cm}
In this section, we focus on distributed versions of Frank-Wolfe algorithm for solving finite-sum problems \eqref{eq:finsum}, where $\{f_i\}_{i = 1}^{n}$ are distributed across $n$ devices connected with parameter-server in a centralized way and each device has an access only to $f_i$. Here we allow different machines to have different data distributions, i.e., we consider the heterogeneous data setting. For such type of problem, the bottleneck commonly is a communication cost \citep{konevcny2016federated}, which motivates to use compressed communication \citep{seide20141}. 
To formally describe compression we introduce the following definition.
\begin{definition}
\label{def:unbcomp}
\textit{Map } $\Q:\mathbb{R}^d\rightarrow\mathbb{R}^d$ \text{ is an } unbiased compression operator, \text{ if there exist a constant } $\omega\geq0$ \text{such that for all } $x\in\mathbb{R}^d$ $$\mathbb{E}[\Q(x)] = x,\;\;\; \mathbb{E}[\|\Q(x) - x\|^2] \leq \omega\|x\|^2.$$ 
\end{definition}
Examples of such operators are random coordinate selection or randomized roundings \citep{beznosikov2022biased}. 
The usage of unbiased compression has been extensively studied. The first and simplest idea that comes to a mind is to apply compression directly on the gradient estimator when forwarding to the server \citep{NIPS2017_6c340f25}. But this kind of approach has a problem, namely for a fixed step it guarantees convergence only to the neighborhood of the solution. Therefore, we propose to consider more advanced techniques that compress some difference that tends to zero during the course of the algorithm. 

\textbf{DIANA FW.} We start with the DIANA technique \citep{mishchenko2019distributed}. Its essence lies in the fact that it maintains the "memory" variables $h_i^k$ at each worker $i$ and compresses gradient differences $\nabla f_i(x^k) - h_i^k$. In particular,
\begin{equation}
\label{eq:diana}
\begin{split}
    &\Delta_i^k = \mathcal{Q}( \nabla f_i^k - h_i^k),\;\; h_i^{k+1} = h_i^k + \alpha\cdot \Hat{\Delta}_i^k,
    \\
    &\textstyle{h^{k+1} = h^k + \alpha \cdot \frac{1}{n}\sum_{i=1}^n \Delta_i^k,} 
    \\
    &\textstyle{g^{k+1} = h^{k+1} + \frac{1}{n}\sum_{i=1}^n \Delta_i^{k+1} ,}
\end{split}
\end{equation}
where first two equations belong to the local computations and last two -- to parameter-server computation. In is important to highlight that we need only compressed differences ${\Delta}_i^k$ for the server updates. 
\vspace{-0.15cm}
\begin{lemma}
    Under Assumptions \ref{as:lip}, \ref{as:conv}, \ref{as:lip:local} the algorithm \eqref{eq:stoch_fw} +\eqref{eq:diana} satisfies Assumption \ref{as:key} with:
$\rho_1 = 1$, $A = \frac{\omega}{n^2}$, $B = \frac{2\omega(\omega+1)\widetilde{L}^2}{n}$, $C = 0$, $\sigma_k^2 = \sum_{i=1}^n\|\nabla f_i(x^k) - h_i^k\|^2$, $\rho_2 = \frac{1}{2(1+\omega)}$, $E = 2(\omega+1)n\widetilde{L}^2$.
\end{lemma}
\begin{corollary}
\label{crl:diana}
For the algorithm \eqref{eq:stoch_fw}+\eqref{eq:diana} in the convex and non-convex cases the following convergences take place:
\begin{eqnarray*}
    &\hspace{-0.3cm}\mathbb{E}[\text{{\small $f(x^{K}) - f^*$}}] =
    \mathcal{O}\Big(\Delta_0\exp\left(-\frac{K}{8\omega}\right)+\frac{LD^2}{K} \left(1 + \frac{\widetilde{L}}{L}\frac{\omega^{\frac{3}{2}}}{\sqrt{n}}\right) \Big),
    \\
    &\hspace{-0.3cm}\EE\left[\min\limits_{0 \leq k \leq K-1} \textbf{gap}(x^k)\right] = \mathcal{O}\left(\frac{\Delta_0}{\sqrt{K}} + \frac{LD^2}{\sqrt{K}} \left(1 + \frac{\widetilde{L}}{L}\frac{\omega^{\frac{3}{2}}}{\sqrt{n}}\right) \right).
\end{eqnarray*}
\end{corollary}
The full proof and convergence results are presented in Section \ref{sec:DIANA}. 

\textbf{MARINA FW.}
Next, as a natural generalization of the idea of compressing gradient differences, we arrive at the fact that gradient estimator could be biased. We have already seen this idea take place with the examples of \texttt{L-SVRG FW} vs \texttt{SARAH FW} (Section \ref{sec:stoch_main}), as well as 
\texttt{SEGA FW} vs \texttt{JAGUAR} (Section \ref{sec:coord_main}). Therefore, we consider the work by \citep{pmlr-v139-gorbunov21a}, where the authors bases their method on SARAH technique but with compression stochasticity instead of \eqref{eq:sarah}. Our algorithm utilizes this idea and performs the update rule \eqref{eq:stoch_fw} with
\begin{align}
     \label{eq:marina}
     \notag
         c_i^{k+1} =&
    \begin{cases}
        \nabla f_i(x^{k+1}) - g_i^{k}, ~\hspace{1.1cm} \text{with probability }p,\\
        \text{{\small $\mathcal{Q}(\nabla f_i(x^{k}) - \nabla f_i(x^{k-1}))$}}
        \text{, otherwise},
    \end{cases}
    \\
        g_i^{k+1} = & g_i^{k} + c_i^{k+1}
    \\
    \notag
    g^{k+1} =& g^k + \textstyle{\frac{1}{n}\sum_{i=1}^n c_i^{k+1},
    }
\end{align}
where $g_i^{k+1}$ are computed on the local devices and $g^{k+1}$ -- on server side. In fact, with $p$ close to zero, one can note that to compute $g^{k+1}$ we typically need only compressed differences: $g^{k+1} = g^{k} + \frac{1}{n}\sum_{i=1}^n \mathcal{Q}(\nabla f_i(x^{k+1}) - \nabla f_i(x^{k}))$. But rarely (with probability $p$) we need to send 
the full uncompressed gradients.
\begin{lemma}\label{lem:marina}
    Under Assumption \ref{as:lip:local} the algorithm \eqref{eq:stoch_fw} +\eqref{eq:marina} satisfies Assumption \ref{as:key} with:
$\rho_1 = p$, $\rho_2 = 1$, $A = 0$, $B = \frac{(1-p)\omega L^2}{n}$, $C = 0$, $\sigma_k = 0$, $E = 0$.
\end{lemma}
\begin{corollary}\label{crl:marina}
    For the algorithm \eqref{eq:stoch_fw}+\eqref{eq:marina} in the convex and non-convex cases respectively the following convergences take place:
\begin{eqnarray*}
    &\hspace{-0.3cm}\mathbb{E}[\text{{\small $f(x^{K}) - f^*$}}] = \textstyle{\mathcal{O}\Big(\Delta_0\exp\left(-\frac{pK}{4}\right) + \frac{LD^2}{K}\Big(1 + \sqrt{\frac{\omega}{pn}} ~\Big)\Big),}
    \\
    &\hspace{-0.3cm}\EE\Bigl[\min\limits_{0 \leq k \leq K-1} \textbf{gap}(x^k)\Bigr] = \textstyle{\mathcal{O}\Big(\frac{\Delta_0}{\sqrt{K}} + \frac{LD^2}{\sqrt{K}}\Big(1 + \sqrt{\frac{\omega}{pn}}\Big)\Big).}
\end{eqnarray*}
\end{corollary}
One can find the full statement with its proof in \ref{sec:MARINA} of Appendix. 

\textbf{EF21 FW.} To complement our results of distributed methods with compression we now introduce biased compressors:
\begin{definition}
\label{def:bcomp}
Map $\C:\mathbb{R}^d\rightarrow\mathbb{R}^d$ is a biased compression operator, if there exist a constant  $\delta \geq 1$ such that for all  $x\in\mathbb{R}^d$ $$\mathbb{E}[\|\C(x) - x\|^2] \leq \left(1-\frac{1}{\delta}\right)\|x\|^2.$$
\end{definition}
This is a broader class than unbiased operators. Here, for example, one can find the greedy choice of coordinates \citep{alistarh2018convergence}, sparse decompositions \citep{vogels2019powersgd} and other operators interesting in practice \citep{beznosikov2022biased}. Intuition suggests that using biased/greedy compressors may improve convergence over unbiased operators. But biased compression are less "suitable" in theory than unbiased ones. Indeed, one can construct a simple convex quadratic problem for which distributed SGD with Top1 compression diverges exponentially fast \citep{beznosikov2022biased}. This issue can be resolved using error compensation technique \citep{seide2014-bit, DBLP:journals/corr/abs-1909-05350, qian2020error}. Then we consider one of the state-of-art algorithms with error compensation techique for biased compression \citep{NEURIPS2021_231141b3}:
\begin{align}
    \label{eq:ef21}
    &g^{k}_i = g^{k-1} + \C(\nabla f_i(x^{k}) - g_i^{k-1}),
    \\
    \notag
    &g^{k} = \textstyle{\frac{1}{n}\sum_{i=1} ^n g_i^k = g^{k-1} + \frac{1}{n}\sum_{i = 1}^{n}\C(\nabla f_i(x^{k}) - g_i^{k-1}).}
\end{align}
Here as in \eqref{eq:marina} the computation of $g^k_i$ takes place on the local devices and the computation of $g^k$ -- on the server and for that only compressed differences $\C(\nabla f_i(x^{k+1}) - g_i^k)$ are needed. Note that \eqref{eq:ef21} does not send uncompressed packages at all.
\begin{lemma}\label{lem:ef21}
    Under Assumptions \ref{as:lip}, \ref{as:conv}, \ref{as:lip:local} the algorithm \eqref{eq:stoch_fw}+\eqref{eq:ef21} satisfies Assumption \ref{as:key} with: $\rho_1 = 1$, $A = 1$, $B = 0$, $C = 0$, $\sigma_k^2 = \frac{1}{n}\sum_{i=1}^n\|g_i^{k} - \nabla f_i(x^{k})\|^2$, $\rho_2 = \frac{1}{2\delta}$, $E = 2\delta\widetilde{L}^2$.
\end{lemma}
\begin{corollary}\label{crl:ef21}
    For the algorithm \eqref{eq:stoch_fw}+\eqref{eq:ef21} in the convex and non-convex cases respectively the following convergences take place:
    \begin{eqnarray*}
        &        \hspace{-0.3cm}\textstyle{
        \mathbb{E}[\text{{\small $f(x^{K}) - f^*$}}] = \mathcal{O}\Big(\Delta_0 \exp(-\frac{K}{8\delta}) + \frac{LD^2}{K}\Big[1 + \frac{\widetilde{L}}{L}\delta\Big]\Big),
        }
        \\
        &        \hspace{-0.3cm}\textstyle{
        \mathbb{E}\Big[\min\limits_{0 \leq k \leq K-1} \textbf{gap}(x^k)\Big] = \mathcal{O}\Big(\frac{\Delta_0}{\sqrt{K}} + \frac{LD^2}{\sqrt{K}}\Big[1 + \frac{\widetilde{L}}{L}\delta\Big]\Big).
        }
    \end{eqnarray*}
\end{corollary}
The full statement together with its proof is provided in \ref{sec:EF21}. One more method with compression will be presented in the next section. A comparison of the presented and already existing methods is provided in Section \ref{sec:comparison_distr}. 
\vspace{-0.2cm}
\subsection{Combinations of different approaches} \label{sec:combination_main}
\vspace{-0.2cm}
In this section, we show that combinations of already presented methods can also be analyzed using Assumption \ref{as:key}.

\textbf{SAGA SARAH FW.}
Here we combine ideas of two approaches: \texttt{SARAH FW} and \texttt{SAGA FW}. Such method is preferable because it obtains benefits of both methods: better rates from SARAH and missing full gradient calculations from SAGA. We get the following gradient estimation:
\begin{align}
        g^{k} =& \textstyle{\frac{1}{b} \sum_{i\in S_k} [\nabla f_i(x^{k}) - \nabla f_i(x^{k-1})] + (1 - \lambda)g^{k-1} +}\notag\\
        &\textstyle{+ \lambda \left( \frac{1}{b}\sum_{i \in S_k} [\nabla f_i(x^{k-1}) - y_i^k] + \frac{1}{n}\sum_{j=1}^n y_j^k\right)},\notag\\
        \label{eq:saga_sarah}
    y_i^{k+1} =& 
    \begin{cases}
        \nabla f_i(x^k), &i\in S_k,\\
        y_i^k, & i\notin S_k,
    \end{cases}
\end{align}
where $\lambda = \tfrac{b}{2n}$ and batches of indexes $S_k$ size of $b$ are generated uniformly and independently.
\begin{lemma}
\label{lem:saga_sarah}
Under Assumptions \ref{as:lip}, \ref{as:conv}, \ref{as:lip:local} the algorithm \eqref{eq:stoch_fw}+\eqref{eq:saga_sarah} satisfies Assumption \ref{as:key} with: 
$\rho_1 = \frac{b}{2n}$, $A = \frac{b}{2n^2}$, $B = \frac{2\widetilde{L}^2}{b}$, $C = 0$, $\sigma_k^2 = \frac{1}{n}\sum_{j=1}^n \|\nabla f_j(x^{k}) - y_j^{k+1}\|^2 $, $\rho_2 = \frac{b}{2n}$, $E = \frac{2n\widetilde{L}^2}{b}$.
\end{lemma}
\begin{corollary}
\label{crl:saga_sarah}
For the algorithm \eqref{eq:stoch_fw}+\eqref{eq:saga_sarah} in the convex and non-convex cases the following convergences take place:
\begin{eqnarray*}
    &\hspace{-0.3cm}\mathbb{E}[\text{{\small $f(x^{K}) - f^*$}}] = \mathcal{O}\left(\Delta_0 \exp\left(-\frac{bK}{8n}\right) + \frac{LD^2}{K}\left[1 + \frac{\widetilde{L}}{L} \frac{\sqrt{n}}{b}\right]\right),
    \\
    &\hspace{-0.3cm}\EE\left[\min\limits_{0 \leq k \leq K-1} \textbf{gap}(x^k)\right] = \mathcal{O}\left(\frac{\Delta_0}{\sqrt{K}} + \frac{LD^2}{\sqrt{K}}\left[1 + \frac{\widetilde{L}}{L} \frac{\sqrt{n}}{b}\right]\right).
\end{eqnarray*}
\end{corollary}
The full statement together with its proof can be found in Section \ref{sec:SAGA SARAH}. The estimates from Corollary \ref{crl:saga_sarah} are the same as for \texttt{SARAH FW} (Corollary \ref{crl:sarah}) and this estimates are the best in the literature on Stochastic Conditional Gradient -- see Section \ref{sec:comparison_stoch} for more details.

Also different combinations of techniques can be considered in the distributed case. In all approaches from Section \ref{sec:distr_main}, we considered that we have a full gradient available for each $f_i$. But for example we can assume that $f_i$ also has the finite-sum form $f_i = \sum_{j=1}^m f_{ij}$. Thus one can add the techniques from Section \ref{sec:stoch_main}. We leave this in Appendix \ref{sec:missing}. In the main part, we consider another illustration.

\textbf{Q-L-SVRG FW.} 
As we noted above, MARINA \eqref{eq:marina} is based on SARAH. Then we can do exactly the same manipulation with L-SVRG and convert it from a method for stochastic sum-type problem into a method with compression:
\begin{align}
\label{eq:lsvrg+c}
        w^{k+1} =& 
            \begin{cases}
                x^k, &\text{ with probability } ~~~ p,\\
                w^k, &\text{ with probability }1-p,
            \end{cases}
        \\g^{k+1} =& \textstyle{\frac{1}{n}\sum_{i = 1}^{n}\mathcal{Q}(\nabla f_i(x^{k+1}) - \nabla f_i(w^{k+1})) + \nabla f(w^{k+1})},\notag
\end{align}
The same way as in \eqref{eq:marina}, to compute $g^{k+1}$ on the server we typically need only compress differences $\mathcal{Q}(\nabla f_i(x^{k+1}) - \nabla f_i(w^{k+1}))$. But with small probability $p$ we need to transfer on the server full packages and compute new $\nabla f(w^{k+1})$. 
\begin{lemma}
\label{lem:lsvrg_c}
Under Assumptions \ref{as:lip}, \ref{as:conv}, \ref{as:lip:local} the algorithm \eqref{eq:stoch_fw}+\eqref{eq:lsvrg+c} satisfies Assumption \ref{as:key} with:
$\rho_1 = 1$, $A = \omega \widetilde{L}^2(1-\frac{p}{2})$, $B = \omega \widetilde{L}^2\left(1 + \frac{8(1-p)}{p}\right)$, $C = 0$, $\sigma_k^2 = \|x^k-w^k\|^2$, $\rho_2 = \frac{p}{2}$, $E =1 + \frac{8(1-p)}{p}$.
\end{lemma}

\begin{corollary}
\label{crl:lsvrg_c}
For the algorithm \eqref{eq:stoch_fw}+\eqref{eq:lsvrg+c} in the convex and non-convex cases accordingly the following convergences take place:
\begin{eqnarray*}
   &\hspace{-0.3cm}\mathbb{E}[\text{{\small $f(x^{K}) - f^*$}}] = \textstyle{\mathcal{O}\left(\Delta_0\exp\left(-\frac{Kp}{8}\right) + \frac{\widetilde{L}D^2}{K+\frac{1}{p}}\left[1 + \frac{\sqrt{\omega}}{p}\right]\right),} 
   \\
   &\hspace{-0.3cm}\EE\Bigl[\min\limits_{0 \leq k \leq K-1} \textbf{gap}(x^k)\Bigr] = \textstyle{\mathcal{O}\left(\frac{\Delta_0}{\sqrt{K}} + \frac{\widetilde{L}D^2}{\sqrt{K}}\left[1 + \frac{\sqrt{\omega}}{p}\right]\right).}
\end{eqnarray*}
\end{corollary}
See the proof in Section \ref{sec:lsvrgc}.
\vspace{-0.2cm}
\section{EXPERIMENTS}\label{sec:exp}
\vspace{-0.2cm}
In this section, we present experimental results to support our theoretical findings.

\textbf{Setup.} We consider the logistic regression problem: 
\begin{equation}
\label{setup}
\hspace{-0.2cm}
    \textstyle{
        \min_{x\in \X}f(x) = \frac{1}{n}\sum_{i=1}^n \log (1 + \exp( - b_i \cdot x^T a_i)),
    }
\end{equation}
where $x$ are weights and $\{a_i, b_i\}_{i=1}^n$ are training data samples with $a_i \in \R^d$ and $b_i \in \{-1, 1 \}$. We choose $\X$ as the $l_1$-ball with radius $r = 2\cdot10^3$. The linear minimization oracle, i.e. $\arg\min_{s\in\X} ~\< g, s>$, for this constraint set can be computed in the closed-form: $s^* = -\text{sign}(g_i) e_i$ with $i = \arg \max_j |g_j|$. We take LibSVM datasets \citep{10.1145/1961189.1961199}.

For further details about the experiments and additional experiments see Section \ref{sec:adexp}.

\begin{figure}[h] 
\begin{minipage}{0.23\textwidth}
    \centering
    \includegraphics[width=1\textwidth]{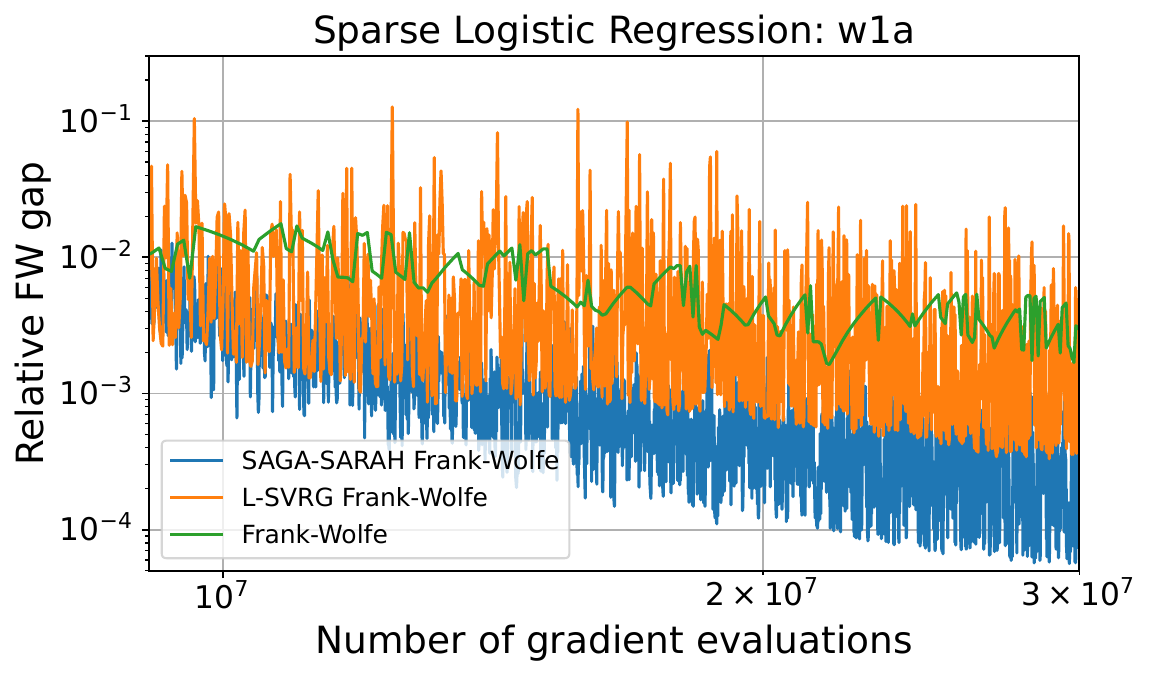}
\end{minipage}
\begin{minipage}{0.23\textwidth}
    \centering
    \includegraphics[width=1\textwidth]{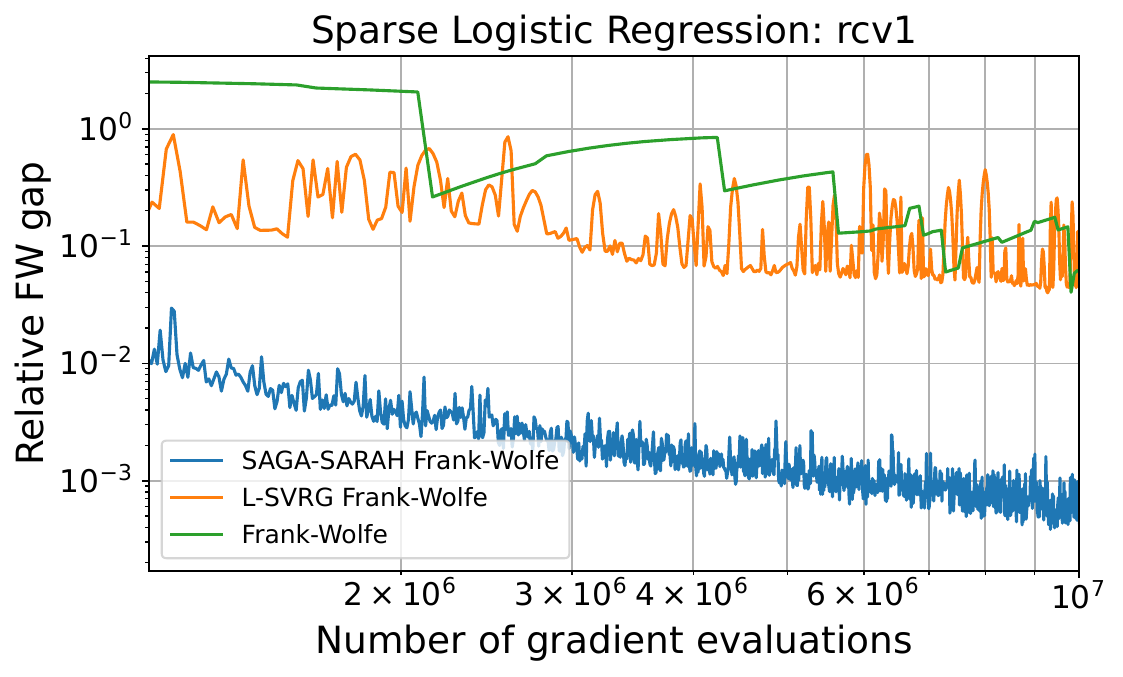}
\end{minipage}
\\
\begin{minipage}{0.01\textwidth}
\quad
\end{minipage}%
\begin{minipage}{0.23\textwidth}
  \centering
\texttt{w1a}
\end{minipage}%
\begin{minipage}{0.23\textwidth}
\centering
\texttt{rcv1}
\end{minipage}%
\vspace{-0.2cm}
\caption{Comparison of methods for solving \eqref{setup} in the stochastic case. \texttt{SARAH FW}, \texttt{SAGA SARAH FW} are considered. The comparison is made on LibSVM datasets \texttt{w1a}, \texttt{rcv1}.}
    \label{fig:stoch}
\end{figure}

\textbf{Stochastic methods.}
In this experiment, we test the performance of proposed stochastic Frank-Wolfe-based methods (\texttt{SAGA SARAH FW}, \texttt{L-SVRG FW}) and compare them to deterministic version of it. The performance is measured in number of full gradients computed. Hyperparameters of the methods are chosen according to theory and all methods starts from zero. The results are provided in Figure \ref{fig:stoch}. One can note that \texttt{SAGA SARAH FW} outperforms competitors. It confirms our theoretical conclusions perfectly.
\newline
In addition, we compared the methods with stochastic baselines chosen from \citep{lu2021generalized}, \citep{mokhtari2020stochastic}, \citep{negiar2020stochastic}:
\begin{figure}[h] 
\begin{minipage}{0.23\textwidth}
    \centering
    \includegraphics[width=1\textwidth]{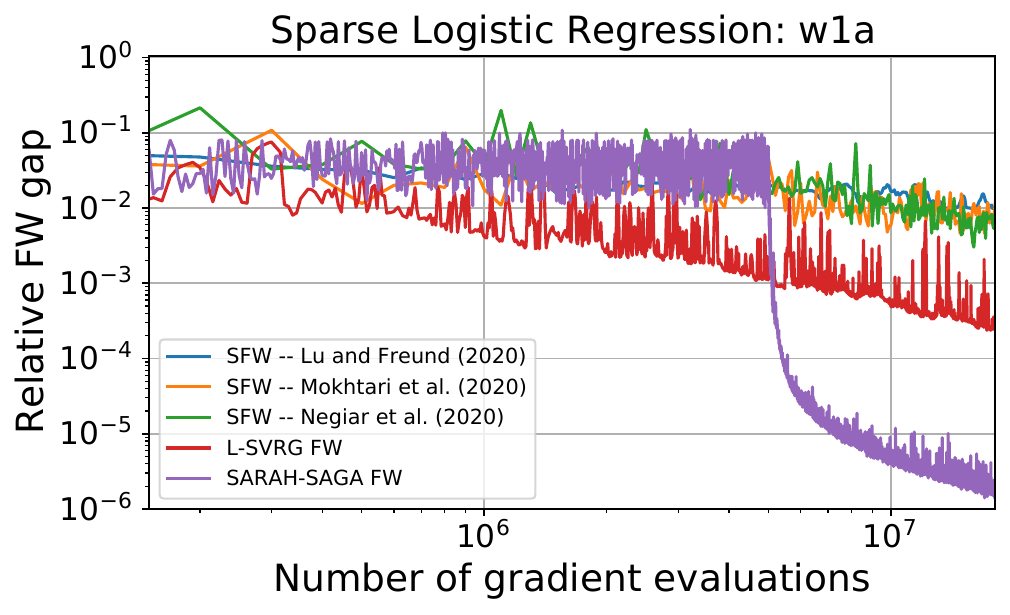}
\end{minipage}
\begin{minipage}{0.23\textwidth}
    \centering
    \includegraphics[width=1\textwidth]{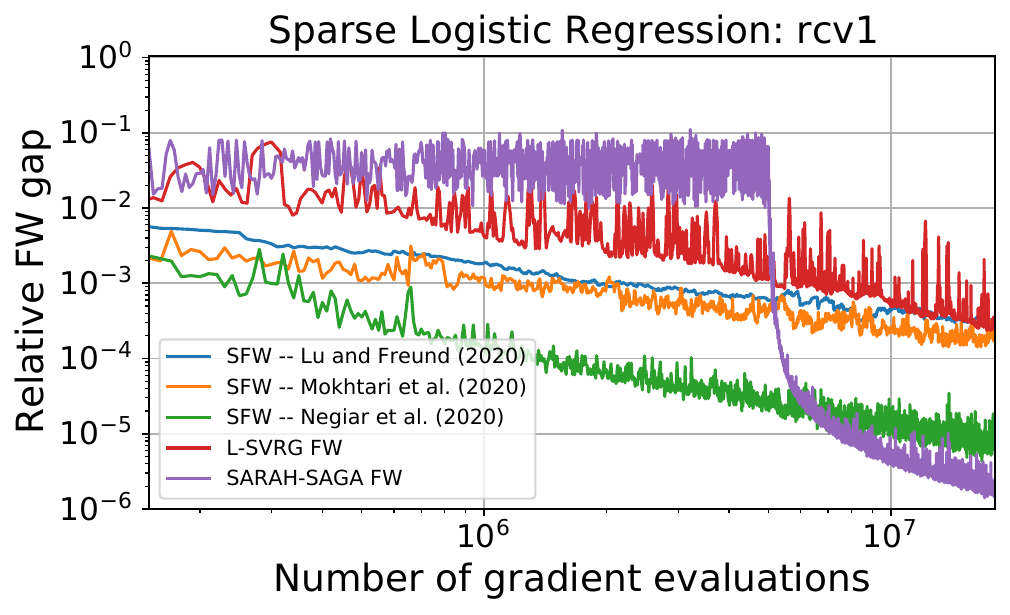}
\end{minipage}
\\
\begin{minipage}{0.01\textwidth}
\quad
\end{minipage}%
\begin{minipage}{0.23\textwidth}
  \centering
\texttt{w1a}
\end{minipage}%
\begin{minipage}{0.23\textwidth}
\centering
\texttt{rcv1}
\end{minipage}%
\vspace{-0.2cm}
\caption{Comparison of stochastic methods with baselines. \texttt{SARAH FW}, \texttt{SAGA SARAH FW} are considered. The comparison is made on LibSVM datasets \texttt{w1a}, \texttt{rcv1}.}
\end{figure}

Finally for this section of methods we introduce a comparison with methods for constrained optimization with Euclidean projection instead of linear minimization. It is presented in Section \ref{sec:adexp} of Appendix.

\begin{figure}[h]
\begin{minipage}{0.23\textwidth}
    \centering
    \includegraphics[width=1\textwidth]{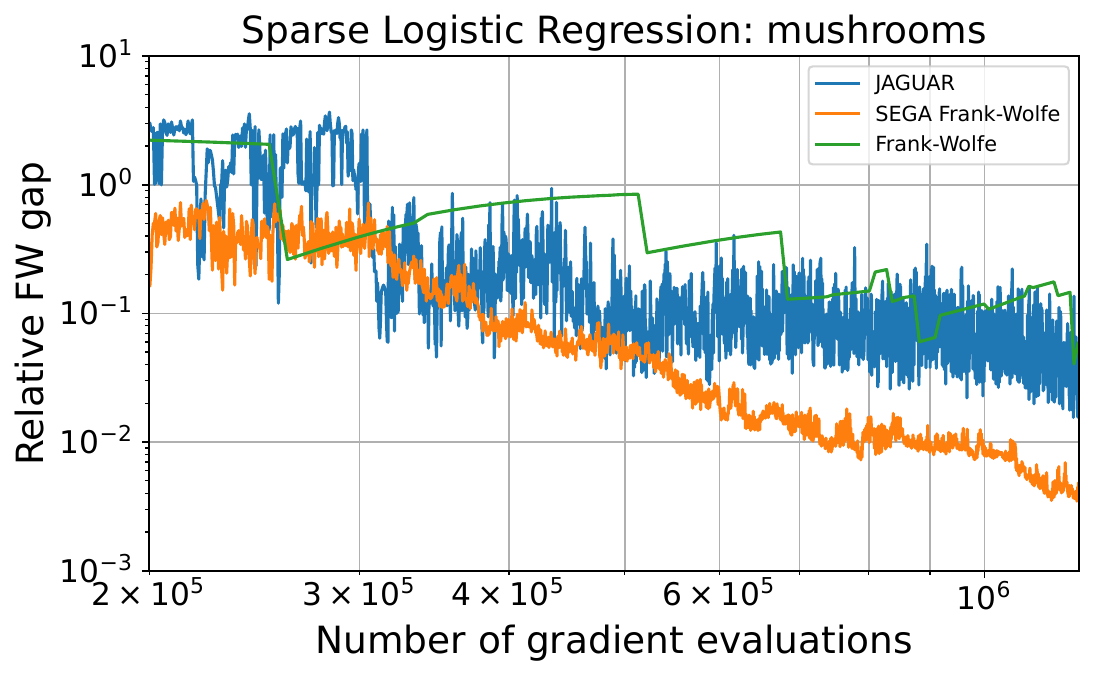}
\end{minipage}
\begin{minipage}{0.23\textwidth}
    \centering
    \includegraphics[width=1\textwidth]{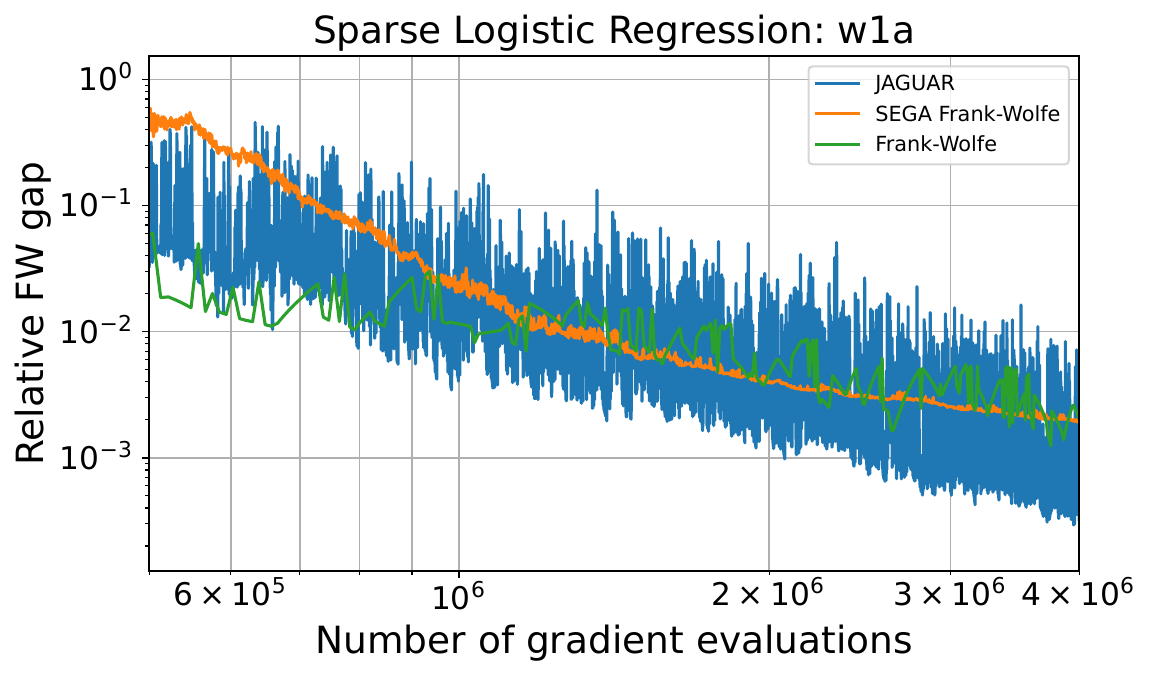}
\end{minipage}
\\
\begin{minipage}{0.01\textwidth}
\quad
\end{minipage}%
\begin{minipage}{0.23\textwidth}
  \centering
\texttt{mushrooms}
\end{minipage}%
\begin{minipage}{0.23\textwidth}
\centering
\texttt{w1a}
\end{minipage}%
\vspace{-0.2cm}
\caption{Comparison of methods for solving \eqref{setup} in the coordinate case. \texttt{SEGA FW}, \texttt{JAGUAR} are considered. The comparison is made on LibSVM datasets \texttt{mushrooms}, \texttt{w1a}.}
    \label{fig:coord}
\end{figure}

\begin{figure}[h]
\begin{minipage}{0.23\textwidth}
    \centering
    \includegraphics[width=1\textwidth]{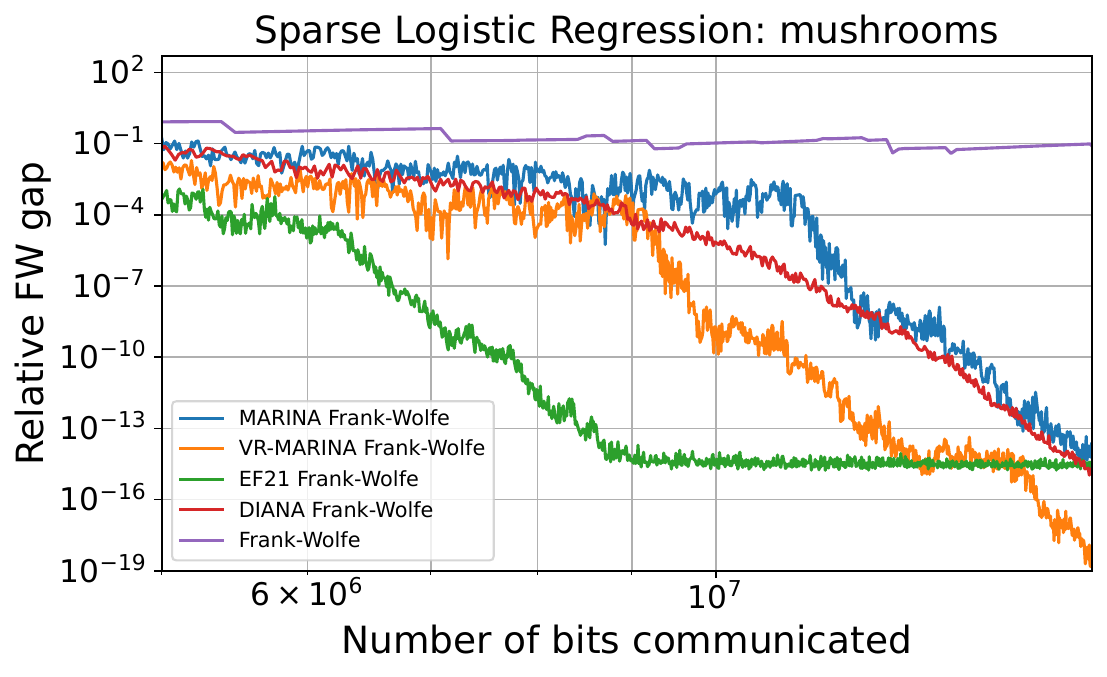}
\end{minipage}
\begin{minipage}{0.23\textwidth}
    \centering
    \includegraphics[width=1\textwidth]{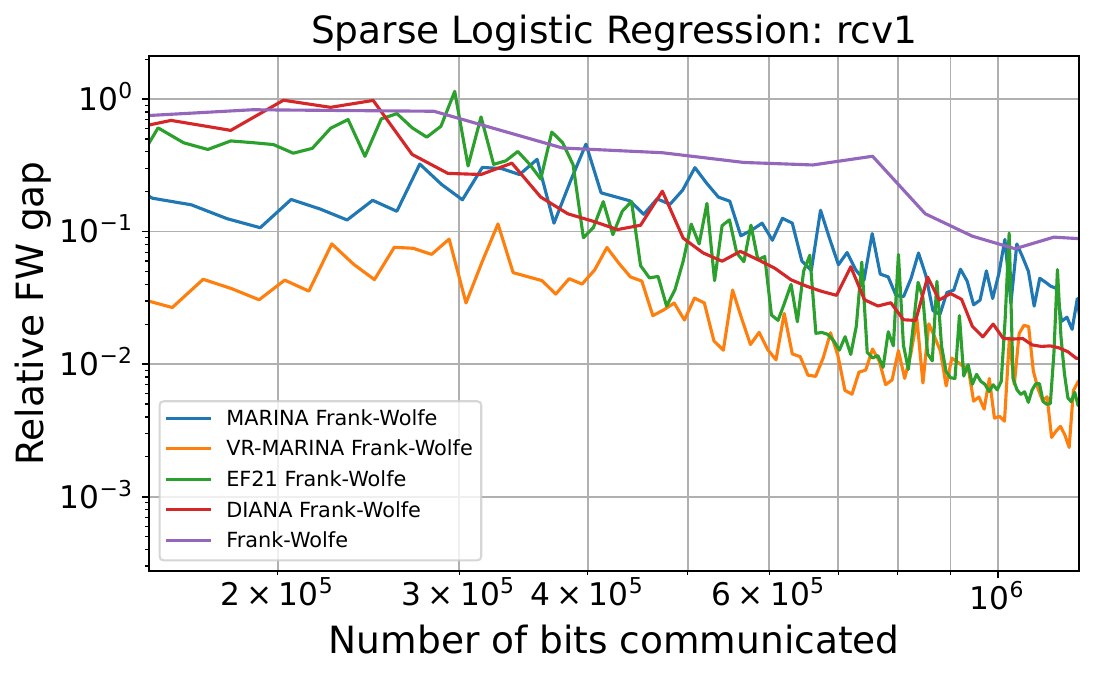}
\end{minipage}
\\
\begin{minipage}{0.01\textwidth}
\quad
\end{minipage}%
\begin{minipage}{0.23\textwidth}
  \centering
\texttt{mushrooms}
\end{minipage}%
\begin{minipage}{0.23\textwidth}
\centering
\texttt{rcv1}
\end{minipage}%
\vspace{-0.2cm}
\caption{Comparison of methods for solving \eqref{setup} in the distributed case. \texttt{MARINA FW}, \texttt{VR-MARINA FW}, \texttt{EF21 FW}, \texttt{DIANA FW} are considered. The comparison is made on LibSVM datasets \texttt{mushrooms}, \texttt{rcv1}.}
    \label{fig:distrib}
\end{figure}

\textbf{Coordinate methods.} In the second experiment, we check the convergence of alleged coordinate methods: \texttt{SEGA FW}, \texttt{JAGUAR} and compare them to the original Frank-Wolfe method. The performance is measured in number of full gradients computed. Hyperparameters of the methods are chosen according to theory and all methods starts from zero. The results are presented in Figure \ref{fig:coord}. 
\vspace{0.4cm}
\newline
\textbf{Distributed methods.}
In the last experiment, we consider the distributed setup of \eqref{setup}, in which we assume that the information about $f_i$ is available for worker $i$ only. Here we compare distributed methods proposed in this paper: \texttt{MARINA Frank-Wolfe}, \texttt{VR-MARINA Frank-Wolfe}, \texttt{DIANA Frank-Wolfe} and \texttt{EF21 Frank-Wolfe}. The performance is measured in number of bits communicated from workers to the server. Hyperparameters of the methods are chosen according to theory and all methods starts from zero. In \texttt{MARINA} and \texttt{DIANA} algorithms RandK (random sparsification) compression is used, while in \texttt{EF21} TopK ("greedy" sparsification) compression is implemented. Convergence performance is shown in Figure \ref{fig:distrib}. The advantage of using compressed communication is clearly observable in every case.

\subsubsection*{Acknowledgements}

The work of R. Nazykov, A. Shestakov and V. Solodkin was supported by a grant for research centers in the field of artificial intelligence, provided by the Analytical Center for the Government of the Russian Federation in accordance with the subsidy agreement (agreement identifier 000000D730321P5Q0002) and the agreement with the Moscow Institute of Physics and Technology dated November 1, 2021 No. 70-2021-00138.

\bibliography{ref}

\appendix
\onecolumn
\newpage
\part*{Supplementary Material}
\tableofcontents

\newpage 

\section{MISSING COMPARISON AND DETAILS} \label{sec:comparison}

\subsection{Stochastic methods} \label{sec:comparison_stoch}

\renewcommand{\arraystretch}{2}
\begin{table*}[h]
    \centering
\captionof{table}{Summary of complexity results for finding an $\varepsilon$-solution \textbf{stochastic finite-sum non-distributed constrained minimization problems} \eqref{eq:finsum} with $n$ terms by projection free methods. Convergence is measured by the functional distance to the solution in the convex case and by the gap function in the non-convex case. Complexities are given in terms of the number of stochastic gradient calls.
\\
\colorbox{bgcolor2}{blue} = results of our paper (in particular, \cite{reddi2016stochastic} do not give a proof in the convex setting).}
\vspace{-0.2cm}
    \label{tab:comparison_stoch}   
    \small
  \begin{threeparttable}
    \begin{tabular}{|c|c|c|c|c|}
    \hline
    \textbf{Method} & \textbf{Link} & \textbf{Reference} & \textbf{CVX} & \textbf{nCVX}
    
    \\\hline
    \texttt{Frank-Wolfe} & Alg. 1 & \citep{frank1956algorithm}, \citep{lacostejulien2016convergence} & $\cO(\frac{n}{\e})$  & $\cO(\frac{n}{\e^2})$ 
    \\\hline
    \texttt{SVRG FW} & Alg. 3 & \cite{reddi2016stochastic} & \ding{55} &  $\mathcal{ O} ( n + \frac{n^{2/3}}{\varepsilon^2} )$
    \\\hline
    \texttt{SVRG FW} & Alg. 1 & \cite{hazan2016variance} & $\cO \left( n + \frac{1}{\varepsilon^2}\right)$ & \ding{55}
    \\\hline
    \texttt{SPIDER FW} & Alg. 2 & \citep{pmlr-v97-yurtsever19b} & $\cO(\frac{1}{\e^2})$ & $\cO(n + \frac{\sqrt{n}}{\e^2})$
    \\\hline
    \texttt{SFW} & Alg. 1 & \cite{negiar2020stochastic} & $\cO \left( \frac{n}{\varepsilon}\right)$ & \ding{55} 
    \\\hline
    \texttt{GSFW} & Alg. 1 & \citep{lu2021generalized} & $\cO(\frac{n}{\e})$ &  \ding{55}
    \\\hline
    \texttt{SVRG FW} & Alg. 3 & \cite{weber2022projection}  & \ding{55} & $\mathcal{O} ( n + \frac{n^{2/3}}{\varepsilon^2} )$
    \\\hline
    \texttt{SAGA FW} & Alg. \ref{alg:saga} & \citep{reddi2016stochastic} & \cellcolor{bgcolor2}{$\widetilde{\cO}(n + \frac{n^{2/3}}{\e})$ \tnote{{\color{blue}(1)}} } & $\cO(n +\frac{n^{2/3}}{\e^2})$ \tnote{{\color{blue}(1)}} 
    \\\hline
    \texttt{SARAH FW} & Alg. \ref{alg:sarah} & \citep{Beznosikov2023SarahFM} & $\widetilde{\cO}(n + \frac{\sqrt{n}}{\e})$ \tnote{{\color{blue}(1)}}  & $\cO(n + \frac{\sqrt{n}}{\e^2})$ \tnote{{\color{blue}(1)}} 
    \\\hline
    \texttt{SAGA SARAH FW} & Alg. \ref{alg:sagasarah} & \citep{Beznosikov2023SarahFM} & $\widetilde{\cO}(n +\frac{\sqrt{n}}{\e})$ \tnote{{\color{blue}(1)}}  & $\cO(n + \frac{\sqrt{n}}{\e^2})$ \tnote{{\color{blue}(1)}} 
    \\\hline
    \cellcolor{bgcolor2}{\texttt{L-SVRG FW}} & \cellcolor{bgcolor2}{Alg. \ref{alg:lsvrg}} & \cellcolor{bgcolor2}{\textbf{NEW}, \citep{pmlr-v117-kovalev20a}} & \cellcolor{bgcolor2}{$\widetilde{\cO}(n + \frac{n^{2/3}}{\e})$ \tnote{{\color{blue}(1)}}  } & \cellcolor{bgcolor2}{$ \cO(n +\frac{n^{2/3}}{\e^2})$ \tnote{{\color{blue}(1)}} }
    \\\hline
    \end{tabular}   
    \begin{tablenotes}
    {\small   
    \item [] \tnote{{\color{blue}(1)}} In the main part of the paper, we give these results with depending on the parameters $p$ and $b$. The optimal choice of them are given in the original paper \citep{Beznosikov2023SarahFM} (\texttt{SARAH FW} and \texttt{SAGA SARAH FW}) or in the corresponding subsections of Section \ref{sec:missing} (\texttt{SAGA FW} and \texttt{L-SVRG FW}).
    }
\end{tablenotes}    
    \end{threeparttable}
\vspace{-0.3cm}
\end{table*}

It is important to note that in Table \ref{tab:comparison_stoch} we only covers works on the finite-sum stochastic optimization. However, there are many papers where the authors also consider a stochastic version of the FW algorithm, but under the (additional) assumption of bounded variance of the stochastic gradients \citep{hazan2012projection, doi:10.1137/140992382, reddi2016stochastic, pmlr-v80-qu18a, pmlr-v97-yurtsever19b, mokhtari2020stochastic, zhang2020one}.

\newpage

\subsection{Coordinate methods} \label{sec:comparison_coord}

\renewcommand{\arraystretch}{2}
\begin{table*}[h]
    \centering
\captionof{table}{Summary of complexity results for finding an $\varepsilon$-solution \textbf{constrained minimization problems} with dimension $d$ by \textbf{coordinate} projection free methods. Convergence is measured by the functional distance to the solution in the convex case and by the gap function in the non-convex case. Complexities are given in terms of the number of coordinates derivatives computed.
\\
\colorbox{bgcolor2}{blue} = results of our paper.}
\vspace{-0.2cm}
    \label{tab:comparison_coord}   
    \small
  \begin{threeparttable}
    \begin{tabular}{|c|c|c|c|c|}
    \hline
    \textbf{Method} & \textbf{Link} & \textbf{Reference} & \textbf{CVX} & \textbf{nCVX}
    \\\hline
    \texttt{BCFW} & Alg. 3 & \citep{pmlr-v28-lacoste-julien13} & ${\cO}(d + \frac{d}{\e})$ & \ding{55}
    \\\hline
    \texttt{SGFFW} & Alg. 1 & \citep{pmlr-v89-sahu19a}\tnote{{\color{blue}(1)}} & ${\cO}(d + \frac{d}{\e})$ & \ding{55}
    \\\hline
    \texttt{SGFFW} & Alg. 2 & \citep{pmlr-v89-sahu19a}\tnote{{\color{blue}(1)}} & ${\cO}(d + \frac{d}{\e^3})$ & ${\cO}(d + \frac{d^{4/3}}{\e^4})$
    \\\hline
    \cellcolor{bgcolor2}{\texttt{SEGA FW}} & \cellcolor{bgcolor2}{Alg. \ref{alg:sega}} & \cellcolor{bgcolor2}{\textbf{NEW}, \citep{NEURIPS2018_fc2c7c47} } & \cellcolor{bgcolor2}{$\widetilde{\cO}(d + \frac{d\sqrt{d}}{\e})$} & \cellcolor{bgcolor2}{$\cO(d + \frac{d^3}{\e^2})$}
    \\\hline
    \cellcolor{bgcolor2}{\texttt{JAGUAR}} & \cellcolor{bgcolor2}{Alg. \ref{alg:aytona}} & \cellcolor{bgcolor2}{\textbf{NEW}} & \cellcolor{bgcolor2}{$\widetilde{\cO}(d + \frac{d}{\e})$} &  \cellcolor{bgcolor2}{$\cO(d +\frac{d^2}{\e^2})$}
    \\\hline
    \end{tabular}   
    \begin{tablenotes}
        {\small   
        \item [] \tnote{{\color{blue}(1)}} Zero-order methods.
        }
    \end{tablenotes}   
    \end{threeparttable}
\vspace{-0.3cm}
\end{table*}

\subsection{Distributed methods} \label{sec:comparison_distr}

\renewcommand{\arraystretch}{2}
\begin{table*}[h]
    \centering
\captionof{table}{Summary of complexity results for finding an $\varepsilon$-solution \textbf{distributed constrained minimization problems} \eqref{eq:finsum} with $n$ devices by projection free methods \textbf{with compression} with parameters $\omega$ and $\delta$ (see Definitions \ref{def:unbcomp} and \ref{def:bcomp}). Convergence is measured by the functional distance to the solution in the convex case and by the gap function in the non-convex case. Complexities are given in terms of the number of transmitted coordinates if we choose RandK and TopK as particular cases of compressors.
\\
\colorbox{bgcolor2}{blue} = results of our paper.}
\vspace{-0.2cm}
    \label{tab:comparison_distrib}   
    \small
\resizebox{\linewidth}{!}{
  \begin{threeparttable}
    \begin{tabular}{|c|c|c|c|c|c|}
    \hline
    \textbf{Method} & \textbf{Link} & \textbf{Reference} & \textbf{Any compression?} & \textbf{CVX} & \textbf{nCVX} 
    \\\hline
    \texttt{dFW} & Alg. 3 & \citep{bellet2015distributed} & \ding{55} & $\cO(\frac{1}{\e})$ &  \ding{55}
    \\\hline
    \cellcolor{bgcolor2}{\texttt{DIANA FW}} & \cellcolor{bgcolor2}{Alg. \ref{alg:diana}} & \cellcolor{bgcolor2}{\textbf{NEW}, \citep{mishchenko2019distributed}} & \cellcolor{bgcolor2}{\ding{51} (unbiased)} & \cellcolor{bgcolor2}{$\widetilde{\cO}(1 + \frac{1}{\e}(\frac{1}{\omega} + \frac{\sqrt{\omega}}{\sqrt{n}}))$} & \cellcolor{bgcolor2}{$\cO(\frac{1}{\e^2}(\frac{1}{\omega} + \frac{\omega^2}{n}))$}
    \\\hline
    \cellcolor{bgcolor2}{\texttt{MARINA FW}} & \cellcolor{bgcolor2}{Alg. \ref{alg:marina}} & \cellcolor{bgcolor2}{\textbf{NEW}, \citep{pmlr-v139-gorbunov21a}} & \cellcolor{bgcolor2}{\ding{51} (unbiased) } & \cellcolor{bgcolor2}{$\widetilde{\cO}(1 + \frac{1}{\e}(\frac{1}{\omega} + \frac{1}{\sqrt{n}}))$} & \cellcolor{bgcolor2}{$\cO(\frac{1}{\e^2}(\frac{1}{\omega} + \frac{\omega}{n}))$}
    \\\hline
    \cellcolor{bgcolor2}{\texttt{VR MARINA FW}} & \cellcolor{bgcolor2}{Alg. \ref{alg:vrmarina}} & \cellcolor{bgcolor2}{\textbf{NEW}, \citep{pmlr-v139-gorbunov21a}} & \cellcolor{bgcolor2}{\ding{51} (unbiased)} & \cellcolor{bgcolor2}{$\widetilde{\cO}(1 + \frac{1}{\e}(\frac{1}{\omega} + \frac{1}{\sqrt{n}}))$} & \cellcolor{bgcolor2}{$\cO(\frac{1}{\e^2}(\frac{1}{\omega} + \frac{\omega}{n}))$}
    \\\hline
     \cellcolor{bgcolor2}{\texttt{EF21 FW}} & \cellcolor{bgcolor2}{Alg. \ref{alg:ef21}} & \cellcolor{bgcolor2}{\textbf{NEW}, \citep{NEURIPS2021_231141b3}} & \cellcolor{bgcolor2}{\ding{51}} & \cellcolor{bgcolor2}{$\widetilde{\cO}(1 + \frac{1}{\e})$} & \cellcolor{bgcolor2}{$\cO(\frac{\delta}{\e^2})$}
    \\\hline
    \cellcolor{bgcolor2}{\texttt{Q-L-SVRG FW}} & \cellcolor{bgcolor2}{Alg. \ref{alg:lsvrgc}} & \cellcolor{bgcolor2}{\textbf{NEW}} & \cellcolor{bgcolor2}{\ding{51} (unbiased)}  & \cellcolor{bgcolor2}{$\widetilde{\cO}(1 + \frac{1}{\e}(\frac{1}{\omega} + \frac{\sqrt{\omega}}{\sqrt{n}}))$} & \cellcolor{bgcolor2}{$\cO(\frac{1}{\e^2}(\frac{1}{\omega} + \frac{\omega^2}{n}))$}
    \\\hline
    \end{tabular}   
    \end{threeparttable}
}
\vspace{-0.3cm}
\end{table*}

Here one can also highlight two papers on distributed Frank-Wolfe algorithms, but without compression \citep{pmlr-v48-wangd16, hou2022distributed}.
\newpage
\section{MISSING METHODS, DETAILS AND PROOFS} \label{sec:missing}
In this section, we first provide complete proofs of our two main Theorems \ref{th:convex} and \ref{th:nonconvex}. Then we discuss the zoo of special cases, in particular we provide the full listing of algorithms, detailed convergence rates and proofs for methods from Section \ref{sec:zoo}. Moreover, we present some statements and algorithms for them that are not encountered in the main part.

\subsection{Technical facts} 

In our proofs, we often apply following inequalities that hold for any $a, b \in \R^d$ and $\alpha > 0$:
\begin{equation}
    \label{eq:young}
    \|a + b\|^2\leq(1 + \alpha)\|a\|^2 + \Big(1 + \frac{1}{\alpha}\Big)\|b\|^2,
\end{equation}
\begin{equation}
    \label{eq:young2}
    2\langle a, b\rangle \leq \frac{1}{\alpha}\|a\|^2 + \alpha\|b\|^2.
\end{equation}\vspace{-0.5cm}
\begin{lemma}
\textup{(Lemma 1.2.3 from \citep{10.5555/2670022})}\textbf{.} Let the Assumption \ref{as:lip} be satisfied. Then for all $x, y \in \mathbb{R}^d$:
\begin{equation}
    \label{Nesterov}
    |f(y) - f(x) -  \langle \nabla f(x),y - x \rangle| \leq \frac{L}{2}\|x-y\|^2.
\end{equation}
\end{lemma}
\begin{lemma}\label{lem:2moment_uni_batch}
\textup{(Lemma A.1 from \citep{lei2017non})}\textbf{.} Let $x_1, \ldots, x_N \in \mathbb{R}^d$ be arbitrary vectors with
$$\sum\limits_{i=1}^N x_i = 0.$$
Further let $S$ be a uniform subset of $\left[N\right]$ with size b. Then
\begin{equation}
    \EE\left\|\frac{1}{b}\sum\limits_{i\in S}x_i\right\|^2 \leq \frac{1}{bN}\sum\limits_{i=1}^N \|x_i\|^2.
\end{equation}
\end{lemma}
Our further analysis relies heavily on work conducted in \citep{stich2019unified}, but it seems to us that there is a typo in the main part. Furthermore, as we will deal with nonreducible "noises", we will need to adjust proofs by a little. For convenience of the proof, we split the result into three following lemmas:
\begin{lemma}\label{lem:b2}
Consider two non-negative sequences $\lbrace r_k \rbrace,\, \lbrace \eta_k \rbrace$, that satisfy\\
\begin{equation}
    \label{Orig_seq}
    r_{k+1} \leq (1 - \eta_k)r_k + a\eta_k^2 + b,
\end{equation} 
\\for all $k \geq 0$, constants $a\geq 0,\, b \geq 0$ and for positive stepsizes $\lbrace \eta_k\rbrace$ with $ \eta_k \leq \frac{1}{d},\, d > 1$. Then there exists $\eta$, such that $\forall k \geq 0,\, \eta_k \equiv \eta$
\begin{equation}
    \label{non_decr_step}
    r_K < r_0 \exp \left[- \frac{K}{d}\right] + \frac{a}{d} + bK.
\end{equation}
\end{lemma}
\textbf{Proof:}
\\Set $\eta_k \equiv \eta,\: \forall k \geq 0$. Unroll inequalities:
\begin{equation*}
    \begin{split}
        r_K &\leq (1-\eta)r_{K-1} + a\eta^2 + b \leq (1-\eta)^Kr_0 + a\eta^2 \sum\limits_{i = 0}^{K-1}(1-\eta)^i + bK \\
        &\leq (1-\eta)^T + a\eta^2\sum\limits_{k = 0}^{\infty}(1-\eta)^k + bK \leq (1-\eta)^Kr_0 + a\eta + bK.
    \end{split}
\end{equation*}\\
Put $\eta = \frac{1}{d}$ and use, that for $\forall x > 0,\; 0 < a < 1 \rightarrow (1-a)^x < \exp(-ax)$. Hence, \\
\begin{equation}
    r_K < r_0 \exp \left[- \frac{K}{d}\right] + \frac{a}{d} + bK.
\end{equation}
\EndProof{}
\begin{lemma}\label{lem:b3}
Let the non-negative sequence $\lbrace r_k \rbrace$ satisfy the conditions of Lemma \ref{lem:b2}. Then there are decreasing stepsizes $\eta_k = \frac{2}{2d + k}$, such that 
\begin{equation}
    \label{simple_dec_step}
    r_{K+1} \leq \frac{4d^2r_0}{(2d+K)^2} + 4a\frac{K+1}{(2d+K)^2} + \frac{b}{6} \frac{(K+1)(24d^2 + 12dK + 2K^2 + K)}{(2d+K)^2}.
\end{equation}
\end{lemma}
\textbf{Proof:}\\
Rearrange (\ref{Orig_seq}):
\begin{equation*}
    0 \leq (1-\eta_k)r_k - r_{k+1} + a\eta_k^2 + b.
\end{equation*}
Divide both parts by $\eta_k$:\\
\begin{equation*}
    0 \leq \frac{1-\eta_k}{\eta_k}r_k - \frac{r_{k+1}}{\eta_k} + a\eta_k + \frac{b}{\eta_k} = (2d+k-2)\frac{r_k}{2} - (2d+k)\frac{r_{k+1}}{2} + a\frac{2}{2d+k} + b\frac{2d+k}{2}.
\end{equation*}
Multiply both parts by $2 \cdot (2d+k)$ and use, that $x \cdot (x-2) \leq (x-1)^2$:\\
\begin{equation*}
    \begin{split}
        0 &\leq (2d+k-2)(2d+k)r_k - (2d+k)^2r_{k+1} + 4a + b(2d+k)^2\\
        &\leq (2d+k-1)^2r_k - (2d+k)^2r_{k+1} + 4a + b(2d+k)^2.
    \end{split}
\end{equation*}
We obtain a telescoping sum, hence:
\begin{equation*}
    \begin{split}
        0 &\leq (2d-1)^2r_0 - (2d+k)^2r_{K+1} + 4a(K+1) + b \sum\limits_{i = 0}^K (2d+i)^2\\
        &\leq 4d^2r_0 - (2d+k)^2r_{K+1} + 4a(K+1) + \frac{b(K+1)(24d^2 + 12dK+2K^2+K)}{6}.
    \end{split}
\end{equation*}\\
\EndProof{Finally, we can gain the original inequality by rearranging the terms.}
\begin{lemma}\label{lem:b4}
    Let $\lbrace r_k \rbrace$ satisfy (\ref{Orig_seq}). Then there exist stepsizes $\eta_k$,
    \begin{align*}
    &\text{if ~K } \leq \text{ d}, && \eta_k = \tfrac{1}{d},
    \\
    &\text{if ~K } > \text{ d ~and ~k } < \text{k}_0, && \eta_k = \tfrac{1}{d},
    \\
    &\text{if ~K } > \text{ d ~and ~k } \geq \text{k}_0, && \eta_k = \tfrac{2}{2d + k - k_0},
    \end{align*}
    where $k_0 = \lfloor\frac{K}{2}\rfloor$, such that\\
    \begin{equation}
        \label{tech_lemma}
        r_{K+1} = \mathcal{O}\left(r_0\exp\left(-\frac{K}{2d}\right) + \frac{a}{d+K} + bK\right).
    \end{equation}
\end{lemma}
\textbf{Proof:}\\ 
If $K \leq d$ we'll take $\eta_k \equiv \frac{1}{d}$, hence using Lemma \ref{lem:b2}  we obtain
\begin{equation*}
    r_K \leq r_0 \exp\left[-\frac{K}{d}\right] + \frac{a}{d} + bK \leq r_0 \exp\left[-\frac{K}{d}\right] + \frac{a}{K} + bK.
\end{equation*}\\
As $\frac{a}{K} / \frac{a}{d+K} = \frac{d+K}{K} \leq \frac{2d}{K} \leq 2d ~\text{ and } \exp\big(-\frac{K}{d}\big) / \exp\big(-\frac{K}{2d}\big) = \exp\big(-\frac{K}{2d}\big) < 1,$ thus \eqref{tech_lemma} is correct for $K\leq d.$

If $K > d$ for $k \leq k_0 := \lceil \frac{K}{2}\rceil$ we'll take $\eta_k \equiv \frac{1}{d}$ and for $k > k_0$: $\eta_k = \frac{2}{2d+k-k_0}$. Therefore, according to \eqref{non_decr_step} we obtain
\begin{equation*}
   r_{k_0} \leq r_0 \exp\left(-\frac{k_0}{d}\right) + \frac{a}{d} + bk_0. 
\end{equation*}
Using Lemma \ref{lem:b3} with $k_1 = K - k_0$ we get\\
\begin{equation*}
       r_{K+1} \leq \frac{4d^2r_{k_0}}{(2d+k_1)^2} + 4a\frac{K+1}{(2d+k_1)^2} + \frac{b}{6} \frac{(k_1+1)(24d^2 + 12dk_1 + 2k_1^2 + k_1)}{(2d+k_1)^2}. 
\end{equation*}\\
Combining inequalities we derive\\
\begin{equation*}
    \frac{4d^2r_{k_0}}{(2d+k_1)^2} \leq \frac{4d^2}{(2d+k_1)^2}\left(r_0 \exp\left(-\frac{k_0}{d}\right) + \frac{a}{d} + bk_0\right) \leq r_0\exp\left(-\frac{k_0}{d}\right) + \frac{4ad}{(2d+k_1)^2} + \frac{4bk_0d^2}{(2d+k_1)^2}.
\end{equation*}\\
Since $\frac{K}{2} \leq \lceil \frac{K}{2} \rceil < K$ and $\frac{K}{4} < \lfloor \frac{K}{2} \rfloor < K$ we can replace all $k_0, k_1$ with K in fractions in O-notation. We gain the replacement in the exponent, as $-\lceil \frac{K}{2}\rceil \leq - \frac{K}{2}$. Thus,\\
\begin{equation*}
    \begin{split}
        r_{K+1} &= \mathcal{O}\left(r_0\exp\left(-\frac{K}{2d}\right) + \frac{ad}{(d+K)^2} + \frac{d^2bK}{(d+K)^2}+ \frac{aK}{(d+K)^2} + bK\frac{d^2+K^2}{(d+K)^2}\right) \\
        &= \mathcal{O}\left(r_0\exp\left(-\frac{K}{2d}\right) + \frac{a}{K+d} + bK\frac{d^2+K^2}{(d+K)^2}\right),
    \end{split}
\end{equation*}\\
\EndProof{where second equality is implied by $d,K > 0$, $d^2+K^2 \leq (d+K)^2$.}
\newpage
\subsection{Unified main theorems}\label{sec:uni}
In this section, we provide complete proofs of our main results.
\begin{lemma}
\label{lem:conv}
If $x^k$ is upgraded due to Algorithm (\ref{eq:stoch_fw}), then for all $\alpha > 0$:
\begin{equation*}
    \begin{split}
        \EE[f(x^{k+1}) - f(x^*)] &\leq (1-\eta_k)\EE[f(x^k) - f(x^*)] + \frac{\alpha}{L} \EE[\|\nabla f(x^k) - g^k \|^2] + \frac{L\eta_k^2}{\alpha}\EE[\|s^k - x^*\|^2] + \frac{L\eta_k^2}{2}\EE[\|s^k - x^k\|^2].
    \end{split}
\end{equation*}
\end{lemma}
\vspace{-0.2cm}
\textbf{Proof:}\\
First, we assume that $f(x)$ satisfies Assumptions \ref{as:lip} and \ref{as:conv}:
\begin{eqnarray*}
    f(x^{k+1}) &\leq& f(x^k) + \<\nabla f(x^k), x^{k+1} - x^k> + \frac{L}{2}\|x^{k+1} - x^k\|^2.  
\end{eqnarray*}
With update of $x^{k+1}$, according to Algorithm \eqref{eq:stoch_fw}:
\begin{eqnarray*}
    f(x^{k+1}) - f(x^*)&\leq& f(x^k) - f(x^*) + \eta_k\<\nabla f(x^k), s^k - x^k> + \frac{L\eta_k^2}{2}\|s^k - x^k\|^2 \\ 
    &=&f(x^k) - f(x^*) + \eta_k\<g^k, s^k - x^k> + \eta_k\<\nabla f(x^k) - g^k, s^k - x^k> + \frac{L\eta_k^2}{2}\|s^k-x^k\|^2.
\end{eqnarray*}
The optimal choice of $s^k$ in Algorithm \eqref{eq:stoch_fw} gives $\<g^k, s^k - x^k>\leq\<g^k, x^* - x^k>$. Then 
\begin{eqnarray*}
    f(x^{+1}) - f(x^*)&\leq& f(x^k) - f(x^*) + \eta_k\<g^k, x^* - x^k> + \eta_k\<\nabla f(x^k) - g^k, s^k - x^k> + \frac{L\eta_k^2}{2}\|s^k-x^k\|^2\\
    &=&f(x^k) - f(x^*) + \eta_k\<\nabla f(x^k), x^* - x^k> + \eta_k\<g^k - \nabla f(x^k), x^* - x^k> \\
    &+&\eta_k\<\nabla f(x^k) - g^k, s^k - x^k> + \frac{L\eta_k^2}{2}\|s^k-x^k\|^2\\
    &=&f(x^k) - f(x^*) + \eta_k\<\nabla f(x^k), x^* - x^k> + \eta_k\<\nabla f(x^k) - g^k, s^k - x^*> + \frac{L\eta_k^2}{2}\|s^k - x^k\|^2.
\end{eqnarray*}
Using Young's inequality here we state for any positive $\alpha$:
\begin{eqnarray*}
    f(x^{k+1}) - f(x^*)&\leq& f(x^k)- f(x^*) + \eta_k\<\nabla f(x^k), x^* - x^k> + \frac{\alpha}{L}\|\nabla f(x^k) - g^k\|^2 \\
    &+& \frac{L\eta_k^2}{\alpha}\|s^k - x^*\|^2 + \frac{L\eta_k^2}{2}\|s^k-x^k\|^2
\end{eqnarray*}
Since $f$ satisfies Assumption \ref{as:conv} we have $\<\nabla f(x^k), x^* - x^k> \leq -(f(x^k) - f(x^*))$. Thus, 
\begin{eqnarray*}
    f(x^{k+1}) - f(x^*) &\leq& f(x^k) - f(x^*)  - \eta_k(f(x^k) - f(x^*)) + \frac{\alpha}{L}\|\nabla f(x^k) - g^k\|^2 \\
    &+& \frac{L\eta_k^2}{\alpha}\|s^k - x^*\|^2 + \frac{L\eta_k^2}{2}\|s^k - x^k\|^2 \\
    &=&(1-\eta_k)(f(x^k) - f(x^*)) + \frac{\alpha}{L}\|\nabla f(x^k) - g^k\|^2 \\
    &+& \frac{L\eta_k^2}{\alpha}\|s^k - x^*\|^2 + \frac{L\eta_k^2}{2}\|s^k - x^k\|^2.
\end{eqnarray*}
\EndProof{Taking the full mathematical expectation finishes the proof.} \\
Now we are ready to proof our main result. We start with the convex case, where $f$ satisfies \ref{as:conv}. For readers convenience, we restate the theorems below.
\begin{theorem}
    \textup{(Theorem \ref{th:convex})}\textbf{.} Let the Assumptions \ref{as:lip}, \ref{as:conv} and \ref{as:key} be satisfied. Then there exist $\eta_k \leq \min(\rho_1, \rho_2)$ for Algorithm \eqref{eq:stoch_fw} and constants $M_1$, $M_2$, $\alpha$ such that:\\
    \begin{eqnarray*}
        r_{K+1} = \mathcal{O}\Biggl(r_0 \exp\left(-\frac{K \min(\rho_1,\rho_2)}{4}\right) + \frac{LD^2}{K + \frac{1}{\min(\rho_1,\rho_2)}} + D^2\sqrt{\frac{B\rho_2 + AE}{\rho_1\rho_2\left(K + \frac{1}{\min(\rho_1,\rho_2)}\right)^2} + \frac{K}{K + \frac{1}{\min(\rho_1,\rho_2)}}\frac{C\rho_2}{\rho_1\rho_2D^2}}\Biggr),
    \end{eqnarray*}\\
where $r_k = \EE_k[f(x^k) - f(x_*) + M_1\|g^k - \nabla f(x^k)\|^2 + M_2\sigma_k^2]$.
\end{theorem}
\vspace{-0.3cm}
\textbf{Proof:}\\
Since Assumption \ref{as:key} holds, we can derive the following inequality:
\begin{eqnarray*}
    r_{k+1}&=&\EE[f(x^{k+1}) - f(x^*) + M_1\|\nabla f(x^{k+1}) - g^{k+1}\|^2 + M_2\sigma_{k+1}^2] \stackrel{(\ref{lem:conv})}{\leq} (1 - \eta_k)\EE[f(x^k) - f(x^*)] \\
    &+& \frac{\alpha}{L}\EE[\|\nabla f(x^k) - g^k\|^2] + \left(\frac{L}{\alpha} + \frac{L}{2}\right)D^2\eta_k^2 + M_1(1-\rho_1)\|g^{k} - \nabla f(x^{k})\|^2 \\
    &+& M_1A\sigma_k^2 + M_1\eta_k^2BD^2 + M_1C + M_2(1-\rho_2)\sigma_k^2 + M_2\eta_k^2ED^2 \\
    &=&(1 - \eta_k)\EE[f(x^k) - f(x^*)] + \left(\frac{\alpha}{L} + M_1(1-\rho_1)\right)\EE\left[\|\nabla f(x^k) - g^k\|^2\right] + \left(M_1A + M_2(1 - \rho_2)\right)\sigma_k^2\\
    &+& \left(\frac{L}{\alpha} + \frac{L}{2} + M_1B + M_2E\right)\eta_k^2D^2 + M_1C.
\end{eqnarray*}
With constants $M_1 = \frac{2\alpha}{\rho_1L}, ~ M_2 = \frac{2M_1A}{\rho_2}$ we get\\
\begin{equation*}
    \begin{split}
        \;\;\;\;\;&\;\EE[f(x^{k+1}) - f(x^*)] + M_1\EE[\|\nabla f(x^{k+1}) - g^{k+1}\|^2] + M_2\sigma_{k+1}^2 ~\leq~ (1 - \eta_k)\EE[f(x^k) - f(x^*)]\\ &+ \left(1-\frac{\rho_1}{2}\right)M_1\EE[\|\nabla f(x^k) - g^k\|^2] + \left(1 - \frac{\rho_2}{2}\right)M_2\sigma_k^2 + \left(\frac{L}{2} + \frac{L}{\alpha} + \frac{2B}{\rho_1L}\alpha + \frac{4AE}{\rho_1\rho_2L}\alpha \right)\eta_k^2D^2 + \frac{2C}{\rho_1L}\alpha. 
    \end{split}
\end{equation*}\\
Using Lemma \ref{lem:b4} with $d = \frac{2}{\min(\rho_1, \rho_2)}$ and $r_k = \EE_k[f(x^k) - f(x_*) + M_1\|\nabla f(x^k) - g^k\|^2 + M_2\sigma_k^2]$ we obtain

\begin{equation*}
    r_{K+1} = \mathcal{O}\left(r_0\exp\left(-\frac{K}{2d}\right) + \frac{LD^2}{2(K+d)} + \frac{\frac{L}{\alpha} + \alpha \frac{2B\rho_2 + 4AE}{\rho_1\rho_2L}}{K+d}D^2 + K\alpha\frac{2C\rho_2}{\rho_1\rho_2L}\right) .  
\end{equation*}

This estimation holds for any $\alpha > 0$, and thus to obtain optimal estimation on the K-th iteration we shall minimize this to $\alpha$. It is easy to see, that minimum of $x\alpha + \frac{y}{\alpha}$ is located at $\alpha = \sqrt{\frac{y}{x}}$ and equals to $2\sqrt{xy}$. Therefore taking optimal $\alpha$ as
\begin{equation*}
    \alpha = \sqrt{\frac{\frac{LD^2}{K+d}}{\frac{2B\rho_2 + 4AE}{\rho_1\rho_2L(K+d)}D^2 + K\frac{2C\rho_2}{\rho_1\rho_2L}}},
\end{equation*}
we get 
$$r_{K+1} = \mathcal{O}\left(r_0\exp\left(-\frac{K}{2d}\right) + \frac{LD^2}{K+d} + D^2\sqrt{\frac{2B\rho_2 + 4AE}{\rho_1\rho_2(K+d)^2} + \frac{K}{K+d}\frac{2C\rho_2}{\rho_1\rho_2D^2}}\right).$$
\EndProof{It completes the proof.}

\begin{theorem}
    \textup{(Theorem \ref{th:nonconvex})}\textbf{.} Let the Assumptions \ref{as:lip} and \ref{as:key} be satisfied. Then, there exist $\eta_k$ for Algorithm \eqref{eq:stoch_fw}, that
\begin{equation}
    \label{non-convex}
    \EE\left[\min\limits_{0\leq k \leq K-1} \textbf{gap}(x^k)\right] = \mathcal{O}\left(\frac{r_0}{\sqrt{K}} + \frac{D^2}{\sqrt{K}}\left[L + \sqrt{\frac{B\rho_2 + AE}{\rho_1\rho_2} + K\frac{C\rho_2}{D^2\rho_1\rho_2}}\right]\right).
\end{equation}
\end{theorem}
\textbf{Proof:}\\
With the Assumption \ref{as:lip} we can use (\ref{Nesterov}):
\begin{eqnarray*}
    f(x^{k+1}) &\leq& f(x^k) + \<\nabla f(x^k), x^{k+1} - x^k> + \frac{L}{2}\|x^{k+1} - x^k\|^2.
\end{eqnarray*}
With update of $x^{k+1}$, according to \ref{eq:stoch_fw}:
\begin{eqnarray*}
    f(x^{k+1}) - f(x^*)&\leq& f(x^k) - f(x^*) + \eta_k\<\nabla f(x^k), s^k - x^k> + \frac{L\eta_k^2}{2}\|s^k - x^k\|^2 \\ 
    &=&f(x^k) - f(x^*) + \eta_k\<g^k, s^k - x^k> + \eta_k\<\nabla f(x^k) - g^k, s^k - x^k> + \frac{L\eta_k^2}{2}\|s^k-x^k\|^2.
\end{eqnarray*}
The optimal choice of $s^k$ in \ref{eq:stoch_fw} gives $\<g^k, s^k - x^k>\leq\<g^k, x - x^k>$ for all $x \in \X$. Then, 
\begin{eqnarray*}
    f(x^{+1}) - f(x^*)&\leq& f(x^k) - f(x^*) + \eta_k\<g^k, x - x^k> + \eta_k\<\nabla f(x^k) - g^k, s^k - x^k> + \frac{L\eta_k^2}{2}\|s^k-x^k\|^2\\
    &=&f(x^k) - f(x^*) + \eta_k\<\nabla f(x^k), x - x^k> + \eta_k\<g^k - \nabla f(x^k), x - x^k> \\
    &+&\eta_k\<\nabla f(x^k) - g^k, s^k - x^k> + \frac{L\eta_k^2}{2}\|s^k-x^k\|^2\\
    &=&f(x^k) - f(x^*) + \eta_k\<\nabla f(x^k), x - x^k> + \eta_k\<\nabla f(x^k) - g^k, s^k - x> + \frac{L\eta_k^2}{2}\|s^k - x^k\|^2.
\end{eqnarray*}
Using Young's inequality \eqref{eq:young2} here we state for any positive $\alpha$:
\begin{eqnarray*}
    f(x^{k+1}) - f(x^*)&\leq& f(x^k)- f(x^*) + \eta_k\<\nabla f(x^k), x - x^k> + \frac{\alpha}{L}\|\nabla f(x^k) - g^k\|^2 \\
    &+& \frac{L\eta_k^2}{\alpha}\|s^k - x\|^2 + \frac{L\eta_k^2}{2}\|s^k-x^k\|^2.
\end{eqnarray*}
After small rearrangements one can get:
\begin{eqnarray*}
    \eta_k\<\nabla f(x^k), x^k - x> &\leq& f(x^k) - f(x^*) - \left(f(x^{k+1}) - f(x^*)\right) + \frac{\alpha}{L}\|\nabla f(x^k) - g^k\|^2 \\
    &+& \frac{L\eta_k^2}{\alpha}\|s^k - x\|^2 + \frac{L\eta_k^2}{2}\|s^k - x^k\|^2.
\end{eqnarray*}
Maximizing over $\X$, taking the full mathematical expectation and bounding the distances by diameter, we get:\\
\begin{eqnarray*}
    \eta_k \EE\Bigl[\max\limits_{x \in \X} \<\nabla f(x^k), x^k - x>\Bigr] &\leq& \EE\Bigl[f(x^k) - f(x^*)\Bigr] - \EE\Bigl[f(x^{k+1}) - f(x^*))\Bigr] + \frac{\alpha}{L}\EE\Bigl[\|\nabla f(x^k) - g^k\|^2\Bigr] \\
    &+& \eta_k^2\frac{LD^2}{\alpha} + \eta_k^2\frac{LD^2}{2}.
\end{eqnarray*}
After multiplying \ref{as:key} by the positive constants $M_1,\,M_2$ (which we will define below) and summarizing with previous inequality we have:
\begin{eqnarray*}
    \eta_k \EE\Bigl[\max\limits_{x \in \X} \<\nabla f(x^k), x^k - x>\Bigr] &\leq& \EE\Bigl[f(x^k) - f(x^*) +  \left(1 - \rho_1 + \frac{\alpha}{M_1L}\right)M_1\|\nabla f(x^k) - g^k\|^2 \\
    &+& \left(1 - \rho_2 + \frac{M_1A}{M_2}\right)M_2\sigma_k^2\Bigr] \\ 
    &-& \EE\Bigl[f(x^{k+1}) - f(x^*)) +M_1\|\nabla f(x^{k+1}) - g^k\|^2 + M_2\sigma_{k+1}^2\Bigr] \\
    &+& D^2\eta_k^2\left(\frac{L}{2} + \frac{L}{\alpha} + M_1B + M_2E\right) + M_1C.
\end{eqnarray*}
With $M_1 = \frac{\alpha}{L\rho_1}, M_2 = \frac{M_1A}{\rho_2}$ and summarizing over all k from 0 to $K-1$ we have:
\begin{eqnarray*}
    \sum\limits_{k = 0}^{K-1} \eta_k \EE\Bigl[\max\limits_{x \in \X} \<\nabla f(x^k), x^k - x>\Bigr] &\leq& f(x^0) - f(x^*) + \|\nabla f(x^0) - g^0\|^2 + \|\sigma_0\|^2 \\
    &+& D^2\left(\frac{L}{2} + \frac{L}{\alpha} + \alpha\left(\frac{B}{\rho_1L} + \frac{AE}{\rho_1\rho_2L}\right)\right) \sum\limits_{k = 0}^{K-1} \eta_k^2 + KM_1C.
\end{eqnarray*}
Assuming $\min_{k \leq K}\eta_k = \eta_{min}$, we gain:
\begin{eqnarray*}
    \sum\limits_{k = 0}^{K-1} \eta_k \EE\Bigl[\max\limits_{x \in \X} \<\nabla f(x^k), x^k - x>\Bigr] &\leq& f(x^0) - f(x^*) + \|\nabla f(x^0) - g^0\|^2 + \|\sigma_0\|^2 \\
    &+& D^2\left(\frac{L}{2} + \frac{L}{\alpha} + \alpha\left(\frac{B}{\rho_1L} + \frac{AE}{\rho_1\rho_2L} + \frac{C}{\rho_1L\eta_{min}^2D^2}\right)\right) \sum\limits_{k = 0}^{K-1} \eta_k^2.
\end{eqnarray*}
With $\alpha = \sqrt{\frac{LD^2}{\frac{1}{\rho_1L}(BD^2 + \frac{C}{\eta^2_{min}}) + \frac{A}{\rho_1\rho_2L}ED^2}}$ one can obtain:
\begin{eqnarray*}
    \sum\limits_{k = 0}^{K-1} \eta_k \EE\Bigl[\max\limits_{x \in \X} \<\nabla f(x^k), x^k - x>\Bigr] &\leq& f(x^0) - f(x^*) + \|\nabla f(x^0) - g^0\|^2 ~~~~~~~~~~~~~~~~\\
    &+& D^2\left(\frac{L}{2} + 2\sqrt{\frac{B\rho_2 + AE}{\rho_1\rho_2} +\frac{1}{\eta_{min}^2}\frac{C\rho_2}{D^2\rho_1\rho_2}}\right)\sum\limits_{k = 0}^{K-1} \eta_k^2.
\end{eqnarray*}
If we take $\eta_k = \frac{1}{\sqrt{K}}$, then $\eta_{min} = \frac{1}{\sqrt{K}}$. Divide both sides by $\sqrt{K}$, then:\\
\begin{eqnarray*}
    \EE\Bigl[\frac{1}{K}\max\limits_{x \in \X} \<\nabla f(x^k), x^k - x>\Bigr] \;~\leq~\; \frac{f(x^0) - f(x^*) + \|\nabla f(x^0) - g^0\|^2 + \|\sigma_0\|^2}{\sqrt{K}} +~~~~~~~~~~~~~~~~~~~~~~~~~~~~~~~~~&\\+ \frac{D^2}{\sqrt{K}}\Bigl(\frac{L}{2} + 2\sqrt{\frac{B\rho_2 + AE}{\rho_1\rho_2} +K\frac{C\rho_2}{D^2\rho_1\rho_2}}\Bigr).
\end{eqnarray*}
\EndProof{Finally, we obtain the needed estimation.}
\subsection{Stochastic methods}
In this section, we provide the detailed convergence rates and proofs for some specific methods (see Section \ref{sec:stoch_main}) in the finite-sum case \eqref{eq:finsum} of the constrained optimization problem \eqref{eq:main}. \\
In the following, we prove that some specific methods, i.e., \texttt{L-SVRG}, \texttt{SARAH} and \texttt{SAGA} satisfy our unified Assumption \ref{as:key} and thus can be captured by our unified analysis. Then, we plug their corresponding parameters (i.e., specific values for $A, B, C, E, \sigma_k^2, \rho_1, \rho_2$) into our unified Theorems \ref{th:convex} and \ref{th:nonconvex} to obtain the detailed convergence rates for these methods.
\subsubsection{L-SVRG Frank-Wolfe}\label{sec:lsvrg}
We first restate our Lemma \ref{lem:lsvrg_conv} for \texttt{L-SVRG FW} method (Algorithm \ref{alg:lsvrg}) and provide its proof. Then we plug its corresponding parameters (i.e., specific values for $A, B, C, E, \sigma_k^2, \rho_1, \rho_2$) into our unified Theorems \ref{th:convex} and \ref{th:nonconvex} to obtain the detailed convergence rate.
\begin{algorithm}[H]
            \caption{\texttt{L-SVRG Frank-Wolfe}}\label{alg:lsvrg}
            \textbf{Input:} initial $x^0$, $w^0 = x^0$, $g^0 = \nabla f(x^0)$ step sizes $\{\eta_k\}_{k\geq0}$, batch size $b$, probability $p\in(0,1]$
            \begin{algorithmic}
            \For {$k = 0, 1, 2,\dots K-1$} 
            \State Compute $s^k = \arg\underset{s\in \mathcal{X}}{\min}\langle s,g^k\rangle$
            \State Update $x^{k+1} = (1 - \eta_k)x^k + \eta_ks^k$
            \State Update $w^{k+1} = 
            \displaystyle\begin{cases}
                x^k, &\text{ with probability } p\\ \vspace{-0.9em}\\
                w^k, &\text{ with probability }1-p
            \end{cases}
            $
            \State Generate batch $S_k$ with size $b$
            \State $g^{k+1} = \tfrac{1}{b}\underset{i\in S_k}{\sum}\left[\nabla f_i(x^{k+1}) - \nabla f_i(w^{k+1})\right] + \nabla f(w^{k+1})$
            
            \EndFor
            \end{algorithmic}
\end{algorithm}
\begin{lemma}
\label{lem:ap:lsvrg}
\textup{(Lemma \ref{lem:lsvrg_conv})}\textbf{.} Under Assumption \ref{as:lip:local} Algorithm \ref{alg:lsvrg} satisfies Assumption \ref{as:key} with
$$\rho_1 = 1, ~ A = \frac{\widetilde{L}^2}{b}\left(1-\frac{p}{2}\right), ~ B = \frac{8\widetilde{L}^2}{pb}, ~C = 0,$$
$$ \sigma_k^2 = \|x^{k}-w^{k}\|^2,~ \rho_2 = \frac{p}{2}, ~ E =\frac{8}{p}.$$
\end{lemma}
\textbf{Proof:}\\
According to Lemma 3 from \citep{li2020unified}  we get an estimation: 
$$
\mathbb{E}_k[\|g^k\|^2]\leq\frac{\widetilde{L}^2}{b}\|x^k - w^k\|^2 + \|\nabla f(x^k)\|^2.
$$
Since $g^k$ is unbiased gradient estimator previous inequality turns to:
\begin{equation*}
\mathbb{E}_k[\|\nabla f(x^k) - g^k\|^2] \leq \frac{\widetilde{L}^2}{b}\|x^k - w^k\|^2.
\end{equation*}
Considering Algorithm \ref{alg:lsvrg} we have:
\begin{eqnarray}
    \EE_k[\|x^{k} - w^{k}\|^2] 
    &=&
    p\EE_k\left[\|x^{k} - x^{k-1}\|^2] + (1-p)\EE_k[\|x^{k} - w^{k-1}\|^2\right]
    \notag\\&=&
    p\eta_{k-1}^2\EE_k \left[\|s^{k-1} - x^{k-1}\|^2] + (1-p)\EE_k[\|x^{k-1} + \eta_k(s^{k-1} - x^{k-1}) - w^{k-1}\|^2\right]
    \notag\\&=&
    \eta_{k-1}^2\EE_k[\|s^{k-1} - x^{k-1}\|^2] + (1-p)\|x^{k-1}-w^{k-1}\|^2 
    \\
    &&+ 2\eta_k(1-p)\EE_k[\langle x^{k-1} - w^{k-1}, s^{k-1} - x^{k-1}\rangle]
    \notag\\
    &=&\eta_{k-1}^2\EE_k[\|s^{k-1} - x^{k-1}\|^2] + (1-p)\|x^{k-1}-w^{k-1}\|^2 
    \\
    &&
    + 2(1-p)\EE_k[\langle x^{k-1} - w^{k-1},\eta_k(s^{k-1} - x^{k-1})\rangle].
    \notag
\end{eqnarray}
According to Young's inequality for any positive $\beta$ there are:
$$\langle x^{k-1} - ^{k-1}, (s^{k-1} - x^{k-1})\eta_k\rangle \leq \beta\|x^{k-1} - w^{k-1}\|^2 + \frac{1}{\beta}\eta_{k-1}^2\|s^{k-1} - x^{k-1}\|^2.$$
Hence,
\begin{eqnarray}\label{for_comp}
    \mathbb{E}_k[\|x^{k} - w^{k}\|^2] 
    &\leq&\eta_{k-1}^2D^2 + (1-p)\|x^{k-1} - w^{k-1}\|^2 + 2(1-p)\beta\|x^{k-1} - w^{k-1}\|^2 \\
    &&+ 2\frac{1-p}{\beta}\|\eta_{k-1}(s^{k-1} - x^{k-1})\|^2
    \notag\\
    &\leq& \left(1 + \frac{2(1-p)}{\beta}\right)\eta_{k-1}^2D^2 + (1-p)(1 + 2\beta)\|x^{k-1} - w^{k-1}\|^2.
\end{eqnarray}
Finally, choose $\beta=\frac{p}{4}$. It leads to:
$$(1-p)(1 +2\beta) \leq \left(1 - \frac{p}{2}\right).$$
Then,
\begin{equation*}
    \EE_k[\|x^{k} - w^{k}\|^2] \leq \frac{8}{p}\eta_{k-1}^2D^2 + \left(1-\frac{p}{2}\right)\|x^{k-1} - w^{k-1}\|^2,
\end{equation*}
and
$$\mathbb{E}_k[\|\nabla f(x^{k}) - g^{k}\|^2]\leq\frac{8\widetilde{L}^2}{pb}\eta_{k-1}^2D^2 + \frac{\widetilde{L}^2}{b}\left(1-\frac{p}{2}\right)\|x^{k-1} - w^{k-1}\|^2.$$
\EndProof{This finishes the proof.}

\begin{corollary}
\textup{(Corollary \ref{crl:lsvrg})}\textbf{.} Suppose that Assumption \ref{as:lip} holds. For Algorithm \ref{alg:lsvrg} in the convex and non-convex cases the following convergences take place:
$$\mathbb{E}[r_{K+1}] =\mathcal{O}\left((f(x_0) - f(x_*))\exp\left(-\frac{Kp}{8}\right) + \frac{LD^2}{K+\frac{1}{p}}\left[1 + \frac{\widetilde{L}}{L}\frac{1}{p\sqrt{b}}\right]\right).$$
$$\EE\Bigl[\min\limits_{0 \leq k \leq K-1} \textbf{gap}(x^k)\Bigr] = \mathcal{O}\left(\frac{f(x_0) - f(x_*)}{\sqrt{K}} + \frac{LD^2}{\sqrt{K}}\left[1 + \frac{\widetilde{L}}{L}\frac{1}{p\sqrt{b}}\right]\right).$$
\end{corollary}
\textbf{Proof:} \\
\EndProof{It suffices to plug parameters from Lemma \ref{lem:ap:lsvrg} into Theorems \ref{th:convex} and \ref{th:nonconvex}}

We proved results for \texttt{L-SVRG FW} depending on the parameters $p$ and $b$ to be tuned. Let us find the optimal choices of $p$ and $b$ for \texttt{L-SVRG FW}. To find the optimal choice of $p$, one can note that on average we call the stochastic gradients $(pn + 2b)$ times. In more details, at each iteration we compute a batch size of $b$ in two points $x^{k}$ and $w^k$ and with probability $p$ we call the full gradient in the new point $w^k$. From Corollary \ref{crl:lsvrg}, we know the estimate on the number of iterations of \texttt{L-SVRG FW}, then we can get an estimate on the number of the stochastic gradient calls  by multiplying this result by $(pn + 2b)$. Then the new estimate can be optimized first by $p$ (in the convex case, we need to minimize $(1 + \tfrac{1}{p\sqrt{b}})(pn + 2b)$) and obtain that the optimal $p \sim b^{1/4} / n^{1/2}$. Then with already optimized $p$, the estimate on the number of the stochastic gradient calls can additionally be optimized by $b$ (actually we need to minimize $b^{1/4} n^{1/2} + b + n / b^{1/2}$) and find the optimal $b \sim n^{2/3}$. The final result for \texttt{L-SVRG FW} is presented in Table \ref{tab:comparison_stoch}. 

\subsubsection{SARAH Frank-Wolfe}\label{sec:SARAH}
We first restate our Lemma \ref{lem:sarah} for \texttt{SARAH FW} method (Algorithm \ref{alg:sarah}) and provide its proof. Then we plug its corresponding parameters (i.e., specific values for $A, B, C, E, \sigma_k^2, \rho_1, \rho_2$) into our unified Theorems \ref{th:convex} and \ref{th:nonconvex} to obtain the detailed convergence rate.
\begin{algorithm}[H]
            \caption{\texttt{SARAH Frank-Wolfe}}\label{alg:sarah}
             \textbf{Input:} initial $x^0$, $g^{0} = \nabla f(x^{0})$ step sizes $\{\eta_k\}_{k\geq0}$, batch size $b$, probability $p\in(0,1]$
            \begin{algorithmic}
            \For {$k = 0, 1, 2,\dots K-1$}
            \State Compute $s^k = \arg\underset{s\in \mathcal{X}}{\min}\langle s, g^k\rangle$
            \State Update $x^{k+1} = (1-\eta_k)x^k + \eta_ks^k$
            \State Generate batch $S_k$ with size $b$
            \State Update $g^{k+1} =
            \displaystyle\begin{cases}
            \nabla f(x^{k+1}), &\text{ with probability }p \\ \vspace{-0.9em}\\
            g^{k1} + \frac{1}{b}\underset{i\in S_k}{\sum}\left[\nabla f_i(x^{k+1}) - \nabla f_i(x^{k})\right], &\text{ with probability }1-p
            \end{cases}
            $
            \EndFor
            \end{algorithmic}
\end{algorithm}
\begin{lemma}\label{lem:ap:sarah}
\textup{(Lemma \ref{lem:sarah})}\textbf{.} Under Assumption \ref{as:lip:local} Algorithm \ref{alg:sarah} satisfies Assumption \ref{as:key} with:
$$\rho_1 = p, ~ A = 0, ~ B = \frac{1-p}{b}\widetilde{L}^2, ~C = 0,$$
$$ \sigma_k = 0,~ \rho_2 = 1, ~ E = 0.$$
\end{lemma}
\textbf{Proof:}\\
Using Lemma 3 from \citep{pmlr-v139-li21a} we can obtain: 
\begin{eqnarray*}
\EE_k\left[\|\nabla f(x^{k}) - g^{k}\|^2\right]&\leq& (1 - p)\|\nabla f(x^{k-1}) - g^{k-1}\|^2 + \frac{1-p}{b}\widetilde{L}^2\|x^{k}-x^{k-1}\|^2 \\
& \leq &(1 - p)\|\nabla f(x^{k-1}) - g^{k-1}\|^2 + \frac{1-p}{b}\widetilde{L}^2\eta_{k-1}^2D^2.
\end{eqnarray*}\EndProof{}
\begin{corollary}
\textup{(Corollary \ref{crl:sarah})}\textbf{.} For Algorithm \ref{alg:sarah} in the convex and non-convex cases the following convergences take place:
$$\EE\left[r_{K+1}\right] = \mathcal{O}\left(\Big(f(x^0) - f(x^*)\Big)\exp\left(-\frac{Kp}{4}\right) + \frac{LD^2}{K+\frac{1}{p}}\left[1 + \frac{\widetilde{L}}{L}\sqrt{\frac{1-p}{pb}}\right]\right). $$
$$\EE\left[\min\limits_{0 \leq k \leq K-1} \textbf{gap}(x^k)\right] = \mathcal{O}\left(\frac{f(x^0) - f(x^*)}{\sqrt{K}} + \frac{LD^2}{\sqrt{K}}\left[1 + \frac{\widetilde{L}}{L}\sqrt{\frac{1-p}{pb}}\right]\right).$$
\end{corollary}
\textbf{Proof:} \\
\EndProof{It suffices to plug parameters from Lemma \ref{lem:ap:sarah} into Theorems \ref{th:convex} and \ref{th:nonconvex}.}

The choices of $p$ and $b$ for \texttt{SARAH FW} are presented in the original paper \citep{Beznosikov2023SarahFM}. The final result for \texttt{SARAH FW} is presented in Table \ref{tab:comparison_stoch}.

\subsubsection{SAGA Frank-Wolfe}\label{sec:SAGA}
We first restate our Lemma \ref{lem:saga} for \texttt{SAGA FW} method (Algorithm \ref{alg:saga}) and provide its proof. Then we plug its corresponding parameters (i.e., specific values for $A, B, C, E, \sigma_k^2, \rho_1, \rho_2$) into our unified Theorems \ref{th:convex} and \ref{th:nonconvex} to obtain the detailed convergence rate.
\begin{algorithm}[H]
            \caption{\texttt{SAGA Frank-Wolfe}}\label{alg:saga}
             \textbf{Input:} initial $x^{0}$, $\forall i \in [n]~ y_i^{0} = \nabla f_i(x^{0})$, $g^0 = \nabla f(x^0)$, step sizes $\{\eta_k\}_{k\geq0}$, batch size $b$
            \begin{algorithmic}
            \For {$k = 0, 1, 2,\dots K-1$}
            \State Compute $s^k = \arg\underset{s\in \mathcal{X}}{\min}\langle s, g^k\rangle$
            \State Update $x^{k+1} = (1-\eta_k)x^k + \eta_ks^k$
            \State Generate batch $S_k$ with size $b$
            \State Update $y_i^{k+1} =
            \displaystyle\begin{cases}
            \nabla f_i(x^{k}), &\text{for } i\in S_k \\ \vspace{-0.9em}\\
            y_i^{k}, &\text{for }i \notin S_k
            \end{cases}
            $
            \State Update $g^{k+1} = \frac{1}{b}\underset{i\in S_k}{\sum}\left[\nabla f_i(x^{k+1}) - y_i^{k+1}\right] + \frac{1}{n}\sum\limits_{j=1}^n  y_j^{k+1}$
            \EndFor
            \end{algorithmic}
\end{algorithm}
\begin{lemma}
\textup{(Lemma \ref{lem:saga})}\textbf{.} Under Assumption \ref{as:lip:local} Algorithm \ref{alg:saga} satisfies Assumption \ref{as:key} with:
$$\rho_1 = 1,~A = \frac{1}{b}\left(1 + \frac{b}{2n}\right), ~ B = \frac{\widetilde{L}^2}{b}\left(1+\frac{2n}{b}\right),~C = 0,$$
$$  \sigma_k^2 = \frac{1}{n}\sum\limits_{j = 1}^n\|\nabla f_j(x^k) - y_j^{k+1}\|^2 ,~\rho_2 = \frac{b}{2n},~ E =\frac{2n}{b}\widetilde{L}^2.$$
\end{lemma}
\textbf{Proof:}\\
We bound the difference between estimator and exact gradient:
\begin{eqnarray*}
    \EE_k\left[\|g^k - \nabla f(x^k)\|^2\right] &=& \EE_k\left[\Biggl\|\frac{1}{b}\sum\limits_{i \in S_k}\left[\nabla f_i(x^{k}) - y_i^{k}\right] + \frac{1}{n}\sum\limits_{j = 1}^n y_j^{k} - \nabla f(x^k) \Biggr\|^2\right]
    \\&=& \EE_k\Biggl[\Biggl\|\frac{1}{b}\left(\sum\limits_{i \in S_k}\left[\nabla f_i(x^{k}) - y_i^{k}\right] - \left(\frac{1}{n}\sum\limits_{j = 1}^n\left[\nabla f_j(x^k)-y_j^{k}\right]\right)\right) \Biggr\|^2\Biggr]
    \\&\stackrel{(\ref{lem:2moment_uni_batch})}{\leq}& \frac{1}{bn}\sum\limits_{j=1}^n\left\|\nabla f_j(x^k) - y_j^k - \left(\frac{1}{n}\sum\limits_{i=1}^n\left[\nabla f_i(x^k)-y_i^k\right]\right)\right\|^2
    \\&\leq&\frac{1}{bn}\sum\limits_{j=1}^n\left\|\nabla f_j(x^k) - y_j^k\right\|^2
    \\&\leq&\frac{1}{bn}\left(1+\alpha\right)\sum\limits_{j=1}^n\|\nabla f_j(x^k) - \nabla f_j(x^{k-1})\|^2 + \frac{1}{bn}\left(1+\frac{1}{\alpha}\right)\sum\limits_{j=1}^n\|\nabla f_j(x^{k-1}) - y_j^k\|^2
    \\&\leq& \frac{\widetilde{L}^2}{b}\left(1+\alpha\right)\eta_{k-1}^2D^2 + \frac{1}{b}\left(1+\frac{1}{\alpha}\right)\sigma_{k-1}^2
\end{eqnarray*}
for $\forall \alpha>0$ (in particular, we can put $\alpha = \frac{2n}{b}$ to obtain the needed estimates). The second inequality holds, since $\frac{1}{n}\sum\limits_{i=1}^n$ can be described, as an expected value. And $\EE\|x - \EE x\|^2 \leq \EE\|x\|^2.$ Then we need to bound the second term:
\begin{eqnarray*}
    \EE_k[\sigma_{k}^2]&=&\EE_k\left[\frac{1}{n}\sum\limits_{j = 1}^n\|\nabla f_j(x^{k}) - y_j^{k+1}\|^2 \right] = \left(1 - \frac{b}{n}\right)\frac{1}{n}\sum\limits_{j = 1}^n\|\nabla f_j(x^{k}) - y_j^{k}\|^2 
    \\ & =& \left(1 - \frac{b}{n}\right)\frac{1}{n}\sum\limits_{j = 1}^n\|\nabla f_j(x^{k}) - \nabla f_j(x^{k-1}) + \nabla f_j(x^{k-1})-  y_j^{k-1}\|^2
    \\ &\leq& \left(1 - \frac{b}{n}\right)(1 + \beta)\frac{1}{n}\sum\limits_{j = 1}^n\|\nabla f_j(x^{k-1}) - y_j^{k-1}\|^2 + \left(1 - \frac{b}{n}\right)\left(1 + \frac{1}{\beta}\right)\widetilde{L}^2\|x^{k} - x^{k-1}\|^2.
\end{eqnarray*}
With $\beta = \frac{b}{2n}$ we have:
\begin{eqnarray*}
    \EE_k[\sigma_{k}^2] \leq \left(1 - \frac{b}{2n}\right)\sigma_{k-1}^2 + \frac{2n}{b}\widetilde{L}^2\eta_{k-1}^2D^2.
\end{eqnarray*}
\EndProof{This finishes the proof}
\begin{corollary}
\textup{(Corollary \ref{crl:saga})}\textbf{.} For Algorithm \ref{alg:saga} in the convex and non-convex cases the following convergences take place:
$$
\mathbb{E}[r_{K+1}]= \mathcal{O}\left(\Big(f(x^0) - f(x^*)\Big)\exp\left(-\frac{Kb}{8n}\right) + \frac{LD^2}{K+\frac{2n}{b}}\left[1 + \frac{\widetilde{L}}{L}\frac{n}{b\sqrt{b}}\right]\right).
$$
$$\EE\Bigl[\min\limits_{0 \leq k \leq K-1} \textbf{gap}(x^k)\Bigr] = \mathcal{O}\left(\frac{f(x^0) - f(x^*)}{\sqrt{K}} + \frac{LD^2}{\sqrt{K}}\left[1 + \frac{\widetilde{L}}{L}\frac{n}{b\sqrt{b}}\right]\right).$$
\end{corollary}
\textbf{Proof:} \\
\EndProof{It suffices to plug parameters from Lemma \ref{lem:ap:sarah} into Theorems \ref{th:convex} and \ref{th:nonconvex}.}

To find the optimal choice of $b$, one can note that we call the stochastic gradients $2b$ times. In more details, at each iteration we compute a batch size of $b$ in two points $x^{k}$ and $w^k$. From Corollary \ref{crl:saga}, we know the estimate on the number of iterations of \texttt{SAGA FW}, then we can get an estimate on the number of the stochastic gradient calls  by multiplying this result by $2b$. We need to minimize $b(1+nb^{-3/2})$ and obtain that the optimal $b \sim n^{2/3}$. One can notice that Algorithm \ref{alg:saga} is required to store $n$ extra vectors $\lbrace y_i\rbrace$ requiring $\mathcal{O}(nd)$ extra memory. The final result for \texttt{SAGA FW} is presented in Table \ref{tab:comparison_stoch}.
\newpage
\subsection{Coordinate methods}
In this section we provide the detailed convergence rates and proofs for specific methods (see Section \ref{sec:coord_main}) of the constrained optimization problem \eqref{eq:main}. These methods use partial derivatives with respect to coordinates instead of taking gradients of terms of finite sums.\\
In the following, we prove that some specific methods, i.e., \texttt{SEGA} and \texttt{JAGUAR} (new proposed method)  satisfy our unified Assumption \ref{as:key} and thus can be captured by our unified analysis. Then, we plug their corresponding parameters (i.e., specific values for $A, B, C, E, \sigma_k^2, \rho_1, \rho_2$) into our unified Theorems \ref{th:convex} and \ref{th:nonconvex} to obtain the detailed convergence rates.
\subsubsection{SEGA Frank-Wolfe}\label{sec:SEGA}
We first restate our Lemma \ref{lem:sega} for \texttt{SEGA FW} method (Algorithm \ref{alg:sega}) and provide its proof. Then we plug its corresponding parameters (i.e., specific values for $A, B, C, E, \sigma_k^2, \rho_1, \rho_2$) into our unified Theorems \ref{th:convex} and \ref{th:nonconvex} to obtain the detailed convergence rate.
\begin{algorithm}[H]
            \caption{\texttt{SEGA Frank-Wolfe}}\label{alg:sega}
             \textbf{Input:} initial  $x^0$, $h^0= \nabla f(x^0)$, $g^0 = \nabla f(x^0)$ step sizes $\{\eta_k\}_{k\geq0}$
            \begin{algorithmic}
            \For {$k = 0, 1, 2,\dots K-1$}
            \State Compute $s^k = \arg\underset{s\in\mathcal{X}}{\min}\langle s, g^k\rangle$
            \State Update $x^{k+1} = (1-\eta_k)x^k + \eta_ks^k$
            \State Sample $i_k\in[d]$ uniformly at random
            \State Set $h^{k+1} = h^k + e_{i_k}(\nabla_{i_k}f(x^k) - h_{i_k}^k)$
            \State Update $g^{k+1} = de_{i_k}(\nabla_{i_k}f(x^{k+1}) - h_{i_k}^{k+1}) + h^{k+1}$
            \EndFor
            \end{algorithmic}
\end{algorithm}
\begin{lemma}\label{lem:ap:sega}
\textup{(Lemma \ref{lem:sega})}\textbf{.} Under Assumptions \ref{as:lip} Algorithm \ref{alg:sega} satisfies Assumption \ref{as:key} with:
$$\rho_1 = 1,~ A = d, ~ B = d^2L^2,~C = 0,$$
$$ \sigma_k^2 =\|h^{k+1} - \nabla f(x^k)\|^2 , ~\rho_2 = \frac{1}{2d},~ E = 3L^2d.$$
\end{lemma}
\textbf{Proof:}\\
We first bound the difference between estimator and exact gradient:
\begin{eqnarray*}
    \EE_{k}\left[\|g^{k} - \nabla f(x^{k})\|^2\right] &=& \EE_{k}\left[\|d e_{i_k} e_{i_k}^T(\nabla f(x^{k}) - h^{k}) + h^{k} - \nabla f(x^{k})\|^2\right] \\
    &=& \EE_{k}\left[\|(I - d e_{i_k} e_{i_k}^T)(h^{k} - \nabla f(x^{k}))\|^2\right]\\
    &=& \EE_{k}\left[(h^{k} - \nabla f(x^{k}))^T(I - d e_{i_k} e_{i_k}^T)^T(I - d e_{i_k} e_{i_k}^T)(h^{k} - \nabla f(x^{k}))\right] \\
    &=& (h^{k} - \nabla f(x^{k}))^T\EE_{k}\left[I - 2d e_{i_k} e_{i_k}^T + d^2e_{i_k} e_{i_k}^T\right](h^{k} - \nabla f(x^{k})) \\
    &=& (h^{k} - \nabla f(x^{k}))^T\left[I - 2\cdot I + d\cdot I\right](h^{k} - \nabla f(x^{k})) \\
    &=& (d-1)\|h^{k} - \nabla f(x^{k})\|^2 \\
    &\leq& (d-1)(1+\alpha)\|h^k - \nabla f(x^{k-1})\|^2 + (d-1)\left(1+\frac{1}{\alpha}\right)\eta_{k-1}^2L^2D^2.
\end{eqnarray*}
Then,
\begin{eqnarray*}
    \EE_k\left[\|h^{k+1} - \nabla f(x^{k})\|^2\right] &=& \EE_k\left[\|h^k + e_{i_k} e_{i_k}^T(\nabla f(x^k) - h^k) - \nabla f(x^{k})\|^2\right]\\
    &=& \EE_k\left[\|(I- e_{i_k} e_{i_k}^T)(h^k - \nabla f(x^{k}))\|^2\right]\\
    &=&\left(1 - \frac{1}{d}\right)\|h^k - \nabla f(x^k)\|^2 \\
    &\leq& \left(1-\frac{1}{d}\right)(1+\beta)\|h^k - \nabla f(x^{k-1})\|^2 + \left(1-\frac{1}{d}\right)\left(1 + \frac{1}{\beta}\right)\eta_{k-1}^2L^2D^2.
\end{eqnarray*}
If $\beta = \frac{1}{2d}$ then $(1 - \frac{1}{d})(1 + \frac{1}{2d}) \leq 1 - \frac{1}{2d}$ and $(1-\frac{1}{d})(1 + 2d) \leq 2d$, then as $d \geq 1$:
\begin{eqnarray*}
     \EE_k\left[\|h^{k+1} - \nabla f(x^{k})\|^2\right] \leq \left(1-\frac{1}{2d}\right)\|h^k - \nabla f (x^{k-1})\|^2 + 3dL^2\eta_{k-1}^2D^2.
\end{eqnarray*}
\EndProof{Taking $\alpha = \frac{1}{d}$, we obtain the needed constants.}
\begin{corollary}
\textup{(Corollary \ref{crl:sega})}\textbf{.} For Algorithm \ref{alg:sega} in the convex and non-convex cases the following convergences take place:
$$\mathbb{E}[r_{K+1}]=\mathcal{O}\left(\Big(f(x^0) - f(x^*)\Big)\exp\left(-\frac{K}{8d}\right) + \frac{LD^2}{K+d}\left[1 + d\sqrt{d}\right]\right).$$
$$\EE\left[\min\limits_{0 \leq k \leq K-1} \textbf{gap}(x^k)\right] = \mathcal{O}\left(\frac{f(x^0) - f(x^*)}{\sqrt{K}} + \frac{LD^2}{\sqrt{K}}\left[1 + d\sqrt{d}\right]\right).$$
\end{corollary}
\textbf{Proof:}\\
\EndProof{It suffices to plug parameters from Lemma \ref{lem:ap:sega} into Theorems \ref{th:convex} and \ref{th:nonconvex}.}
\subsubsection{JAGUAR}\label{sec:JAGUAR}
We first restate our Lemma \ref{lem:yaguar} for \texttt{JAGUAR} method (Algorithm \ref{alg:aytona}) and provide its proof. Then we plug its corresponding parameters (i.e., specific values for $A, B, C, E, \sigma_k^2, \rho_1, \rho_2$) into our unified Theorems \ref{th:convex} and \ref{th:nonconvex} to obtain the detailed convergence rate.
\begin{algorithm}[H]
            \caption{\texttt{JAGUAR}}\label{alg:aytona}
            \textbf{Input:} initial $x^0$, $g^0 = \nabla f(x^0)$, step sizes $\{\eta_k\}_{k\geq0}$ 
            \begin{algorithmic}
            \For {$k = 0, 1, 2,\dots K-1$}
            \State Compute $s^k = \arg\underset{s\in\mathcal{X}}{\min}\langle s, g^k\rangle$
            \State Update $x^{k+1} = (1-\eta_k)x^k + \eta_ks^k$
            \State Sample $i_{k+1}\in[d]$ uniformly at random
            \State Update $g^{k+1} = e_{i_{k+1}}(\nabla_{i_{k+1}}f(x^{k}) - g_{i_{k+1}}^{k}) + g^{k}$
            \EndFor
            \end{algorithmic}
\end{algorithm}
\begin{lemma}\label{lem:ap:yaguar}
\textup{(Lemma \ref{lem:yaguar})}\textbf{.} Under Assumptions \ref{as:lip} Algorithm \ref{alg:aytona} satisfies Assumption \ref{as:key} with:
$$\rho_1 = \frac{1}{2d}, A = 0, ~ B = 3dL^2,~C = 0,$$
$$ \sigma_k^2 = 0 ,~\rho_2 = 1, E =0.$$
\end{lemma}
\textbf{Proof:}\\
We first bound the difference between estimator and exact gradient:
\begin{eqnarray*}
    \EE_k\left[\|g^{k} - \nabla f(x^{k})\|^2\right] &=& \EE_k\left[\| e_{i_k} e_{i_k}^T(\nabla f(x^{k-1}) - g^{k-1}) + g^{k-1} - \nabla f(x^{k})\|^2\right] \\
    &=& \EE_k\left[\| e_{i_k} e_{i_k}^T(\nabla f(x^{k-1}) - g^{k-1}) + g^{k-1} - \nabla f(x^{k}) + \nabla f(x^{k-1}) - \nabla f(x^{k-1})\|^2\right]\\
    &=& \EE_k\left[\|(I - e_{i_k} e_{i_k}^T)(\nabla f(x^{k-1}) - g^{k-1})  + \nabla f(x^{k-1}) - \nabla f(x^{k})\|^2\right]\\
    &\leq& (1+\beta)\EE_k\left[\|(I - e_{i_k} e_{i_k}^T)(g^{k-1} - \nabla f(x^{k-1}))\|^2\right] + \left(1 + \frac{1}{\beta}\right)\eta_{k-1}^2L^2D^2\\
    &=& (1+\beta)\left(1 - \frac{1}{d}\right)\|g^{k-1} - \nabla f(x^{k-1})\| ^2+ \left(1 + \frac{1}{\beta}\right)\eta_{k-1}^2L^2D^2.
\end{eqnarray*}

If $\beta = \frac{1}{2d}$ then $(1 - \frac{1}{d})(1 + \frac{1}{2d}) \leq 1 - \frac{1}{2d}$ and then as $d \geq 1:$
\begin{center}
    $\EE_k\left[\|g^{k} - \nabla f(x^{k})\|^2\right] \leq \left(1 - \frac{1}{2d}\right)\|g^{k-1}-\nabla f(x^{k-1})\|^2 + 3d\eta_{k-1}^2L^2D^2.$
\end{center}
\EndProof{This finishes the proof.}
\begin{corollary}
\textup{(Corollary \ref{crl:yaguar})}\textbf{.} For Algorithm \ref{alg:aytona} in the convex and non-convex cases the following convergences take place:
$$\mathbb{E}[r_{K+1}]= \mathcal{O}\left(\Big(f(x^0) - f(x^*)\Big)\exp\left(-\frac{K}{8d}\right) + \frac{LD^2}{K+d}\left[1 + d\right]\right).$$
$$\EE\left[\min\limits_{0 \leq k \leq K-1} \textbf{gap}(x^k)\right] = \mathcal{O}\left(\frac{f(x^0) - f(x^*)}{\sqrt{K}} + \frac{LD^2}{\sqrt{K}}\left[1 + d\right]\right).$$
\end{corollary}
\textbf{Proof:}\\
\EndProof{It suffices to plug parameters from Lemma \ref{lem:ap:yaguar} into Theorems \ref{th:convex} and \ref{th:nonconvex}.}

\subsubsection{ZOJA}
We introduce \texttt{ZOJA} (\texttt{Zero-Order JAGUAR}) method and provide proof of its convergence. The essence of this method is that in some setting we cannot compute the direction derivatives, but only approximate through zero-order information:
$$
\nabla_{i_{k+1}}f(x^{k}) \approx \widetilde{\nabla}_{i_k}f(x^k) = \frac{f(x^k + \tau e_{i_{k+1}}) - f(x^k)}{\tau}.
$$
This can be used instead of the real directional derivative in the \texttt{JAGUAR FW} method, but it is worth considering the error that arises due to the approximation. 
We derive corresponding parameters (i.e., specific values for $A, B, C, E, \sigma_k^2, \rho_1, \rho_2$) into our unified Theorems \ref{th:convex} and \ref{th:nonconvex} to obtain the detailed convergence rate.
\begin{algorithm}[H]
            \caption{\texttt{ZOJA} (\texttt{Zero-Order JAGUAR})}\label{alg:zoja}
            \textbf{Input:} initial $x^0$, $g^0 = \sum_{i=1}^d \frac{f(x^0 + \tau e_{i}) - f(x^0)}{\tau}$, step sizes $\{\eta_k\}_{k\geq0}, \tau > 0$ 
            \begin{algorithmic}
            \For {$k = 0, 1, 2,\dots K-1$}
            \State Compute $s^k = \arg\underset{s\in\mathcal{X}}{\min}\langle s, g^k\rangle$
            \State Update $x^{k+1} = (1-\eta_k)x^k + \eta_ks^k$
            \State Sample $i_{k+1}\in[d]$ uniformly at random
            \State Compute $\widetilde{\nabla}_{i_k}f(x^k) = \frac{f(x^k + \tau e_{i_{k+1}}) - f(x^k)}{\tau}$
            \State Update $g^{k+1} = e_{i_{k+1}}(\widetilde{\nabla}_{i_{k+1}}f(x^{k}) - g_{i_{k+1}}^{k}) + g^{k}$
            \EndFor
            \end{algorithmic}
\end{algorithm}
\begin{lemma}\label{lem:ap:zoja}
Under Assumptions \ref{as:lip} Algorithm \ref{alg:zoja} satisfies Assumption \ref{as:key} with:
$$\rho_1 = \frac{1}{4d}, A = 0, ~ B = ,3dL^2~C = \frac{5dL^2\tau^2}{4},$$
$$ \sigma_k^2 = 0 ,~\rho_2 = 1, E =0.$$
\end{lemma}
\textbf{Proof:}\\
We bound the difference between estimator and exact gradient:
\begin{eqnarray*}
    \EE_k\left[\|g^{k} - \nabla f(x^{k})\|^2\right] &=& \EE_k\left[\| e_{i_k}(\widetilde{\nabla}_{i_{k-1}} f(x^{k-1}) - g_{i_{k-1}}^{k-1}) + g^{k-1} - \nabla f(x^{k})\|^2\right] \\
    &=& \EE_k\left[\| e_{i_k}(\widetilde{\nabla}_{i_{k-1}} f(x^{k-1}) - g_{i_{k-1}}^{k-1}) + g^{k-1} - \nabla f(x^{k}) + \nabla f(x^{k-1}) - \nabla f(x^{k-1})\|^2\right] \\
    &\leq&  (1+\beta)\EE_k\left[\| e_{i_k}(\widetilde{\nabla}_{i_{k-1}} f(x^{k-1}) - g_{i_{k-1}}^{k-1}) + g^{k-1} - \nabla f(x^{k-1})\|^2\right]  + \left(1 + \frac{1}{\beta}\right)\eta_{k-1}^2L^2D^2\\
    &=& (1+\beta)\EE_k\left[\|(I - e_{i_k}e_{i_k}^T)(g^{k-1} - \nabla f(x^{k-1})) + e_{i_k}(\widetilde{\nabla}_{i_k}f(x^{k-1}) - \nabla_{i_k} f(x^{k-1}))\|\right]\\
    &+& \left(1 + \frac{1}{\beta}\right)\eta_{k-1}^2L^2D^2\\
    &\leq& (1+\beta)(1+\alpha)\left(1 - \frac{1}{d}\right)\|g^{k-1} - \nabla f(x^{k-1})\|^2 \\
    &+& (1+\beta)\left(1+\frac{1}{\alpha}\right)\frac{1}{\tau^2}\|f(x^{k-1} + \tau e_{i_{k-1}}) - f(x^{k-1}) - \langle \nabla f(x^{k-1}), \tau e_{i_k}\rangle\|^2\\
    &+& \left(1 + \frac{1}{\beta}\right)\eta_{k-1}^2L^2D^2\\
    &\stackrel{\ref{Nesterov}}\leq& (1+\beta)(1+\alpha)\left(1 - \frac{1}{d}\right)\|g^{k-1} - \nabla f(x^{k-1})\|^2 + (1+\beta)\left(1+\frac{1}{\alpha}\right)\frac{L^2\tau^2}{4} \\
    &+& \left(1 + \frac{1}{\beta}\right)\eta_{k-1}^2L^2D^2.
\end{eqnarray*}

If $\beta = \frac{1}{2d}$, then $(1 - \frac{1}{d})(1 + \frac{1}{2d}) \leq 1 - \frac{1}{2d}$. And with $\alpha = \frac{1}{4d}$, we get
\begin{center}
    $\EE_k\left[\|g^{k} - \nabla f(x^{k})\|^2\right] \leq \left(1 - \frac{1}{4d}\right)\|g^{k-1}-\nabla f(x^{k-1})\|^2 + 3d\eta_{k-1}^2L^2D^2 + \frac{5dL\tau^2}{4}.$
\end{center}
\EndProof{This finishes the proof.}
\begin{corollary}
For Algorithm \ref{alg:zoja} in the convex and non-convex cases the following convergences take place:
$$\mathbb{E}[r_{K+1}]= \mathcal{O}\left(\Big(f(x^0) - f(x^*)\Big)\exp\left(-\frac{K}{8d}\right) + \frac{LD^2}{K+d}\left[1 + d\sqrt{1 + \frac{K}{K+d}\frac{\tau^2}{D^2}} \right]\right).$$
$$\EE\left[\min\limits_{0 \leq k \leq K-1} \textbf{gap}(x^k)\right] = \mathcal{O}\left(\frac{f(x^0) - f(x^*)}{\sqrt{K}} + \frac{LD^2}{\sqrt{K}}\left[1 + d\sqrt{1 + K\frac{\tau^2}{D^2}}\right]\right).$$
\end{corollary}
\textbf{Proof:}\\
\EndProof{It suffices to plug parameters from Lemma \ref{lem:ap:zoja} into Theorems \ref{th:convex} and \ref{th:nonconvex}.}

\newpage 
\subsection{Distributed methods}
In this section, we provide the detailed convergence rates and proofs for specific methods (see Section \ref{sec:distr_main}) solving constrained optimization problem \eqref{eq:main} with finite-sum form \eqref{eq:finsum} in distributed/federated setting, i.e., 
$$\min\limits_{x\in\X}\Bigg\{f(x) := \frac{1}{n}\sum\limits_{i = 1}^{n}f_i(x)\Bigg\}.$$
Each device $i$ has an access only to $f_i$. Furthermore, we allow that different machine can have different data distribution, i.e., heterogeneous data setting. \\
As we have already emphasized (see Section \ref{sec:distr_main}), the bottleneck for such type of problem usually is communication cost. Therefore, we focus on methods with compressed communication.\\
In the following, we prove that some specific methods satisfy our unified Assumption \ref{as:key} and thus can be captured by our unified analysis. Then, we plug their corresponding parameters into our unified Theorems \ref{th:convex} and \ref{th:nonconvex} to obtain the detailed convergence rates for these methods.
\subsubsection{DIANA Frank-Wolfe}\label{sec:DIANA}
We first restate our main convergence lemma for \texttt{DIANA FW} method (Algorithm \ref{alg:diana}) and provide its proof for various stochastic gradients. Then we plug its corresponding parameters (i.e., specific values for $A, B, C, E, \sigma_k^2, \rho_1, \rho_2$) into our unified Theorems \ref{th:convex} and \ref{th:nonconvex} to obtain the detailed convergence rate, depending on the stochastic gradient we use. After that we provide readers with convergence rate for \eqref{eq:stoch_fw}+\eqref{eq:diana} method.
\begin{algorithm}[H]
            \caption{\texttt{DIANA Frank-Wolfe}}\label{alg:diana}
             \textbf{Input:} initial point $x^0$, $\forall i \in [n]~ h_i^{0} = \nabla f_i(x^{0})$, $h^0 = \frac{1}{n}\sum_{i=1}^nh_i^0$, step sizes $\{\eta_k\}_{k\geq0}$ , $\alpha > 0$
            \begin{algorithmic}
            \For {$k = 0, 1, 2,\dots K-1$}
            \State Compute $s^k = \arg\underset{s\in \mathcal{X}}{\min}\langle s, g^k\rangle$
            \State Update $x^{k+1} = (1 - \eta_k)x^k + \eta_ks^k$
            \State Update $h^{k+1} = h^k + \alpha\cdot\frac{1}{n}\sum_{i=1}^n \Delta_i^k$
            \For{  $i = 1, \dots, n$}
                \State Compress shifted local gradient $\Delta_i^{k+1} = \mathcal{Q}(\nabla f_i(x^{k+1}) - h_i^{k+1}) $ and send $\Delta_i^{k+1}$ to the server 
                    \State Update local shift $h_i^{k+1} = h_i^k + \alpha\cdot\mathcal{Q}(\nabla f_i(x^{k+1}) - h_i^{k+1})$
            \EndFor
            \State Aggregate received compressed gradient information $g^{k+1} = h^{k+1} + \frac{1}{n}\sum_{i=1}^n\Delta_i^{k+1}$
            \EndFor
            \end{algorithmic}
\end{algorithm}
\begin{lemma}\label{lem:ap:diana}
Under Assumption \ref{as:lip:local} with $\alpha = \frac{1}{1+\omega}$ Algorithm \ref{alg:diana} satisfy Assumption \ref{as:key} with:
$$\rho_1 = 1,~ A = \frac{\omega}{n^2}, ~ B = \frac{2\omega(\omega+1)\widetilde{L}^2}{n}, ~ C = 0,$$
$$\sigma_k^2 = \sum_{i=1}^n\|\nabla f_i(x^{k}) - h_i^k\|^2,~\rho_2 = \frac{1}{2(1+\omega)}, ~E = 2(\omega+1)n\widetilde{L}^2.$$   
\end{lemma}
\textbf{Proof:} Deriving inequalities from the proof of Theorem 7 from \citep{li2020unified}, we get
\begin{eqnarray*}
    \mathbb{E}_k\left[\|g^k-\nabla f(x^k)\|^2\right] &\leq&
    \frac{\omega}{n^2}\mathbb{E}_k \left[ \sum\limits_{i=1}^n\|\nabla f_i(x^k) - h_i^k\|^2\right]
\end{eqnarray*}
\begin{eqnarray*}
    \mathbb{E}_k \left[ \sum\limits_{i = 1}^n\|\nabla f_i(x^{k}) - h_i^{k}\|^2 \right] &\leq& \Bigl( 1 - 2\alpha + \frac{(1 - \alpha)^2}{\beta} + \alpha^2(1+\omega) \Bigl) \sum\limits_{i = 1}^n \mathbb{E}_k\Bigl[ \|\nabla f_i(x^{k-1}) - h_i^{k-1}\|^2\Bigl]\\
    &+& (1 + \beta)\sum\limits_{i = 1}^n \mathbb{E}_k\Bigl[ \|\nabla f_i(x^{k}) - \nabla f_i(x^{k-1})\|^2 \Bigl]
\end{eqnarray*}
for $\forall \beta > 0$. Choose $\beta = \frac{2\omega^2}{1+\omega}$, then
\begin{eqnarray*}
    \mathbb{E}_k \left[ \sum\limits_{i = 1}^n\|\nabla f_i(x^{k}) - h_i^{k}\|^2 \right] &\leq& \frac{\omega + \frac{1}{2}}{\omega + 1}\sum\limits_{i = 1}^n \mathbb{E}_k\Bigl[ \|\nabla f_i(x^{k-1}) - h_i^{k-1}\|^2\Bigl] + \frac{2\omega^2+\omega+1}{\omega+1}n\widetilde{L}^2\eta_{k-1}^2D^2,\\
    &\leq&\left(1 - \frac{1}{2(1+\omega)}\right)\sum\limits_{i = 1}^n \mathbb{E}_k\Bigl[ \|\nabla f_i(x^{k-1}) - h_i^{k-1}\|^2\Bigl] + 2(\omega+1)n\widetilde{L}^2\eta_{k-1}^2D^2.
\end{eqnarray*}
\begin{eqnarray*}
    \mathbb{E}_k\left[\|g^k-\nabla f(x^k)\|^2\right] &\leq& \frac{\omega}{n^2}\sum\limits_{i = 1}^n \mathbb{E}_k\Bigl[ \|\nabla f_i(x^{k-1}) - h_i^{k-1}\|^2\Bigl] + 2\frac{\omega}{n}(\omega+1)\widetilde{L}^2\eta_{k-1}^2D^2.
\end{eqnarray*}
\EndProof{This finishes the proof.}
\begin{corollary}
\textup{(Corollary \ref{crl:diana})}\textbf{.} For the algorithm \eqref{eq:stoch_fw}+\eqref{eq:diana} in the convex and non-convex cases the following convergences take place: 
\begin{eqnarray*}
    \EE[r_{K+1}] &=& \cO\left(\Big(f(x^0) - f(x^*)\Big)\exp\left(-\frac{K}{8(\omega+1)}\right) + \frac{LD^2}{K+\omega}\left[1 + \frac{\widetilde{L}}{L}\frac{(\omega+1)\sqrt{\omega}}{\sqrt{n}}\right]\right).
\end{eqnarray*}
\begin{eqnarray*}
    \EE\left[\min\limits_{0 \leq k \leq K-1} \textbf{gap}(x^k)\right] = \mathcal{O}\left(\frac{f(x^0) - f(x^*)}{\sqrt{K}} + \frac{LD^2}{\sqrt{K}}\left[1 + \frac{\widetilde{L}}{L}\frac{(\omega+1)\sqrt{\omega}}{\sqrt{n}}\right]\right).
\end{eqnarray*}
\end{corollary}
\subsubsection{MARINA Frank-Wolfe}\label{sec:MARINA}
We first restate our Lemma \ref{lem:marina} for \texttt{MARINA FW} method (Algorithm \ref{alg:marina}) and provide its proof. 
Then we plug its corresponding parameters (i.e., specific values for $A, B, C, E, \sigma_k^2, \rho_1, \rho_2$) into our unified Theorems \ref{th:convex} and \ref{th:nonconvex} to obtain the detailed convergence rate.
\begin{algorithm}[H]    
            \caption{\texttt{MARINA Frank-Wolfe}}\label{alg:marina}
             \textbf{Input:} initial $x^0$, $\forall i \in [n]~ g_i^{0} = \nabla f_i(x^{0})$, $g^0 = \nabla f(x^0)$, step sizes $\{\eta_k\}_{k\geq0}$, probability $p\in(0,1]$
            \begin{algorithmic}
            \For {$k = 0, 1, 2,\dots K-1$}
            \State Compute $s^k = \arg\underset{s\in \mathcal{X}}{\min}\langle s,g^k\rangle$
            \State Update $x^{k+1} = (1 - \eta_k)x^k + \eta_ks^k$
            \State Broadcast $g^k$ to all workers
                \For{$i =1, \dots, n$}
                \State Update $c_i^{k+1} = 
                \displaystyle\begin{cases}
                    \nabla f_i(x^{k+1}) - g_i^{k}, &\text{ with probability } p\\ \vspace{-0.9em}\\
                    \mathcal{Q}(\nabla f_i(x^{k+1}) - \nabla f_i(x^k)), &\text{ with probability }1-p
                \end{cases}
                $
                \State Send $c_i^{k+1}$ to the server
                \State Update $g_i^{k+1} = g_i^{k} + c_i^{k+1}$
                \EndFor
            \State Set $g^{k+1} = g^k + \frac{1}{n}\sum_{i=1}^nc_i^{k+1}$
            \EndFor
            \end{algorithmic}
\end{algorithm}
\begin{lemma}
\textup{(Lemma \ref{lem:marina})}\textbf{.} Under Assumption \ref{as:lip:local} Algorithm \ref{alg:marina} satisfy Assumption \ref{as:key} with:
$$\rho_1 = p, ~ A = 0, ~ B = \frac{(1-p)\omega L^2}{n},~C = 0,$$ 
$$ \sigma_k = 0, ~ \rho_2 = 1, ~ E = 0.$$
\end{lemma}
\textbf{Proof:}\\
Using the Theorem 2.1 from \citep{pmlr-v139-gorbunov21a} we can obtain:
\begin{eqnarray*}
    \EE\left[\| g^{k+1} - \nabla f(x^{k+1})\|^2\right] &\leq& 
    \frac{(1 - p)\omega L^2}{n}\EE\left[\| x^{k+1} - x^k \|^2\right] + (1 - p)\EE\left[\| g^k - \nabla f(x^k) \|^2\right]\\ 
    & = &
    \frac{(1 - p)\omega L^2}{n} \eta^{2}_{k}\EE\left[\| s^k - x^k \|^2\right] + (1 - p)\EE\left[\| g^k - \nabla f(x^k) \|^2\right]\\
    & = &
    \frac{(1 - p)\omega L^2}{n} \eta^{2}_{k}D^2 + (1 - p)\EE\left[\| g^k - \nabla f(x^k) \|^2\right].
\end{eqnarray*}
\EndProof{}
\begin{corollary}
\textup{(Corollary \ref{crl:marina})}\textbf{.} For algorithm \ref{alg:marina} in the convex and non-convex cases the following convergences take place:
$$\EE\left[r_{K+1}\right] = \mathcal{O}\left(\Big(f(x^0) - f(x^*)\Big)\exp\left(-\frac{pK}{4}\right) + \frac{LD^2}{K+\frac{1}{p}}\left[1 + \sqrt{\frac{(1-p)\omega}{np}} \right]\right).$$
$$\EE\Bigl[\min\limits_{0 \leq k \leq K-1} \textbf{gap}(x^k)\Bigr] = \mathcal{O}\left(\frac{f(x^0) - f(x^*)}{\sqrt{K}} + \frac{LD^2}{\sqrt{K}}\left[1 + \sqrt{\frac{(1-p)\omega}{pn}}\right]\right).$$
\end{corollary}
\textbf{Proof:}\\ 
\EndProof{It suffices to plug parameters from Lemma \ref{lem:ap:sarah} into Theorems \ref{th:convex} and \ref{th:nonconvex}.}

We proved results for \texttt{MARINA FW} depending on the parameters $p$ to be tuned. To find the optimal choice of $p$, we consider the particular case of RandK operator. One can note that on average we send  $pd + (1-p)k$ coordinates from a single worker per iteration. In more details, at each iteration with probability $p$ we call the full gradient in the new point $x^{k+1}$ or send the compressed difference of gradients. For $\text{RandK}$ compressor number of transmitted coordinates per iteration is equal to $k$. It is optimal to choose $p$ from the condition $pd = (1-p)k$, hence $p = \frac{k}{d+k}$. From Corollary \ref{crl:marina}, we know the estimate on the number of iterations of \texttt{MARINA FW}, then we can get an estimate on the number of the stochastic gradient calls  by multiplying this result by $pd + (1-p)k$. The final result for \texttt{MARINA FW} is presented in Table \ref{tab:comparison_distrib}. 

\subsubsection{EF21 Frank-Wolfe}\label{sec:EF21}
We first restate our Lemma \ref{lem:ef21} for \texttt{EF21 FW} method (Algorithm \ref{alg:ef21}) and provide its proof. 
Then we plug its corresponding parameters (i.e., specific values for $A, B, C, E, \sigma_k^2, \rho_1, \rho_2$) into our unified Theorems \ref{th:convex} and \ref{th:nonconvex} to obtain the detailed convergence rate.
\begin{algorithm}[H]
            \caption{\texttt{EF21 Frank-Wolfe}}\label{alg:ef21}
            \textbf{Input:} initial $x^0$, $\forall i \in [n]~ g_i^{0} = \nabla f_i(x^{0})$, $g^{0} = \frac{1}{n}\sum_{i=1}^ng_i^0$, step sizes $\{\eta_k\}_{k\geq0}$
            \begin{algorithmic}
            \For {$k = 0, 1, 2,\dots K-1$}
            \State Compute $s^k = \arg\underset{s\in \mathcal{X}}{\min}\langle s, g^{k}\rangle$
            \State Update $x^{k+1} = (1-\eta_k)x^k + \eta_ks^k$
            \State Broadcast $x^{k+1}$ to all workers
                \For{$i = 1, \dots, n$}
                \State Compress $c_i^{k+1} = \C(\nabla f_i(x^{k+1}) - g_i^{k})$ and send to the server
                \State Update $g^{k+1}_i = g_i^{k} +c_i^{k+1}$
                \EndFor
            \State Update $g^{k+1} = g^k + \frac{1}{n}\sum_{i=1}^n c_i^{k+1}$
            \EndFor
            \end{algorithmic}
\end{algorithm}
\begin{lemma}\label{lem:ap:ef21}
\textup{(Lemma \ref{lem:ef21})}\textbf{.} Under Assumption \ref{as:lip:local} Algorithm \ref{alg:ef21} satisfies Assumption \ref{as:key} with:
$$\rho_1 = 1,~ A = 1, ~ B = 0, ~ C = 0,$$
$$\sigma_k^2 = \frac{1}{n}\sum_{i=1}^n\|g_i^{k} - \nabla f_i(x^{k})\|^2,~\rho_2 = \frac{\delta+1}{2\delta^2}, ~E = 2\delta\widetilde{L}^2.$$   
\end{lemma}
\textbf{Proof:}\\
First, let us notice:
\begin{eqnarray*}
    \mathbb{E}_k\Big[\|g^{k} - \nabla f(x^{k})\|^2\Big] = \mathbb{E}_k\left[\left\|\frac{1}{n}\sum_{i=1}^n\Big(g_i^{k} - \nabla f_i(x^{k})\Big)\right\|^2\right]\leq \frac{1}{n}\sum_{i=1}^n\mathbb{E}_k\Big[\left\|g_i^{k} - \nabla f_i(x^{k})\right\|^2\Big].
\end{eqnarray*}
Similar to the Proof of Theorem 1 from \citep{NEURIPS2021_231141b3}, we can derive:
\begin{eqnarray*}
    \frac{1}{n}\sum_{i=1}^n\mathbb{E}_k\Big[\|g_i^{k} - \nabla f_i(x^{k})\|^2\Big] &=& \frac{1}{n}\sum_{i=1}^n\mathbb{E}_k\Big[\|g_i^{k-1} + \C(\nabla f_i(x^{k}) - g_i^{k-1}) - \nabla f_i(x^{k})\|^2\Big]\\
    &\leq&\left(1 - \frac{1}{\delta}\right)\frac{1}{n}\sum_{i=1}^n\|g_i^{k-1} - \nabla f_i(x^{k})\|^2 \\
    &\leq& \left(1 - \frac{1}{\delta}\right)(1 + \alpha)\frac{1}{n}\sum_{i=1}^n\|g_i^{k-1} - \nabla f_i(x^{k-1})\|^2 + \left(1-\frac{1}{\delta}\right)\Big(1 + \frac{1}{\alpha}\Big)\widetilde{L}^2\eta_{k-1}^2D^2.
\end{eqnarray*}
for any $\alpha>0$. Choose $\alpha = \frac{1}{2\delta}$, hence
$$\frac{1}{n}\sum_{i=1}^n\mathbb{E}_k\Big[\|g_i^{k} - \nabla f_i(x^{k})\|^2\Big] \leq\left(1 - \frac{\delta+1}{2\delta^2}\right)\frac{1}{n}\sum_{i=1}^n\|g_i^{k-1} - \nabla f_i(x^{k-1})\|^2 + 2\delta\widetilde{L}^2\eta_{k-1}^2D^2.$$
\EndProof{}
\begin{corollary}
\textup{(Corollary \ref{crl:ef21})}\textbf{.} For Algorithm \ref{alg:ef21} in the convex and non-convex cases the following convergences take place:
$$\mathbb{E}\Big[r_{K+1}\Big] = \mathcal{O}\Bigg(\Big(f(x^0) - f(x^*)\Big)\exp\Big(-\frac{K}{8\delta}\Big) + \frac{LD^2}{K + \delta}\Bigg[1 + \frac{\widetilde{L}}{L}\delta\Bigg]\Bigg).$$
$$\EE\Big[\min\limits_{0 \leq k \leq K-1} \textbf{gap}(x^k)\Big] = \mathcal{O}\Bigg(\frac{f(x^0) - f(x^*)}{\sqrt{K}} + \frac{LD^2}{\sqrt{K}}\Bigg[1 + \frac{\widetilde{L}}{L}\delta\Bigg]\Bigg).$$
\end{corollary}
\textbf{Proof:}\\ 
\EndProof{It suffices to plug parameters from Lemma \ref{lem:ap:ef21} into Theorems \ref{th:convex} and \ref{th:nonconvex}.}
\newpage
\subsection{Combinations of different approaches}
In this section, we provide the detailed convergence rates and proofs for specific methods (see Section \ref{sec:combination_main}) solving constrained optimization problem \eqref{eq:main} using combinations of methods, presented earlier. This approach may outperform existing methods, as they combine advantages of both algorithms.
\subsubsection{SAGA SARAH Frank-Wolfe}\label{sec:SAGA SARAH}
We first restate our Lemma \ref{lem:saga_sarah} for \texttt{SAGA SARAH FW} method (Algorithm \ref{alg:sagasarah}) and provide its proof. 
Then we plug its corresponding parameters (i.e., specific values for $A, B, C, E, \sigma_k^2, \rho_1, \rho_2$) into our unified Theorems \ref{th:convex} and \ref{th:nonconvex} to obtain the detailed convergence rate.
\begin{algorithm}[H]
            \caption{\texttt{SAGA SARAH Frank-Wolfe}}\label{alg:sagasarah}
             \textbf{Input:} initial $x^{0}$, $y_i^0 = \nabla f_i(x^0)$, $g^{0} = \nabla f(x^0)$, step sizes $\{\eta_k\}_{k\geq0}$, momentum $\lambda$, batch size $b$
            \begin{algorithmic}
            \For {$k = 0, 1, 2,\dots K-1$}
            \State Generate batch $S_k$ with size $b$
            \State Compute $s^k = \arg\underset{s\in \mathcal{X}}{\min}\langle s, g^k\rangle$
            \State Update $x^{k+1} = (1-\eta_k)x^k + \eta_ks^k$
            \State Update $y_i^{k+1} = 
            \displaystyle\begin{cases}
            \nabla f_i(x^k), &\text{for }i\in S_k,\\ \vspace{-0.9em}\\
            y_i^k, &\text{for }i\notin S_k
            \end{cases}
            $
            \State Update $g^{k+1} = \frac{1}{b}\underset{i\in S_k}{\sum}[\nabla f_i(x^{k+1}) - \nabla f_i(x^{k})] + (1 - \lambda)g^{k} + \lambda \left( \frac{1}{b}\underset{i \in S_k}{\sum} [\nabla f_i(x^{k}) - y_i^{k+1}] + \frac{1}{n}\sum\limits_{j=1}^n y_j^{k+1}\right)$
            \EndFor
            \end{algorithmic}
\end{algorithm}
\begin{lemma}\label{lem:ap:saga_sarah}
\textup{(Lemma \ref{lem:saga_sarah})}\textbf{.} Under Assumption \ref{as:lip:local} Algorithm \ref{alg:sagasarah} satisfies Assumption \ref{as:key} with: 
$$\rho_1 = \frac{b}{2n}, ~ A = \frac{b}{2n^2}, ~ B = \frac{2\widetilde{L}^2}{b},~C = 0,$$
$$\sigma_k^2 = \frac{1}{n}\sum_{j=1}^n\EE[\|\nabla f_j(x^{k}) - y_j^{k+1}\|^2],~ \rho_2 = \frac{b}{2n}, ~ E = \frac{2n\widetilde{L}^2}{b}.$$
\end{lemma}
\textbf{Proof:}\\
Using Lemma 2 from \citep{pmlr-v139-li21a} we can obtain:
\begin{eqnarray*}
    \EE_k\left[\|\nabla f(x^{k}) - g^{k}\|^2\right] &\leq& (1 - \lambda)^2\|\nabla f(x^{k-1}) - g^{k-1}\|^2 +\frac{2\lambda^2}{b}\frac{1}{n}\sum_{j=1}^n\|\nabla f_j(x^{k-1}) - y_j^k\|^2 \\
    &+& \frac{2\widetilde{L}}{b}\|x^{k} - x^{k-1}\|^2\\
    &\leq&  (1 - \lambda)^2\|\nabla f(x^k) - g^k\|^2 +\frac{2\lambda^2}{b}\frac{1}{n}\sum_{j=1}^n\|\nabla f_j(x^{k-1}) - y_j^{k}\|^2 + \frac{2\widetilde{L}}{b}\eta_{k-1}^2D^2.
\end{eqnarray*}
Additionally Lemma 3 from \citep{pmlr-v139-li21a} with $\beta_k = \frac{b}{2n}$ gives us:
\begin{eqnarray*}
    \frac{1}{n}\sum_{j=1}^n\|\nabla f_j(x^{k}) - y_j^{k+1}\|^2 &\leq& \left(1 - \frac{b}{2n}\right)\frac{1}{n}\sum_{j=1}^n\|\nabla f_j(x^{k-1}) - y_j^k\|^2 + \frac{2n\widetilde{L}^2}{b}\|x^{k} - x^{k-1}\|^2\\
    &\leq& \left(1 - \frac{b}{2n}\right)\frac{1}{n}\sum_{j=1}^n\|\nabla f_j(x^{k-1}) - y_j^k\|^2 + \frac{2n\widetilde{L}^2}{b}\eta_{k-1}^2D^2.
\end{eqnarray*}
With $\lambda = \frac{b}{2n}$ we have:
\begin{eqnarray*}
    \EE_k\left[\|\nabla f(x^{k}) - g^{k}\|^2\right] &\leq& \left(1 - \frac{b}{2n}\right)\|\nabla f(x^{k-1}) - g^k\|^2 +\frac{b}{2n^2}\frac{1}{n}\sum_{j=1}^n\|\nabla f_j(x^{k-1}) - y_j^k\|^2 + \frac{2\widetilde{L}}{b}\eta_{k-1}^2D^2.
\end{eqnarray*}
\EndProof{That finishes the proof.}
\begin{corollary}
\textup{(Corollary \ref{crl:saga_sarah})}\textbf{.} For Algorithm \ref{alg:sagasarah} in the convex and non-convex cases the following convergences take place:
$$\EE[r_{K+1}] = \mathcal{O}\left(\Big(f(x^0) - f(x^*)\Big)\exp\left(-\frac{bK}{8n}\right) + \frac{LD^2}{K+\frac{n}{b}}\left[1 + \frac{\widetilde{L}}{L}\frac{\sqrt{n}}{b}\right]\right).$$
$$\EE\left[\min\limits_{0 \leq k \leq K-1} \textbf{gap}(x^k)\right] = \mathcal{O}\left(\frac{f(x^0) - f(x^*)}{\sqrt{K}} + \frac{LD^2}{\sqrt{K}}\left[1 + \frac{\widetilde{L}}{L}\frac{\sqrt{n}}{b}\right]\right).$$
\end{corollary}
\textbf{Proof:}\\ 
\EndProof{It suffices to plug parameters from Lemma \ref{lem:ap:saga_sarah} into Theorems \ref{th:convex} and \ref{th:nonconvex}.}

The choice $b$ for \texttt{SAGA SARAH FW} is presented in the original paper \citep{Beznosikov2023SarahFM}. 

\subsubsection{Q-L-SVRG Frank-Wolfe with compression}\label{sec:lsvrgc}
We first restate our Lemma \ref{lem:lsvrg_c} for \texttt{Q-L-SVRG FW} method (Algorithm \ref{alg:lsvrgc}) and provide its proof. 
Then we plug its corresponding parameters (i.e., specific values for $A, B, C, E, \sigma_k^2, \rho_1, \rho_2$) into our unified Theorems \ref{th:convex} and \ref{th:nonconvex} to obtain the detailed convergence rate.
\begin{algorithm}[H]
            \caption{\texttt{Q-L-SVRG Frank-Wolfe}}\label{alg:lsvrgc}
            \textbf{Input:} initial $x^0$, $w^0 = x^0$, step sizes $\{\eta_k\}_{k\geq0}$, batch size $b$, probability $p\in(0,1]$
            \begin{algorithmic}
            \For {$k = 0, 1, 2,\dots K-1$} 
            \State Compute $s^k = \arg\underset{s\in \mathcal{X}}{\min}\langle s,g^k\rangle$
            \State Update $x^{k+1} = (1 - \eta_k)x^k + \eta_kg^k$
            \State Update $w^{k+1} = 
            \displaystyle\begin{cases}
                x^k, &\text{with probability } p\\ \vspace{-0.9em}\\
                w^k, &\text{with probability } 1-p
            \end{cases}
            $
            \State Broadcast $x^{k+1}$ to all workers
                \For{$i = 1, \dots, n$}
                \State Compress $c_i^{k+1} = \C(\nabla f_i(x^{k+1}) - g_i^{k})$ and send to the server
                \State Update $g^{k+1}_i = g_i^{k} +c_i^{k+1}$
                \EndFor
            \State Update
            $g^{k+1} = \dfrac{1}{n}\sum\limits_{i = 1}^{n}\Q(\nabla f_i(x^{k+1}) - \nabla f_i(w^{k+1})) + \nabla f(w^{k+1})$
            \EndFor
            \end{algorithmic}
\end{algorithm}
\begin{lemma}\label{lem:ap:lsvrgc}
\textup{(Lemma \ref{lem:lsvrg_c})}\textbf{.} Under Assumptions \ref{as:lip:local} Algorithm \ref{alg:lsvrgc} satisfy Assumption \ref{as:key} with:
$$\rho_1 = 1, ~ A = \frac{\omega\widetilde{L}^2}{n}\left(1-\frac{p}{2}\right), ~ B = \frac{\omega\widetilde{L}^2}{n}\left(1 + \frac{8(1-p)}{p}\right), ~ C = 0, $$
$$ \sigma_k^2 = \|x^{k}-w^{k}\|^2,~ \rho_2 = \frac{p}{2}, ~ E =1 + \frac{8(1-p)}{p}.$$
\end{lemma}
\textbf{Proof:}\\
\begin{eqnarray}\label{ap:eq:lsvrg+c} 
    \mathbb{E}_k\left[\| g^k - \nabla f(x^k)\|^2\right] &=& \mathbb{E}_k\left[\left\|\frac{1}{n}\sum\limits_{i = 1}^{n}\Q(\nabla f_i(x^k) - \nabla f_i(w^k)) + \nabla f(w^k) - \nabla f(x^k) \right\|^2\right] \notag \\ 
    &=& \frac{1}{n^2}\sum\limits_{i=1}^n \EE_k\left[\left\|\Q(\nabla f_i(x^k) - \nabla f_i(w^k)) + \nabla f_i(w^k) - \nabla f_i(x^k)\right\|^2\right] \notag\\
    &+& \frac{2}{n^2}\sum\limits_{i<j} \EE_k\Bigl[\Big\langle\Q(\nabla f_i(x^k) - \nabla f_i(w^k)) + \nabla f_i(w^k) - \nabla f_i(x^k), \notag \\
    && \hspace{1.9cm}\Q(\nabla f_j(x^k) - \nabla f_j(w^k)) + \nabla f_j(w^k) - \nabla f_j(x^k) \Big\rangle \Bigr] \notag \\
    &\leq& 
    \frac{\omega}{n^2}\sum\limits_{i = 1}^{n} \mathbb{E}_k \left[\|\nabla f_i(x^k) - \nabla f_i(w^k) \|^2\right] \notag
    \\[2pt]
    &\leq &  
    \frac{\omega\widetilde{L}^2}{n} \|x^k - w^k\|^2,
\end{eqnarray}
since $\Q(\nabla f_i(x^k) - \nabla f_i(w^k))$ and $\Q(\nabla f_j(x^k) - \nabla f_j(w^k))$ are independent. According to \eqref{for_comp} we derive:
\begin{eqnarray}
    \mathbb{E}_k[\|x^{k} - w^{k}\|^2] 
    &\leq& \notag
    \left(1 + \frac{2(1-p)}{\beta}\right)\eta_{k-1}^2D^2 + (1-p)(1 + 2\beta)\|x^{k-1} - w^{k-1}\|^2 .
\end{eqnarray}
Finally substituting it in \eqref{ap:eq:lsvrg+c} we get
\begin{eqnarray}
    \mathbb{E}_k[\|\nabla f(x^{k}) - g^{k}\|^2]&\leq&\frac{\omega\widetilde{L}^2}{n}\left(1 + \frac{2(1-p)}{\beta}\right)\eta_{k-1}^2D^2 + \frac{\omega\widetilde{L}^2}{n}(1 - p)(1 + 2\beta)\|x^{k-1} - w^{k-1}\|^2. \notag
\end{eqnarray}
With $\beta=\frac{p}{4}$ we have
\begin{eqnarray}
    \mathbb{E}_k[\|x^{k} - w^{k}\|^2]\leq\left(1 + \frac{8(1-p)}{p}\right)\eta_{k-1}^2D^2 + \left(1-\frac{p}{2}\right)\|x^{k-1} - w^{k-1}\|^2,
\end{eqnarray}
and
\begin{eqnarray*}
    \mathbb{E}_k[\|\nabla f(x^{k}) - g^{k}\|^2]\leq\frac{\omega\widetilde{L}^2}{n}\left(1 + \frac{8(1-p)}{p}\right)\eta_{k-1}^2D^2 + \frac{\omega\widetilde{L}^2}{n}\left(1-\frac{p}{2}\right)\|x^{k-1} - w^{k-1}\|^2.
\end{eqnarray*}
\EndProof{}
\begin{corollary}
\textup{(Corollary \ref{crl:lsvrg_c})}\textbf{.} For Algorithm \ref{alg:lsvrgc} in the convex and non-convex cases the following convergences take place:
$$\mathbb{E}[r_{K+1}] =\mathcal{O}\left(\Big(f(x^0) - f(x^*)\Big)\exp\left(-\frac{Kp}{8}\right) + \frac{LD^2}{K+\frac{1}{p}}\left[1 + \frac{\widetilde{L}}{L}\frac{\sqrt{\omega}}{p\sqrt{n}}\right]\right).$$
$$\EE\Bigl[\min\limits_{0 \leq k \leq K-1} \textbf{gap}(x^k)\Bigr] = \mathcal{O}\Bigg(\frac{f(x^0) - f(x^*)}{\sqrt{K}} + \frac{LD^2}{\sqrt{K}}\Bigg[1 + \frac{\widetilde{L}}{L}\frac{\sqrt{\omega}}{p\sqrt{n}}\Bigg]\Bigg).$$
\end{corollary}
\textbf{Proof:}\\ 
\EndProof{It suffices to plug parameters from Lemma \ref{lem:ap:lsvrgc} into Theorems \ref{th:convex} and \ref{th:nonconvex}.}

We proved results for \texttt{Q-L-SVRG FW} depending on the parameters $p$ to be tuned. Let us find the optimal choice of $p$  for \texttt{Q-L-SVRG FW}. To find the optimal choice of $p$, we consider the particular case of RandK operator. One can note that on average we send  $pd + k$ coordinates from a single worker per iteration. In more details, at each iteration with probability $p$ we call the full gradient in the new point $w^{k+1}$ and also send the compressed difference of gradient and its estimator. For $\text{RandK}$ compressor number of transmitted coordinates per iteration is equal to $k$. It is optimal to choose $p$ from the condition $pd = k$, hence $p = \frac{k}{d}$. From Corollary \ref{crl:lsvrg_c}, we know the estimate on the number of iterations of \texttt{Q-L-SVRG FW}, then we can get an estimate on the number of the stochastic gradient calls  by multiplying this result by $(pd + k)$.  
The final result for \texttt{Q-L-SVRG FW} is presented in Table \ref{tab:comparison_distrib}. 

\subsubsection{VR-MARINA Frank-Wolfe}\label{sec:vrmarina}
Developing the idea of combining different approaches, we present a modification of the \texttt{MARINA FW}, which is an adaptation from \citep{pmlr-v139-gorbunov21a}. We first present the algorithm, and then provide detailed convergence result together with its proof. 
\begin{algorithm}[H]            
            \caption{\texttt{VR-MARINA Frank-Wolfe}}\label{alg:vrmarina}
             \textbf{Input:} initial $x^0$, $g^0 = \nabla f(x^0)$, step sizes $\{\eta_k\}_{k\geq0}$, batch size $b$, probability $p\in(0,1]$
            \begin{algorithmic}
            \For {$k = 0, 1, 2,\dots, K-1$}
            \State Compute $s^k = \arg\underset{s\in \mathcal{X}}{\min}\langle s,g^k\rangle$
                \State Update $x^{k+1} = (1 - \eta_k)x^k + \eta_ks^k$
            \State Broadcast $g^k$ to all workers
                \For{$i =1, \dots, n$}
                \State Generate batch $|S_k^i| = b$
                \State Update $c_i^{k+1} = 
                \displaystyle\begin{cases}
                    \nabla f_i(x^{k+1}) - g_i^{k}, &\text{with probability } p\\ \vspace{-0.9em}
                    \Q\left(\underset{j\in S_k^i}{\sum}(\nabla f_{ij}(x^{k+1}) - \nabla f_{ij}(x^k))\right), &\text{with probability }1-p
                \end{cases}
                $
                \newline
                \newline
                \vspace{-0.3cm}
                \State Send $c_i^{k+1}$ to the server
                \State Update $g_i^{k+1} = g_i^{k} + c_i^{k+1}$
                \EndFor
            \State Set $g^{k+1} = g^k + \frac{1}{n}\sum_{i=1}^nc_i^{k+1}$
            \EndFor
            \end{algorithmic}
\end{algorithm}
In this section, we assume that the local loss on
each node has either a finite-sum form:
$$f_i(x) = \frac{1}{m}\sum\limits_{j = 1}^{m}f_{ij}(x).$$
\begin{assumption}\label{as:vrmarina}
\textup{(Average }$\mathcal{L}$\textup{-smoothness)}\textbf{.} For all $k \geq 0$ and $i \in [n]$ the minibatch stochastic gradients difference  $\widetilde{\Delta}_i^k = \frac{1}{b}\sum\limits_{S_b^k}\left(\nabla f_{ij}(x^{k+1}) - \nabla f_{ij}(x^k)\right)$ computed on the i-th machine satisfies:
\begin{equation*}
    \EE\left[\widetilde{\Delta}_i^k ~|~ x^k,x^{k+1}\right] = \Delta_i^k,
\end{equation*}
\begin{equation*}
    \EE\left[\|\widetilde{\Delta}_i^k - \Delta_i^k\|^2~|~x^k,x^{k+1}\right] \leq \frac{\mathcal{L}_i}{b}\|x^{k+1}-x^k\|^2,
\end{equation*}
\\
with some $\mathcal{L}_i > 0$ and $\Delta_i^k = \nabla f_i(x^{k+1}) - \nabla f_i(x^k)$.
\end{assumption}
This assumption is satisfied in many cases. In particular, if $S_{i,k}^{'} = \{1\ldots, m\}$, then $\mathcal{L}_i = 0$, and if $S_{i,k}^{'}$ consists of $b^{'}$ i.i.d. samples from the uniform distributions on $\{1\ldots, m\}$ and $f_{ij}$ are $L_{ij}$-smooth, then $\mathcal{L}_i \leq \max_{j\in[m]}L_{ij}$.
Under this and the previously introduced assumptions, we
derive the following result.
\begin{lemma}\label{lem:vrmarina}
Under Assumptions \ref{as:lip:local}, \ref{as:vrmarina} Algorithm \ref{alg:vrmarina} satisfies Assumption \ref{as:key} with:
$$\rho_1 = p, ~ A = 0, ~ B = \frac{(1-p)}{n}(\omega L^2 + \frac{(1 + \omega)\mathcal{L}^2}{b}),~ C = 0,$$
$$\sigma_k = 0,~ \rho_2 = 1, ~ E = 0.$$
\end{lemma}
\textbf{Proof:}\\
Using the Theorem 3.1 from \citep{pmlr-v139-gorbunov21a}, we can obtain:
\begin{eqnarray*}
    \EE_k\left[\| g^{k+1} - \nabla f(x^{k+1})\|^2\right] &\leq& 
    \frac{(1 - p)}{n}(\omega L^2 + \frac{(1 + \omega)\mathcal{L}^2}{b})\EE\left[\| x^{k+1} - x^k \|^2\right] + (1 - p)\EE\left[\| g^k - \nabla f(x^k) \|^2\right]\\ 
    & = &
    \frac{(1 - p)}{n}(\omega L^2 + \frac{(1 + \omega)\mathcal{L}^2}{b})  \eta^{2}_{k}\EE\left[\| s^k - x^k \|^2\right] + (1 - p)\EE\left[\| g^k - \nabla f(x^k) \|^2\right]\\
    & = &
    \frac{(1 - p)}{n}(\omega L^2 + \frac{(1 + \omega)\mathcal{L}^2}{b}) \eta^{2}_{k}D^2 + (1 - p)\EE\left[\| g^k - \nabla f(x^k) \|^2\right].
\end{eqnarray*}
\EndProof{}
\begin{corollary}
For Algorithm \ref{alg:vrmarina} in the convex and non-convex cases the following convergences take place:
\end{corollary}
\begin{equation*}
    \EE\left[r_{K+1}\right] = \mathcal{O}\left(\Big(f(x^0) - f(x^*)\Big)\exp\left(-\frac{pK}{4}\right) +
    \frac{D^2}{K+\frac{1}{p}}\left[L + \sqrt{\frac{(1-p)\left(\omega L^2 + \frac{1+\omega}{b}\mathcal{L}^2\right)}{np}} \right]\right).
\end{equation*}
\begin{equation*}
    \EE\Bigl[\min\limits_{0 \leq k \leq K-1} \textbf{gap}(x^k)\Bigr] = \mathcal{O}\left(\frac{f(x^0) - f(x^*)}{\sqrt{K}} + \frac{D^2}{\sqrt{K}}\left[L + \sqrt{\frac{(1-p)\left(\omega L^2 + \frac{1+\omega}{b}\mathcal{L}^2\right)}{np}}\right]\right).
\end{equation*}
\textbf{Proof:}\\
\EndProof{It suffices to plug parameters from Lemma \ref{lem:vrmarina} into Theorems \ref{th:convex} and \ref{th:nonconvex}.}

The optimal choice of parameter $p$ for \texttt{VR-MARINA FW} remains the same as in case of \texttt{MARINA FW}, hence $p = \frac{k}{d+k}$. The final result for \texttt{VR-MARINA FW} is presented in Table \ref{tab:comparison_distrib}.

\subsubsection{PP-L-SVRG Frank-Wolfe}
Combining ideas of \texttt{L-SVRG} and distributed methods, one might introduce \texttt{PP-L-SVRG}, decentralized method, that is similar to \texttt{SAGA} when size of batch equals to 1, with only exception - there is no fixed size of batch. The key idea behind Algorithm \ref{alg:lsvrg_pp} is that by choosing a random index $i_k$ at each iteration we choose the number of the device that will communicate with the server at current step. In this manner we not only update the point where the gradient is calculated with probability p, utilizing conception of the classical version of \texttt{L-SVRG}, but also reduce the number of communications to one device per iteration.

\begin{algorithm}[H]
            \caption{\texttt{PP-L-SVRG Frank-Wolfe}}\label{alg:lsvrg_pp}
            \textbf{Input:} initial $x^0=w^0$, $g^0 = \nabla f(x^0)$, step sizes $\{\eta_k\}_{k\geq0}$, batch size $b$, probability $p\in(0,1]$
            \begin{algorithmic}
            \For {$k = 0, 1, 2,\dots K-1$}
            \State Compute $s = \arg\underset{s\in \mathcal{X}}{\min}\langle s,g^k\rangle$
            \State Update $x^{k+1} = (1 - \eta_k)x^k + \eta_kg^k$
            \State Update $w^{k+1} = 
            \begin{cases}
                x^k, &\text{with probability } p\\ \vspace{-1em}\\
                w^k, &\text{with probability }1-p
            \end{cases}
            $
            \If {$w^{k+1} = x^k$}
            \For {each device} 
            \State Compute $\nabla f_i(w^{k+1})$ and send to the server 
            \EndFor
            \State Compute $\nabla f(w^{k+1}) = \frac{1}{n}\sum\limits_{i = 1}^{n} \nabla f_i(w^{k+1})$
            \EndIf
            \State Sample $i_k\in[d]$ uniformly at random
            \State Compute $g^{k+1} = \nabla f_{i_k}(x^{k+1}) - \nabla f_{i_k}(w^{k+1}) + \nabla f(w^k)$
            \EndFor
            \end{algorithmic}
\end{algorithm}
\begin{lemma}
\label{lem:ap:lsvrg_pp}
Under Assumption \ref{as:lip:local} Algorithm \ref{alg:lsvrg_pp} satisfies Assumption \ref{as:key} with
$$\rho_1 = 1, ~ A = 1+\frac{p}{2}, ~ B = \widetilde{L}^2\left(1 + \frac{2}{p}\right), ~C = 0,$$
$$ \sigma_k^2 = \sum\limits_{i=1}^n\|\nabla f_i(x^k)-w^{k+1}\|^2,~ \rho_2 = \frac{p}{2}, ~ E =\widetilde{L}^2\frac{2}{p}.$$
\end{lemma}
\textbf{Proof:}\\
We bound the difference between estimator and exact gradient:
\begin{eqnarray*}
    \EE_k\left[\|g^k - \nabla f(x^k)\|^2\right] &=& \EE_k\left[\Biggl\|\nabla f_{i_k}(x^{k}) - \nabla f_{i_k}(w^{k}) + \frac{1}{n}\sum\limits_{j = 1}^n \nabla f_j(w^{k}) - \nabla f(x^k) \Biggr\|^2\right]
    \\&\stackrel{(\ref{lem:2moment_uni_batch})}{\leq}& \frac{1}{n}\sum\limits_{j=1}^n\left\|\nabla f_j(x^k) - \nabla f_j(w^k) - \left(\frac{1}{n}\sum\limits_{i=1}^n\left[\nabla f_i(x^k)-\nabla f_i(w^k)\right]\right)\right\|^2
    \\&\leq&\frac{1}{n}\sum\limits_{j=1}^n\left\|\nabla f_j(x^k) - \nabla f_j(w^k)\right\|^2
    \\&\leq&\frac{1}{n}\left(1+\alpha\right)\sum\limits_{j=1}^n\|\nabla f_j(x^k) - \nabla f_j(x^{k-1})\|^2 + \frac{1}{n}\left(1+\frac{1}{\alpha}\right)\sum\limits_{j=1}^n\|\nabla f_j(x^{k-1}) - \nabla f_j(w^k)\|^2
    \\&\leq& \widetilde{L}^2\left(1+\alpha\right)\eta_{k-1}^2D^2 + \left(1+\frac{1}{\alpha}\right)\sigma_{k-1}^2
\end{eqnarray*}
for $\forall \alpha>0$. The second inequation holds, since $\frac{1}{n}\sum\limits_{i=1}^n$ can be described, as an expected value. And $\EE\|x - \EE x\|^2 \leq \EE\|x\|^2.$ Then we need to bound the second term:
\begin{eqnarray*}
    \EE_k[\sigma_{k}^2]&=&\EE_k\left[\frac{1}{n}\sum\limits_{j = 1}^n\|\nabla f_j(x^{k}) - \nabla f_j(w^{k+1})\|^2 \right] = \left(1 - p\right)\frac{1}{n}\sum\limits_{j = 1}^n\|\nabla f_j(x^{k}) - \nabla f_j(w^{k})\|^2 
    \\ & =& \left(1 - p\right)\frac{1}{n}\sum\limits_{j = 1}^n\|\nabla f_j(x^{k}) - \nabla f_j(x^{k-1}) + \nabla f_j(x^{k-1})-  \nabla f_j(w^{k})\|^2
    \\ &\leq& \left(1 - p\right)(1 + \beta)\frac{1}{n}\sum\limits_{j = 1}^n\|\nabla f_j(x^{k-1}) - \nabla f_j(w^{k})\|^2 + \left(1 - p\right)\left(1 + \frac{1}{\beta}\right)\widetilde{L}^2\|x^{k} - x^{k-1}\|^2.
\end{eqnarray*}
With $\beta = \frac{p}{2}$ we have:
\begin{eqnarray*}
    \EE_k[\sigma_{k}^2] \leq \left(1 - \frac{p}{2}\right)\sigma_{k-1}^2 + \frac{2}{p}\widetilde{L}^2\eta_{k-1}^2D^2.
\end{eqnarray*}
\EndProof{Taking $\alpha = \frac{2}{p}$, we obtain the needed constants.}

\begin{corollary}
\label{crl:pp}
Suppose that Assumption \ref{as:lip} holds. For Algorithm \ref{alg:lsvrg_pp} in the convex and non-convex cases the following convergences take place:
$$\mathbb{E}[r_{K+1}] =\mathcal{O}\left(\Big(f(x^0) - f(x^*)\Big)\exp\left(-\frac{Kp}{8}\right) + \frac{LD^2}{K+\frac{1}{p}}\left[1 + \frac{\widetilde{L}}{L}\frac{1}{p}\right]\right).$$
$$\EE\Bigl[\min\limits_{0 \leq k \leq K-1} \textbf{gap}(x^k)\Bigr] = \mathcal{O}\left(\frac{f(x^0) - f(x^*)}{\sqrt{K}} + \frac{LD^2}{\sqrt{K}}\left[1 + \frac{\widetilde{L}}{L}\frac{1}{p}\right]\right).$$
\end{corollary}
\textbf{Proof:} \\
\EndProof{It suffices to plug parameters from Lemma \ref{lem:ap:lsvrg_pp} into Theorems \ref{th:convex} and \ref{th:nonconvex}.}

We proved results for \texttt{PP-L-SVRG FW} depending on the parameters $p$ to be tuned. Let us find the optimal choice of $p$  for \texttt{PP-L-SVRG FW}. To find the optimal choice of $p$, one can note that on average we send  $pd + 1 - p$ coordinates from a single worker per iteration. In more details, at each iteration with probability $p$ we call the full gradient in the new point $w^{k+1}$ leveraging each worker and otherwise we choose a single device with probability $\frac{1}{d}$, which sends the update. It is optimal to choose $p$ from the condition $pd = 1-p$, hence $p = \frac{1}{d+1}$. From Corollary \ref{crl:pp}, we know the estimate on the number of iterations of \texttt{PP-L-SVRG FW}, then we can get an estimate on the number of the stochastic gradient calls  by multiplying this result by $(pd + 1 - p)$. The final result for \texttt{PP-L-SVRG FW} is presented in Table \ref{tab:comparison_distrib}. 

\newpage
\section{ADDITIONAL EXPERIMENTS}\label{sec:adexp}

In this section, we provide additional experiments comparing convergence results for several methods from Tables \ref{tab:comparison_stoch}, \ref{tab:comparison_coord}, and \ref{tab:comparison_distrib}. Since one of the main goals of our experiments is to justify the theoretical findings of the paper, in the experiments, we use the stepsizes from the corresponding theoretical results for our methods. As already described in Section \ref{sec:exp}, we assume the particular case of finite-sum constrained optimization problem \ref{eq:finsum} with $f(x)$ specified as:
$$f(x) = \frac{1}{n}\sum_{i=1}^n \log (1 + \exp( - b_i \cdot x^T a_i)),$$
where $\{a_i, b_i\}_{i=1}^n$ is $i$-th data-label pair with $a_i \in \R^d$ and $b_i \in \{-1, 1 \}$. We choose $\X$ as the $l_1$ norm ball with radius $r = 2\cdot10^3$. One can show that for given $\X$ the linear minimization oracle, i.e., $\arg\min\limits_{s\in\X} ~\< g, s>$, can be computed as: $$s^* = -\text{sign}(g_i) e_i\text{, with }i = \arg \max_j |g_j|.$$The data and labels are obtained from LibSVM datases \texttt{w1a, mushrooms} and \texttt{rcv1}.

\subsection{Point projection}
Here we introduce missing experiment of comparison \texttt{L-SVRG FW} and \texttt{SAGA SARAH FW} with effective $O(n)$ time point Euclidean projection method from \cite{duchi2008efficient} combied with SGD, as well as SVRG:

\vspace{-0.3cm}
\begin{figure}[h] 
\begin{minipage}{0.45\textwidth}
    \centering
    \includegraphics[width=1\textwidth]{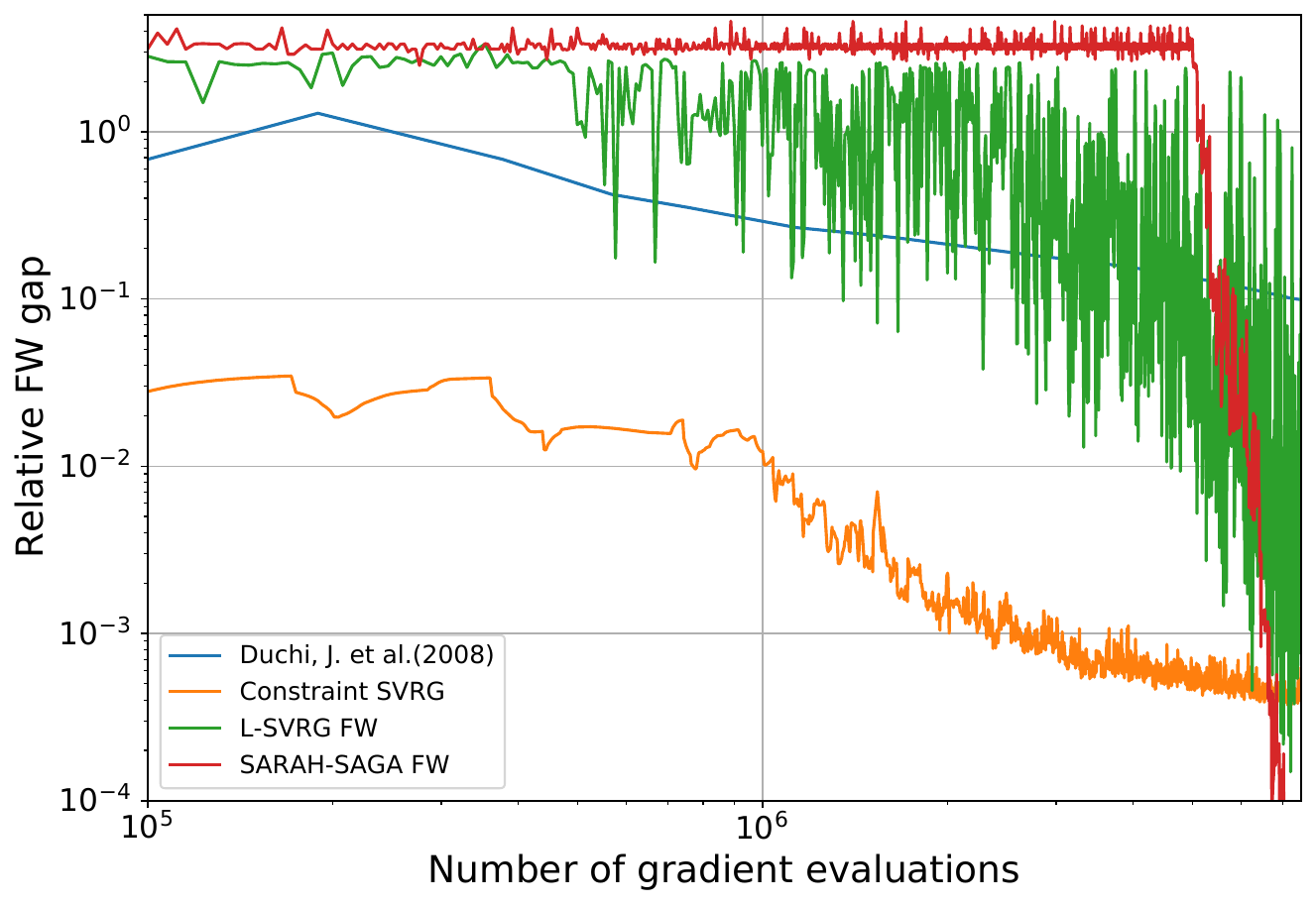}
\end{minipage}
\begin{minipage}{0.45\textwidth}
    \centering
    \includegraphics[width=1\textwidth]{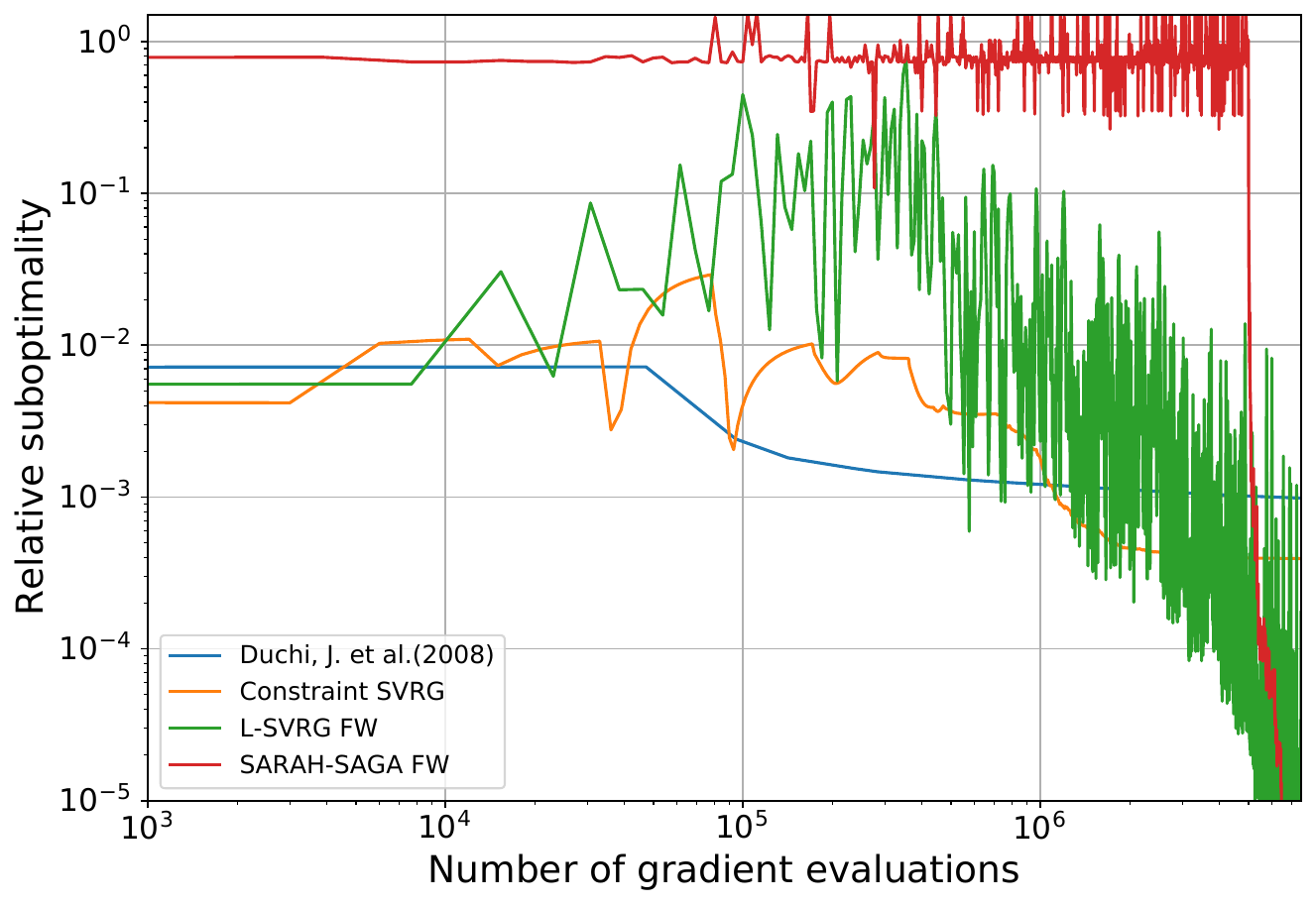}
\end{minipage}
\\
\begin{minipage}{0.01\textwidth}
\quad
\end{minipage}%
\vspace{-0.2cm}
\caption{Comparison of Frank-Wolfe-based algorithms and methods with projection for the stochastic problem. \texttt{L-SVRG FW}, \texttt{SAGA SARAH FW} as well as SGD and SVRG with projection are considered. The comparison is made on LibSVM dataset \texttt{mushrooms}.}
\end{figure}
According to conducted experience, regardless long time convergence proposed algorithms \texttt{L-SVRG}  and \texttt{SARAH-SAGA} are comparatively better than direct competitors. Despite the modest success in comparison with effective projection algorithms, the result can be considered significant.

\subsection{Additional runs}
In order to create an even more complete picture of our research, we provide additional runs below, comparing our methods with some others suitable for the tasks. In particular, we do a comparison with additional well-performing competitive approaches: SPIDER FW \citep{YU2015614} for the stochastic setting and the Block-Coordinate FW algorithm  \citep{lacoste2013block} for the coordinate setting. Comparison is made within all datasets used in the main part. 
\vspace{0.1cm}
\begin{figure}[h] 
\begin{minipage}{0.32\textwidth}\label{mush}
    \centering
    \includegraphics[width=1\textwidth]{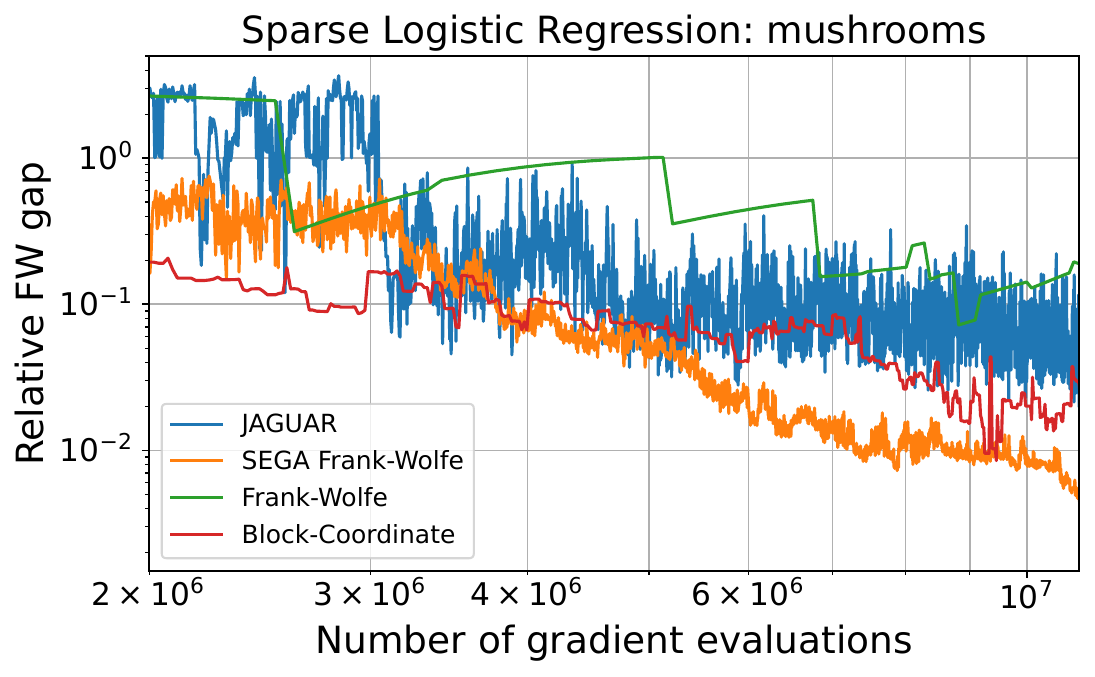}
\end{minipage}
\begin{minipage}{0.32\textwidth}
    \centering
    \includegraphics[width=1\textwidth]{mush_p.pdf}
\end{minipage}
\begin{minipage}{0.32\textwidth}
    \centering
    \includegraphics[width=1\textwidth]{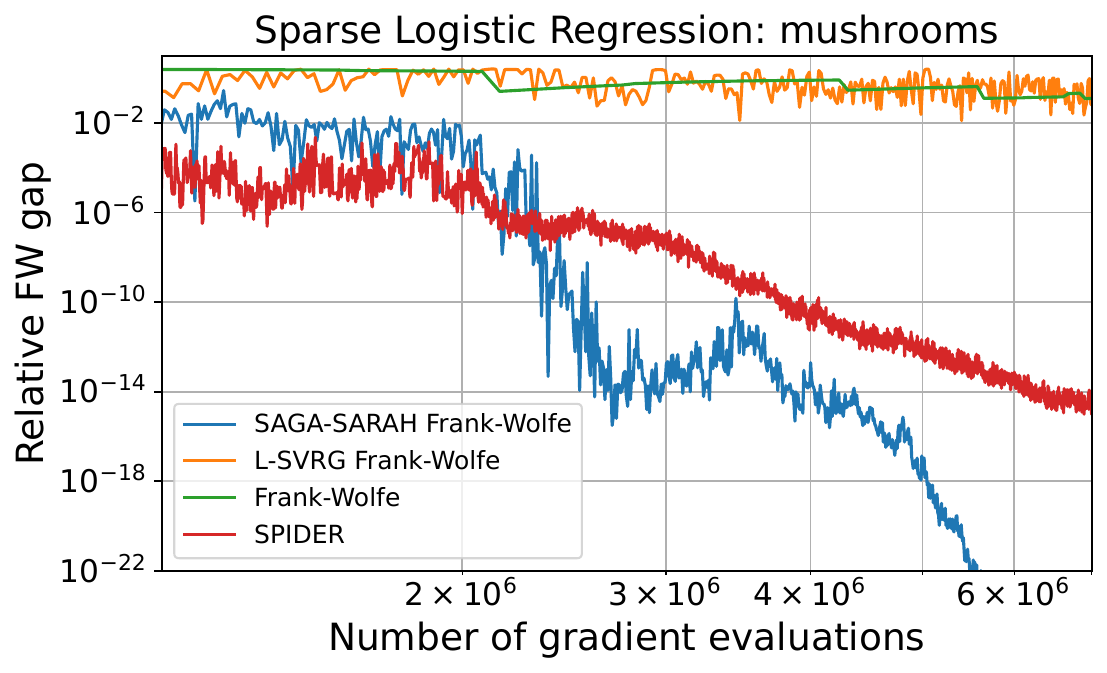}
\end{minipage}
\vspace{-0.2cm}
\caption{Comparison of methods for solving \eqref{setup} made on LibSVM dataset \texttt{mushrooms}.}
\end{figure}
\begin{figure}[h]\label{w1a}
\begin{minipage}{0.32\textwidth}
    \centering
    \includegraphics[width=1\textwidth]{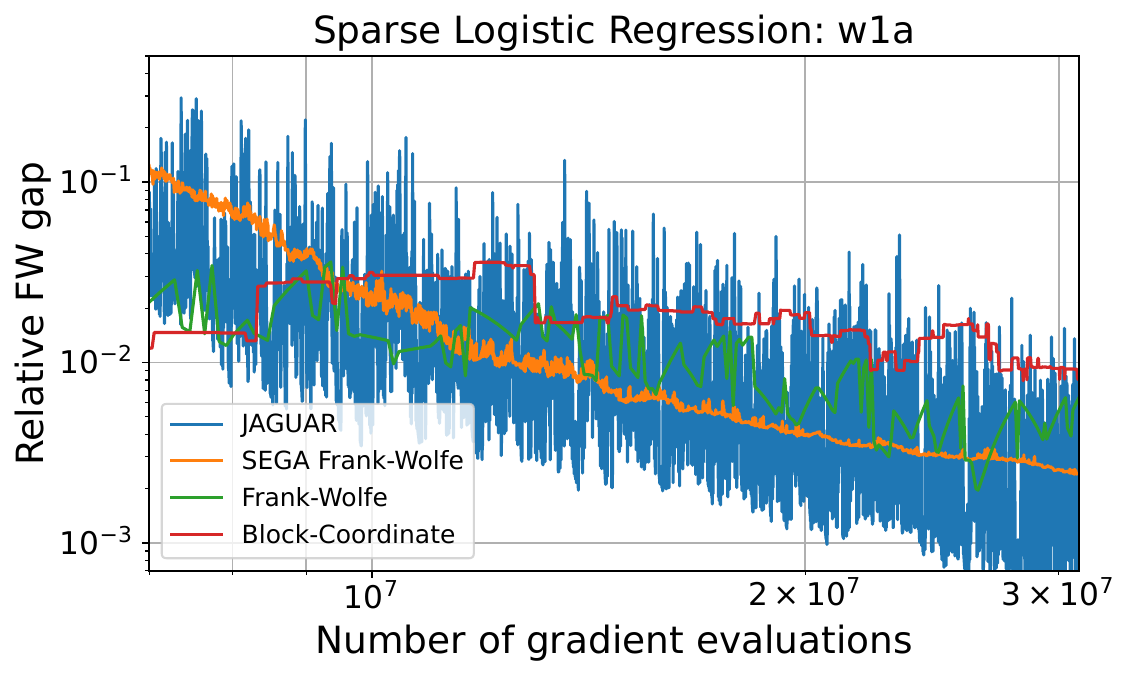}
\end{minipage}
\begin{minipage}{0.32\textwidth}
    \centering
    \includegraphics[width=1\textwidth]{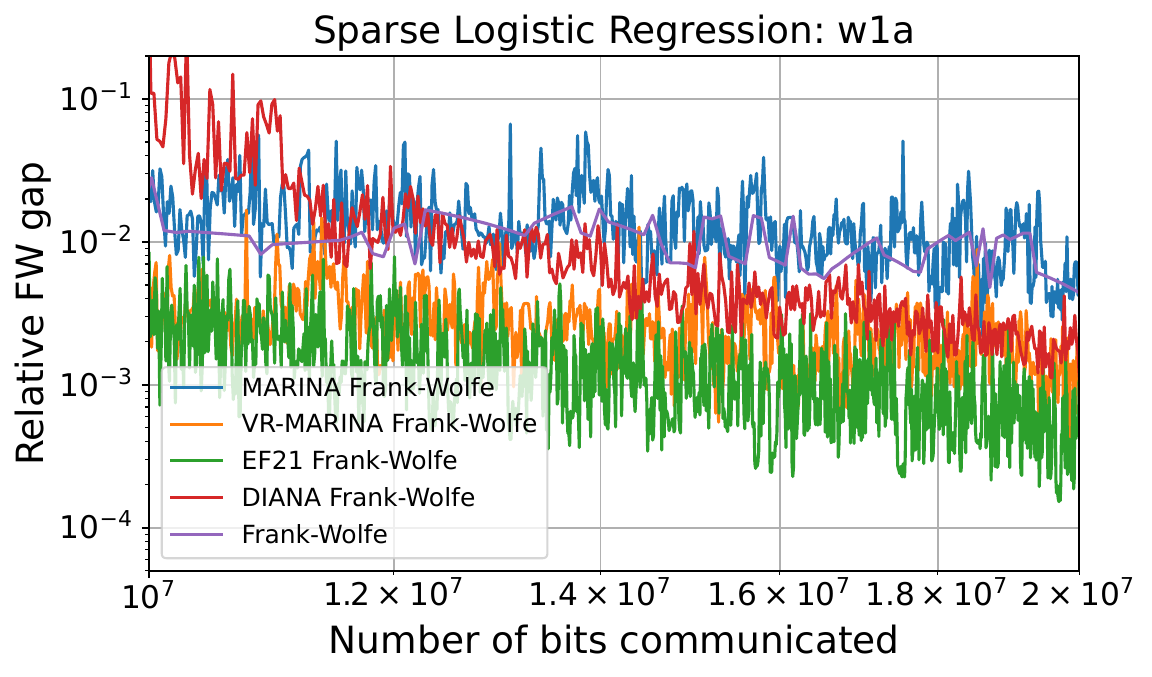}
\end{minipage}
\begin{minipage}{0.32\textwidth}
    \centering
    \includegraphics[width=1\textwidth]{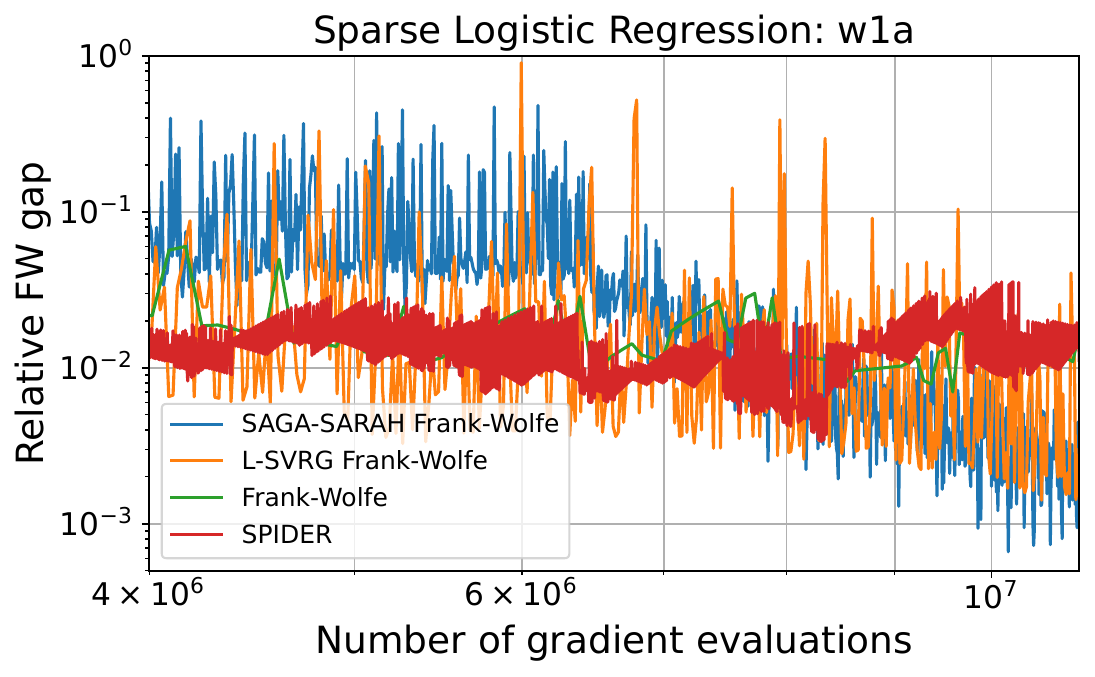}
\end{minipage}

\caption{Comparison of methods for solving \eqref{setup} made on LibSVM dataset \texttt{w1a}.}
\end{figure}
\vspace{-0.5cm}
\begin{figure}[h]\label{rcv1}
\begin{minipage}{0.32\textwidth}
    \centering
    \includegraphics[width=1\textwidth]{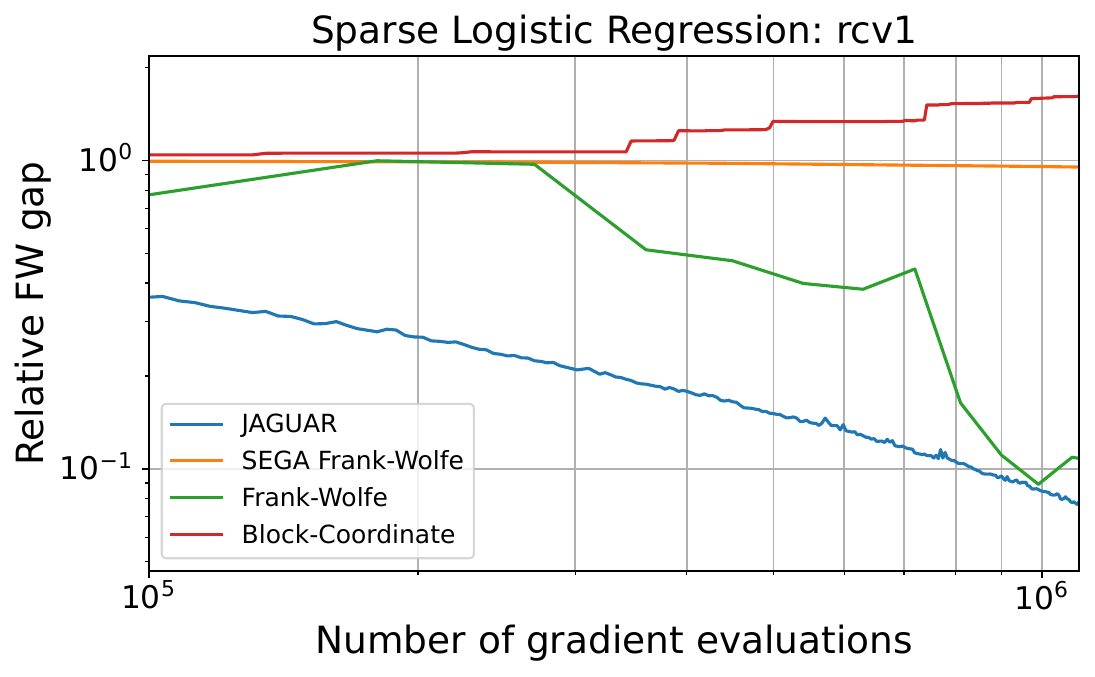}
\end{minipage}
\begin{minipage}{0.32\textwidth}
    \centering
    \includegraphics[width=1\textwidth]{rcv1_p.pdf}
\end{minipage}
\begin{minipage}{0.32\textwidth}
    \centering
    \includegraphics[width=1\textwidth]{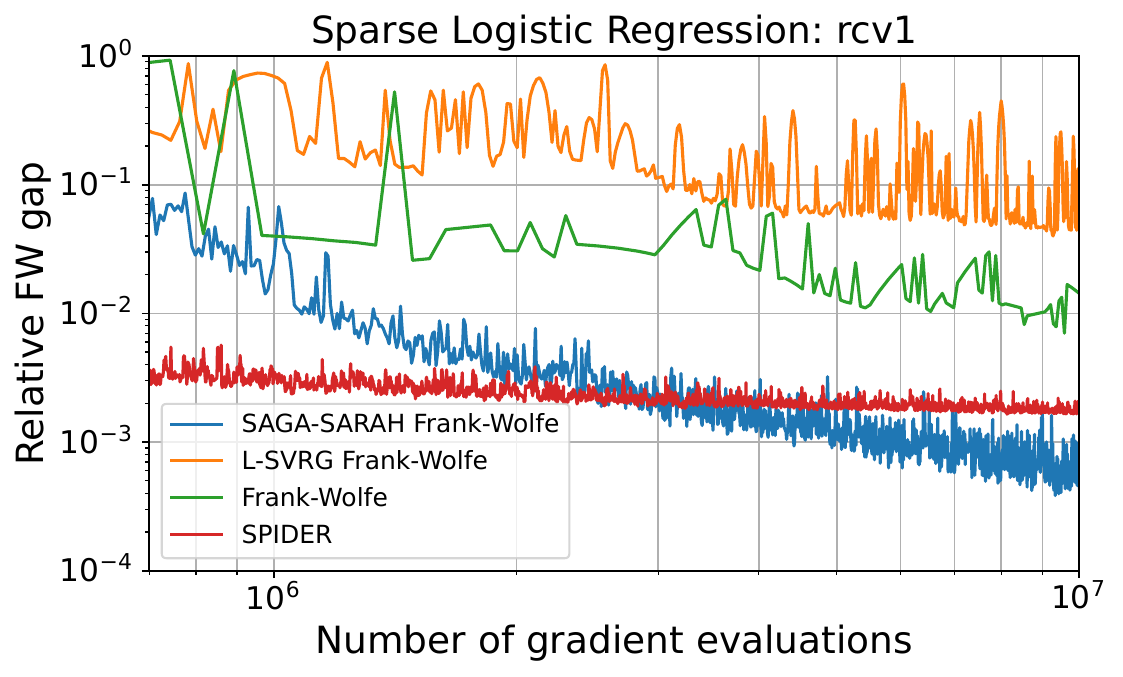}
\end{minipage}
\vspace{-0.2cm}
\caption{Comparison of methods for solving \eqref{setup} made on LibSVM dataset \texttt{rcv1}.}
\end{figure}

\end{document}